# A COMPUTATIONAL ALGORTHM FOR THE HARDY FUNCTION $Z(t)$, UTLISING SUB-SEQUENCES OF GENERALISED QUADRATIC GAUSS SUMS, WITH AN OVERALL OPERATIONAL COMPLEXITY $O((t/\varepsilon_t)^{1/3}\{log(t)\}^{2+o(1)})$

D. M. Lewis

*Department of Mathematics, University of Liverpool, M&O Building, Peach St, Liverpool L69 7ZL*

(Tel: *+44 151 794 4014*.     Email: *D.M.Lewis@liverpool.ac.uk*)

**Abstract**

This paper describes a practical methodology for computing the Hardy function $Z(t)$, using just $O\big((t/\varepsilon_t)^{1/3}\{log(t)\}^{2+o(1)}\big)$ standard computational operations, to a tolerance of $\varepsilon_t$ in the relative error. The methodology is analogous to work published relatively recently by G. H. Hiary, although the details are very different. Initially the Hardy function is formulated into sub-sequences of generalised quadratic Gaussian sums of ever increasing length $N$. Adapting a theoretical framework formulated by R. B. Paris, an algorithm for the computation of quadratic Gaussian sums using just $O\big(log(N)\big)$ operations is developed and tested. This algorithmic methodology is itself incorporated (as a sub-program) into a somewhat larger algorithm (*algorithm ZT13*) designed for the computation of $Z(t)$, which is considerably more efficient than classical $O\big(\sqrt{t}\big)$ methods, as exemplified by variants of the Riemann-Siegel formula. Sample computations using *algorithm ZT13* are presented, the results of which conform to technical theoretical predictions made regarding the limits of its performance (both in terms of its operational count and degree of precision). Speculative ideas briefly floated in the conclusions suggest that adaptations to this methodology have the potential to deliver significant further computational savings in the future.

# 1. Introduction

Numerous evaluation methodologies have been put forward for calculating the *Riemann zeta-function* $\zeta(s = \sigma + it)$, formally defined by

$$\zeta(s) = \sum_{n=1}^{\infty} n^{-s} = \prod_{all\ primes\ p} (1 - p^{-s})^{-1} \qquad (\sigma > 1), \tag{1}$$

with extension for all $s \in \mathbb{C}$ (apart from the simple pole at $s = 1$) by means of analytic continuation ([41], [12], [23]). An excellent modern review of various computationally efficient means of calculating $\zeta(s)$ over different parts of the complex plane is given in [7]. The connection between $\zeta(s)$ and the primes through the Euler product formula in (1), enabled Riemann in his paper of 1859 [35] to establish an exact formulation for the prime counting function $\pi(x)$ (not to be confused with the number $\pi$)

$$\pi(x) = \sum_{n=1}^{N} \frac{\mu(n)}{n} J(x^{1/n}), \qquad (x^{1/N} \geq 2 > x^{1/(N+1)}) \tag{2a}$$

$$J(x) = \text{Li}(x) - \sum_{Im\ \rho > 0} [\text{Li}(x^\rho) + \text{Li}(x^{1-\rho})] - \log(2) + \int_x^\infty \frac{dt}{t(t^2 - 1)\log(t)}, \tag{2b}$$

where $\text{Li}(x)$ is the logarithmic integral [15, eq. 8.240(2)] and $\mu(n)$ is the Möbius mu function e.g. [22 p32]. For an account of this work see [12]; the first proof of (2b) for the $J$ function was published in [29]. The process of Möbius inversion inherent in (2a) is discussed in a general form by [9, Section 10]. The numbers $\rho$ represent the non-trivial zeros of $\zeta(s)$. The first proofs of the Prime Number Theorem $\pi(x) \sim \text{Li}(x)$, published by [16] and [43], demonstrated that $0 < Re(\rho) < 1$, which is sufficient to show the error term in (2b) is $o(\text{Li}(x))$ as $x \to \infty$. For more modern proofs of the PNT see [9]. The famous Riemann Hypothesis (RH) asserts that $(\rho) = 1/2\ \forall \rho$, which would imply

$$\pi(x) = \text{Li}(x) + O\left(x^{1/2} \log(x)\right). \tag{3}$$

However, the current best estimate for the error

$$\pi(x) = \text{Li}(x) + O\left(x \exp\left\{-C\left[(\log(x))^3/\log(\log(x))\right]^{1/5}\right\}\right) \quad C > 0, \tag{4}$$

(see [22]), is very much larger than (3).

Given the fundamental importance of the RH to the distribution of the primes, much of the computational work devoted to $\zeta(s)$ has focused on calculating $\zeta(s)$ along the critical line $s = 1/2 + it\ (t > 0)$, where [see 12]

$$\zeta(1/2 + it) = e^{-i\theta(t)} Z(t) = e^{-i\theta(t)} 2Re\left[e^{i\theta(t)} \text{RSI}(1/2 + it)\right], \tag{5}$$

$$\theta(t) = Im\{\text{Log}\Gamma(1/4 + it/2)\} - \frac{t}{2}\log(\pi) \approx \frac{t}{2}\left\{\log\left(\frac{t}{2\pi}\right) - 1\right\} - \frac{\pi}{8} + \frac{1}{48t} + \cdots \tag{6}$$

In the above $Z(t)$ is *Hardy's Z*-function [23], formally defined by

$$Z(t) = \zeta(1/2 + it) \left[\frac{\Gamma(1/4-it/2)}{\Gamma(1/4+it/2)} \pi^{it}\right]^{-1/2}, \qquad (7)$$

while $\theta(t)$ is the Riemann-Siegel theta function. The Riemann-Siegel Integral RSI(s) (see 9) can be used in place of (1) as an alternative means of defining $\zeta(s) \,\forall\, s \neq 1$. It is relatively straightforward [23] to show that that $Z(t) = Z(-t)$ and $\theta(t) = -\theta(-t)$. The RH is equivalent to statement that $Z(Im(\rho)) = 0 \,\forall \rho$. The Riemann-Siegel Formula (RSF), discovered by C. L. Siegel from Riemann's notes in the Göttingen University library [39, 14, 7], is a computationally efficient method for calculating the Hardy function for large $t$. It states that

$$Z(t) = 2 \sum_{N=1}^{N_t} \frac{cos(\theta(t) - tlog(N))}{\sqrt{N}} - (-1)^{N_t} \left(\frac{2\pi}{t}\right)^{1/4} \sum_{r=0}^{M} (-1)^r \left(\frac{2\pi}{t}\right)^{r/2} \Psi_r(p) + R_M(t), \qquad (8)$$

where $N_t = \lfloor\sqrt{t/2\pi}\rfloor, p = frac\left(\sqrt{t/2\pi}\right)$, $\Psi_0(p) = cos\{2\pi(p^2 - p - 1/16)\}/cos\{2\pi p\}$ and $\Psi_r(p)$ are combinations of the derivatives of the function $\Psi_0(p)$. (Bounds on the cut-off remainder term $R_M(t)$ are discussed in [14], but typically $|R_2(t)| < 0.011 t^{-7/4} \,\forall\, t > 200$.) This formula forms the basis of most of the numerical calculations of $\zeta(s)$ on the critical line and has been used to show that at least the first $10^{13}$ zeros with imaginary parts up to $t = 2.5 \times 10^{12}$ do indeed satisfy the RH [13]. This type of ongoing work employs the ingenious idea of [31] of utilising a non-uniform FFT algorithm [10, 11] to compute multiple evaluations of $\zeta(s)$ along $s = 1/2 + i(t + k\tau)$ with $k = 0,1..K-1$ and $K\tau \leq \sqrt{t}$. This allows one to make the necessary $K$ evaluations of (8) in just $O\left(t^{1/2} log(K)\right)$ operations each to within an accuracy of $\pm t^{-|O(1)|}$, as opposed to the $O(Kt^{1/2})$ operations needed if the RSF was applied sequentially. An alternative methodology has recently been devised by [34] to isolate the first 103,800,788,359 zeros on the line with imaginary parts $t < 30,610,046,000$ to an accuracy of $\pm 2^{-102}$.

The RSF has been known for more than eighty years and it, along with its related variants [5], had long been thought to provide the most computationally efficient means of computing $Z(t)$ and hence $\zeta(1/2 + it)$ along the critical line. Indeed its characteristics led E. Bombieri to propose in [6] that the evaluation *of a single arbitrary* value of $\zeta(1/2 + it)$ requires a calculation which is fundamentally of $O_\varepsilon(t^{1/2-o(1)})$ in operational complexity, in order to achieve an error bounded by a fixed $\varepsilon$. However, the much more recent ground breaking work of G. A. Hiary [18, 19] proposes a number of algorithmic methods that effectively rebut Bombieri's proposition. In particular [18] re-formulates $\zeta(1/2 + it)$ into collections of either quadratic or cubic Gaussian sums of lengths $O(t^{1/6})$ and $O(t^{5/26})$ respectively. Using the principle of quadratic reciprocity, the former can be estimated recursively using $\sim O(log(t))$ operations (as discussed in section 3 and *algorithm QGS*) for arbitrary exponents, whilst the latter can be estimated in a somewhat analogous manner, provided the cubic coefficient in the exponent is not too large. This means that in order to compute $\zeta(1/2 + it)$ to an accuracy of $\pm t^{-\lambda}$ one requires only $O(t^{1/3+o_\lambda(1)})$ operations in the quadratic case and, in principle, as low as $O(t^{4/13+o_\lambda(1)})$ operations in the cubic case. Both methods incur a necessary pre-computational costing of $O(t^{1/3+o_\lambda(1)})$ and $O(t^{4/13+o_\lambda(1)})$ operations respectively, and both require similar orders of bits of storage before any actual calculations can commence. The cubic method utilises a standard Fast Fourier Transform (FFT) subroutine.

Hiary's paper [18] consists mainly of theoretical formulations and does not contain much analysis of the actual performance of the final algorithms. However, some sample computations of $\zeta(1/2 + it)$ for $t \in [10^{24-31}]$ do appear on a related website [20]. This paper seeks to build upon the algorithmic ideas of [18] (although the actual methodology is quite distinctive), in order to produce, essentially, an $O(t^{1/3+o_\lambda(1)})$ operational variant for the computation of $\zeta(1/2 + it)$, but with more emphasis on the details underlying the $O(t^{o_\lambda(1)})$ terms. To be more specific, it is the aim of this paper to demonstrate that a single evaluation of $\zeta(1/2 + it)$ for arbitrarily large $t$, requires no more than $O((t/\varepsilon_t)^{1/3}\{log(t)\}^{2+o(1)})$ *elementary computational* operations (specifically elementary *cosine* evaluations, see section 3.6) in order to yield an approximation with absolute *relative error* (not a fixed bound) less than $\varepsilon_t$. Here $\varepsilon_t > 0$ is a parameter satisfying $\varepsilon_t \xrightarrow[t \to \infty]{} 0$, although doing so more slowly than any power of $t$ (specifically $\varepsilon_t^{-1} = o(t^\mu)$ for any $\mu > 0$). For practical calculations there is a good deal of flexibility about the precise choice $\varepsilon_t$, centring around the trade-off between the greater accuracy obtained when $\varepsilon_t$ is small at the expense of a higher contribution to the operation count, and less accurate calculations (larger $\varepsilon_t$) obtained at greater speed. The final algorithm (*algorithm ZT13* of section 5) does not require any pre-computational costings or the use of FFT, but does require $O\left((log(t))^2\right)$ bits of storage for numbers of $O(log(t))$ bits.

Formulation of the methodology is founded upon the alternative definition (5) of $\zeta(1/2 + it)$ on the critical line in terms of the Riemann-Siegel Integral. The integral, first found by Bessel-Hagen from Riemann's notes [39], is defined by the formula

$$\text{RSI}(s) = \int_{0\searrow 1} \frac{e^{-i\pi z^2} z^{-s}}{e^{i\pi z} - e^{-i\pi z}} dz = \frac{1}{2i} \int_{0\searrow 1} \frac{e^{-i\pi z^2} z^{-s}}{\sin(\pi z)} dz, \qquad (9)$$

where the symbol $0\searrow 1$ denotes a path of integration along a line of slope $-1$ crossing the real axis between 0 and 1 directed from upper left to lower right. This integral was first studied in detail by Turing [42] with a view to utilising it to calculate $Z(t)$. Turing's method, as discussed in [23], was to shift the line of integration to the right to cross the real axis between $N_t$ and $N_t + 1$ to create a semi-infinite parallelogram. The integral along this new line is much easier to estimate than the original because the main bulk of the integrand is confined to a relatively small region of the integration range. The drawback is that the poles of the integrand at the integers $k = 1, 2, .., N_t$ inside the parallelogram yield the same residues as the main sum of the RSF (8). Hence Turing's method offers no significant computational improvements over the RSF and so has rarely been utilised for that purpose. (The remainder terms are also more difficult to formulate.)

However, the RSI contains far more information about $\zeta(1/2 + it)$ hidden within its structure than can be discerned using Turing's methodology. In an unpublished manuscript [27], this author carried out a thorough analysis of (9) which revealed a new and distinctive representation for $Z(t)$. This was achieved by integrating (9) directly along the line $0 \searrow 1$ (specifically along the line crossing the real axis at $z = 1/2$) with $s = 1/2 + it$ for $t > 0$. The main difficulty with this direct approach is that the $z^{-s}$ term grows in competition with the $e^{-i\pi z^2}$ term when $Re(z) < 0$, which forces one to consider (9) as integral over an essentially infinite range (the bulk of the integral covers a range $\propto t^{1/3}$ as $t \to \infty$) to calculate it effectively. But with some lengthy analysis [27 - Appendix A] it is possible to circumvent this problem. The result is a representation of the Hardy function in terms of a sum of the confluent hypergeometric or Kummer's function [32, Chap 13], given by

$$Z(t) = \frac{2e^{i\theta(t) - 3\pi t/4 - i\pi/8} \pi^{5/4 + it}}{\Gamma(1/4 + it/2)(1 + e^{-2\pi t})} \sum_{\alpha \in 2\mathbb{N}+1} \alpha \Phi\left(\frac{3}{4} - \frac{it}{2}, \frac{3}{2}, \frac{\pi \alpha^2 i}{4}\right)$$

$$= \frac{2e^{-\pi t/2 - i\pi/8}}{(1 + e^{-2\pi t})\sqrt{2\pi}} \left(\frac{t\pi^5}{2}\right)^{1/4} e^{-[1/32t^2 + \varpi]} \sum_{\alpha \in 2\mathbb{N}+1} \alpha \Phi\left(\frac{3}{4} - \frac{it}{2}, \frac{3}{2}, \frac{\pi \alpha^2 i}{4}\right). \quad (10)$$

Essentially (10) is equation (A40) of [27] for $RSI(1/2 + it)$, combined with (5) to obtain $Z(t)$. Here $\Phi(a, b, z) \equiv {}_1F_1(a, b, z)$ is Kummer's function and $\varpi = 5/256t^4 + O(t^{-6})$ (arising from *Sterling's* series for the Gamma function). Formally the sum is taken over all positive odd integers $\alpha \geq 1$. However, because it turns out that its terms decay only as $1/\sqrt{\alpha}$, the infinite sum does not converge and in practice it must be terminated at some odd integer $N_\alpha \geq INT_O(t/\pi) + 2$. (The notation $N/INT_{O,E}$ will be used to denote the round and floor functions, to the nearest *O*dd and *E*ven integers respectively.) Beyond this cut-off the remainder of the series can be continued analytically and estimated to very high accuracy using Euler-Maclaurin summation (in just the same way as (1) can be continued analytically for $\sigma \leq 1$). These observations concerning the truncation of (10) can be established from estimates [27] of the respective Kummer's functions, generated using some standard methods of contour integration. An adapted version of this estimation methodology will be discussed in some detail in Section 2, because it is fundamental to the proposition that the calculation $\zeta(1/2 + it)$ is an $O\left((t/\varepsilon_t)^{1/3}\{log(t)\}^{2+o(1)}\right)$ computational operation for large $t$. But for the moment applying the methodology on a term by term basis [27], in which each Kummer function in (10) is estimated separately, one obtains the following somewhat simpler expression.

$$Z(t) = \mathcal{H}(t) 2\sqrt{2} \left\{ \sum_{\substack{\alpha > a \\ odd}}^{N_\alpha} \frac{\cos\left(\frac{t}{2}\left\{\log(pc) + \frac{1}{pc}\right\} + \frac{t}{2} + \frac{\pi}{8}\right)}{(\alpha^2 - a^2)^{\frac{1}{4}}} [1 + O(\Lambda_\alpha^{-1})] + R_{EM} \right\} + T_{a+\varepsilon}, \quad (11)$$

where $\quad \mathcal{H}(t) = \dfrac{e^{-[1/32t^2 + O(t^{-4})]} e^{2/(12t^2+3) + O(t^{-4})}}{(1 + e^{-2\pi t})(1 + 1/4t^2)^{1/4} e^{[1/24t^2 + O(t^{-4})]}} \approx 1 + \dfrac{1}{32t^2}. \quad (11a)$

Here $a = \sqrt{8t/\pi}$, the relative error $\Lambda_\alpha = MIN\{t^{1/8}(\alpha - a)^{3/4}, \sqrt{t}\}$ and $R_{EM}$ represents the Euler-Maclaurin remainder terms [see specifically eqs. A73-74 of 27]. The "portcullis" $pc$ variable (so designated because it acts as a gateway from the summands in (11) to the summands in (8), see below) is defined for all real $\alpha \geq a$ by

$$pc(\alpha) = 2(\alpha/a)^2\left[1 + (1 - (a/\alpha)^2)^{1/2}\right] - 1 \in [1, \infty),$$

$$\Rightarrow \pi\alpha^2/2 = t(pc+1)^2/pc \Rightarrow \alpha/a = (pc+1)/2\sqrt{pc}. \tag{12}$$

Actually $pc(\alpha)$ fixes the position of a saddle point of a certain integral (see eq. 21) representation of the respective Kummer functions. Note that for reasons discussed in section 2.1, the approximate sum (11) commences at $\alpha = INT_O(a) + 2$, the nearest odd integer above $a$, not at $\alpha = 1$ as in the exact sum (10). A slight modification to (11) is required whenever $a$ itself happens to be very close to an odd integer, specifically when $NINT_O(a) - a = \varepsilon$, where $\varepsilon \in [\pm t^{-1/6}]$. If $\varepsilon > 0$ the main sum begins at $\alpha = INT_O(a) + 4$ and a transition term $T_{a+\varepsilon}$, representing the contribution at $\alpha = INT_O(a) + 2$, must be included. For $\varepsilon = gt^{-1/6} > 0$ and small $g < 1/4$ this takes the form q.v. [27 Appendix A, section A3.6, eq. A65]

$$T_{a+\varepsilon} = \mathcal{H}(t) \frac{2^{3/4}\Gamma(1/3)\exp(-(32\pi^3)^{1/4}g^{3/2}/3)}{3^{2/3}\pi^{1/4}t^{1/12}} \left[\cos\left(t + \sqrt{\frac{\pi}{2}}t^{1/3}g + \frac{\pi}{24}\right) + O(g)\right]. \tag{13}$$

For larger $g \in [1/4, 1]$ and for the case when $\varepsilon < 0$ it is straightforward to calculate $T_{a+\varepsilon}$ numerically (in the latter scenario the main sum recommences at $\alpha = INT_O(a) + 2$). However, it is difficult to formulate any relatively simple approximations for $T_{a+\varepsilon}$ that smoothly transform themselves into the first summation term of (11) in these instances. Dropping all the relatively small $O(\Lambda_\alpha^{-1})$ contributions in (11) and concentrating on the sum of the main terms alone, leads to the following theorem.

*Theorem 1*

Let $t > t_0 = 30$. Define $N_t = INT[(t/2\pi)^{1/2}]$, let $\theta(t)$ be the Riemann-Siegel theta function (6) and $\theta_C(t) = \frac{t}{2}\left\{\log\left(\frac{t}{2\pi}\right) - 1\right\} - \frac{\pi}{8} = \theta(t) - \Delta(t)$, where $\Delta(t)$ represents those terms $< O(1)$ in (6). Then there exists an upper bound $E(t) > 0$, satisfying the $\lim_{t\to\infty} E(t) = 0$, such that

$$\left|2\sum_{N=1}^{N_t} \frac{\cos\{\theta_C(t) - t\log(N)\}}{\sqrt{N}} - 2\sqrt{2}\sum_{\substack{\alpha > a \\ odd}}^{N_\alpha} \frac{\cos\left(\frac{t}{2}\left\{\log(pc) + \frac{1}{pc}\right\} + \frac{t}{2} + \frac{\pi}{8}\right)}{(\alpha^2 - a^2)^{\frac{1}{4}}}\right| < E(t) = \frac{6.15}{t^{1/12}}. \tag{14}$$

The proof of this theorem is the main focus of [27], in which $N_\alpha = INT_O(2t/\pi - a)$ for convenience. However, there is nothing special about this value and the second sum can be terminated at any $N_\alpha \geq INT_O(t/\pi) + 2$ and then continued analytically using Euler-Maclaurin summation [e.g. 59 eq. 2.10.1]. (The proof for any general $N_\alpha$ proceeds virtually identically, provided the integration contour of [27 eq. 8] is defined to be the circle $C(R)$ centred at $z = (N_\alpha + a)/2$ with maximum radius $R_{max} \approx (N_\alpha - a)/2$.) Notice also the appearance of $\theta_C(t)$ rather than $\theta(t)$ which appears in the definition of (8). A precise definition of the difference $\Delta(t) = \theta(t) - \theta_C(t) \approx 1/48t$ is given by [23 eq. 5.2]. It is also important to point out that the transition term (13) is not included in this estimate of the discrepancy between the main sums of (8) & (11). In fact it is this omission which is the main source of error in (14), because $|T_{a+\varepsilon}| = O(t^{-1/12})$. It the difficulty of establishing the precise transition between (13) and the main terms of (11) that is reflected in the conservative value assigned to $E(t)$. But the most interesting result is not (14) *per se*, but rather the connection established in the course of its proof between the respective summands $\alpha$ and $N$ by means of the portcullis variable (12). The terms in (14) corresponding to the main sum of the RS formula arise from the saddle points of certain integrals defined by [27 eq. 26] situated at $N = a/4\sqrt{pc_{sad}} \Rightarrow pc_{sad} = t/2\pi N^2$. Eliminating $pc$ from (12) establishes a direct connection (for real $\alpha \geq a$) between the summands $N$ in (8) and $\alpha$ in (11) given by

$$\alpha = 2N\left(\frac{t}{2\pi N^2}+1\right) \Rightarrow N = \frac{\alpha}{4}\left[1-\sqrt{1-\left(\frac{a}{\alpha}\right)^2}\right]. \tag{15}$$

This connection leads to some fascinating results. Utilising the ideas used to prove *Theorem 1*, [27] was able to demonstrate (essentially the idea is simply to reduce the radius of $C(R)$ systematically, to exclude the terms $\alpha = INT_O(a) + 2$, $INT_O(a) + 4, \ldots$ one at a time, and reapply the methodology to the smaller circle) that moving between from $\alpha = a$ to $\alpha = INT_O(a) + 2$ to corresponds to a change from $N_t$ to a point $N = N_t - \Delta N$ where $\Delta N = O[(t/8\pi)^{1/4}]$. Hence the first term alone of (11), the calculation of which an $O(1)$ computational operation, corresponds to the result of summing the last $O[(t/8\pi)^{1/4}]$ terms of the main sum RS formula, an $O(t^{1/4})$ operation. More specifically let $\alpha_1 = INT_O(a) + 2$, $\alpha_2 = \alpha_1 + 2$, and denote the corresponding nearest integer values satisfying (15) by $N_1$ and $N_2$ respectively. Then

$$\frac{2\sqrt{2}\cos\left(\frac{t}{2}\left\{\log(pc(\alpha_1))+\frac{1}{pc(\alpha_1)}\right\}+\frac{t}{2}+\frac{\pi}{8}\right)}{(\alpha_1^2-a^2)^{1/4}} \approx \sum_{N=(N_1+N_2)/2}^{N_t} \frac{2\cos\{\theta_C(t)-t\log(N)\}}{\sqrt{N}}. \tag{16}$$

So for example if $t = 10^{24}$, then $\alpha_1 = 1595769121607$, $(N_1 + N_2)/2 = 398941625041$ and $N_t = 398942280401$. The sum of the 655,360 terms on the right of (16) computes to $4.727 \times 10^{-4}$ compared to the predicted value of $4.744 \times 10^{-4}$ on the left. The *relative* error of $0.0034 \equiv 2.8\Lambda_\alpha^{-1}$, lies in the range specified by (11).

This kind of dramatic computation saving cannot be expected to continue indefinitely because using the entire sum (11) to estimate $Z(t)$ is manifestly an $O(t)$ computational operation compared to the $O(\sqrt{t})$ operations needed for the RSF. As $\alpha$ increases each individual term of (11) corresponds to ever smaller portions of the main sum of the RSF, until one reaches a boundary point situated at $= 2 \equiv \alpha = 3\sqrt{t/\pi} \equiv N = \frac{1}{2}\sqrt{t/\pi}$, when the terms are in a rough one-to-one correspondence. Moving on into the $pc > 2$ regime, one finds the roles of terms making up (8) and (11) swap over, so that now each individual term of the RSF now corresponds to ever larger portions of the sum over $\alpha$. But the swap over boundary at $pc = 2$ is still sufficiently large to realise a significant reduction in the operational count necessary to compute $Z(t)$. This can be achieved by utilising a hybrid of the two main sums of (8) and (11). The crux of this idea is to substitute the last $N \in \left[\frac{1}{2}\sqrt{t/\pi}, N_t\right]$ terms of the RSF with the first $\alpha \in \left[INT_O(a) + 2, INT_O\left(3\sqrt{t/\pi}\right)\right]$ terms of (11), to give

$$\sum_{N=\lceil\sqrt{t/4\pi}\rceil}^{N_t} \frac{2\cos\left\{\frac{t}{2}\log\left(\frac{t}{2\pi}\right) - \frac{t}{2} - \frac{\pi}{8} - t\log(N)\right\}}{\sqrt{N}}$$

$$= \sum_{\alpha=INT_O(a)+2}^{INT_O(3\sqrt{t/\pi})} \frac{2\sqrt{2}\cos((t/2)\{\log(pc) + 1/pc\} + t/2 + \pi/8)}{(\alpha^2 - a^2)^{\frac{1}{4}}} + T_{a+\varepsilon} + error.$$

(17)

Setting the cut-off at $pc = 2 \equiv \alpha = 3\sqrt{t/\pi} \equiv N = \frac{1}{2}\sqrt{t/\pi}$ as in (17), reduces the total number of terms needed compute $Z(t)$ from the $N_t$ using (8) alone, to a minimum of $2(\sqrt{2} - 1)N_t$ using (8) and (17) in combination, a reduction of $\approx 17.16\%$. In computer time the saving is closer to 15% because the terms on the right of (17) are slightly more complicated to evaluate (but still elementary) compared to those of (8). Although theoretically the error in (17) is no better than $E(t)$ found in (14), in practice millions of computations [Lew] show that it bounded above by a term $< 1.01(64\pi/t)^{1/4}$. The size of this bound corresponds to the maximum magnitude of the terms of both (8) and (11) at the cut-off point $pc = 2$.

In light of these computational savings, an obvious question to ask is whether one can do any better? There seems to be scope for such a scenario, simply because (11) is a sum of $O(t)$ elementary terms, and as such contains huge tranches of computational redundancy (in the regime when $pc > 2$), as well as the computational efficiencies highlighted in (16) and (17). Understanding the nature of the computational redundancy hidden within the mathematical structure of (11), provides the first step towards the goal reducing the computational complexity of $Z(t)$ much further.

## 2. Formulating Hardy's Z-Function into Generalised Quadratic Gauss Sums

To make progress one must first appreciate how the relatively simple approximation (11) arises from the exact representation (10) for $Z(t)$. All the key details are discussed in [27, Appendix A] so what follows is only a brief summary. Kummer's function $\Phi(3/4 - it/2, 3/2; \pi\alpha^2 i/4)$ can be expressed in terms of Euler's integral formulation [28, eq. 4.2(1)]

$$\Phi\left(\frac{3}{4} - \frac{it}{2}, \frac{3}{2}; \frac{\pi\alpha^2 i}{4}\right) = \frac{\Gamma(3/2)}{\Gamma(3/4 + it/2)\Gamma(3/4 - it/2)} \int_0^1 e^{i\pi\alpha^2 x/4} \frac{\exp[i(t/2)\log\{(1-x)/x\}]}{[x(1-x)]^{1/4}} dx,$$

$$= \frac{\pi^{-1/2}e^{\pi t/2}e^{2/(12t^2+3)+O(t^{-4})}}{4(t^2/4 + 1/16)^{1/4}e^{1/24t^2+O(t^{-4})}} \int_0^1 e^{i\pi\alpha^2 x/4} \frac{\exp[i(t/2)\log\{(1-x)/x\}]}{[x(1-x)]^{1/4}} dx, \quad (18)$$

after applying Stirling's asymptotic expansion to both complex gamma functions. Substituting into (10) gives

$$Z(t) = \mathcal{H}(t)e^{-i\pi/8}\left(\frac{\pi}{2^5 t}\right)^{\frac{1}{4}} \sum_{\alpha \in 2\mathbb{N}+1} \alpha \int_0^1 e^{i\pi\alpha^2 x/4} \frac{\exp[i(t/2)\log\{(1-x)/x\}]}{[x(1-x)]^{1/4}} dx. \quad (19)$$

The integral in (19) can be split into two and written as $B(\alpha, t) + e^{i\pi\alpha^2/4}\overline{B(\alpha, t)}$, where

$$B(\alpha, t) = \int_0^{1/2} e^{i\pi\alpha^2 x/4} \frac{\exp[i(t/2)log\{(1-x)/x\}]}{[x(1-x)]^{1/4}} dx = \int_1^\infty \frac{\exp[i\pi\alpha^2/4(w+1) + i(t/2)log(w)]}{w^{1/4}(w+1)^{3/2}} dw. \quad (20)$$

The second form of $B(\alpha, t)$ is obtained following the substitution $x = (w+1)^{-1}$ and provides the easiest means of estimating its asymptotic value.

2.1 *Asymptotic Analysis of the integral* $B(\alpha, t)$

This section summarises the main points of [27 eq. A3.3], concerning the asymptotic behaviour of $B(\alpha, t)$ needed for the improvements of sections 2.2-4. The basic idea is to treat (20) as a contour integral, let $w \in \mathbb{C}$ and devise an appropriate path that is equivalent to the real line $w \geq 1$ and then invoke Cauchy's Theorem. This approach is valid provided the entire path lies in the $\text{Re}(w) > 0$ plane, where the integrand in (21) is analytic. Two suitable paths are shown in Figs. 1a & b, one for $\alpha < a$ and one for $\alpha > a$ respectively. Writing $w = Re^{i\phi}$ with $R \geq 1$ and $\phi \in [-\pi/2, \pi/2]$. The phase of the integrand (20) is given by

$$Re\left[i\left\{\frac{\pi\alpha^2}{4(w+1)} + \left(\frac{t}{2}\right)log(w)\right\}\right] = \left[\frac{\pi\alpha^2 R\sin(\phi)}{4(R^2+1+2R\cos(\phi))} - \frac{t\phi}{2}\right]. \quad (21)$$

So if $\alpha < a$, then (21) is negative $\forall R > 1$ and $\phi > 0$. Hence integrating along any path in this region starting from $|w| = 1$ and finishing at an arbitrary $w = Re^{i\phi>0}$ gives a negligibly small $O(e^{-O(t)})$ contribution to (20). On the circle $|w| = R$ the integrand is $O(e^{-O(t)}/R^{7/4})$, which means the integral vanishes as the radius $R \to \infty$. The consequence of these two observations is that for $\alpha < a$ the only significant contribution to (20) along the path in Fig. 1a, will be the result of the integration on the small portion of the unit circle beginning at $w = 1$. But in [27 eq. A45-46] it was established that this contribution always takes the form

$$B_{integrating\ on\ unit\ circle}(\alpha, t) = ie^{i\pi\alpha^2/8} \times real\ number, \quad (22)$$

irrespective of the direction of rotation from $w = 1$. The value of the *real number* is immaterial, because when (22) is substituted into $B(\alpha, t) + e^{i\pi\alpha^2/4}\overline{B(\alpha, t)}$ the result is automatically zero. So all the terms in (19) associated with odd integers $\alpha < a$ can only be of $O(e^{-O(t)})$ in size, and the sum of these is negligible in comparison with the $\alpha > a$ terms. This is the reason for commencing the summation at $\alpha \geq a$ in (11) and (14). When $\alpha > a$, the integration contour is as shown in Fig. 1b. One moves around the unit circle to the point $w = -i$, which again contributes nothing to (19) because of (22), and then along two mutually perpendicular lines, the second of which crosses the real axis at the point $w = pc(\alpha)$ given by (12). This defines the saddle point of the integrand's (20) phase and it is at this saddle point when the overwhelming bulk of the integral is concentrated. Using standard saddle point methods [2] one finds that

$$B(\alpha > a, t) = \frac{e^{i\pi/4} 2\sqrt{\pi}(pc)^{3/4} \exp[i\{\pi\alpha^2/4(pc+1) + tlog(pc)/2\}]}{(pc+1)\sqrt{t(pc-1)}} [1 + O(\Lambda_\alpha^{-1})] \Rightarrow$$

$$\int_0^1 \frac{e^{i\pi\alpha^2 x/4 + i(t/2)log\left\{\frac{(1-x)}{x}\right\}}}{[x(1-x)]^{1/4}} dx = \frac{e^{\frac{i\pi}{8}}(pc)^{3/4} 4\sqrt{\pi} \cos(\pi\alpha^2/4(pc+1) + (t/2)log(pc) + \pi/8)}{(pc+1)\sqrt{t(pc-1)}}[1 + O(\Lambda_\alpha^{-1})].$$

(23a, b)

Substituting these results into (19) along with the simplification (an exact consequence of 12)

$$\frac{2^{5/4}\pi^{1/4}(pc)^{1/4}}{t^{1/4}\sqrt{pc-1}} = \frac{2\sqrt{2}}{(\alpha^2 - a^2)^{1/4}},$$ (24)

gives (11). The case when $\alpha = a + \varepsilon$, where $\varepsilon \in [\pm t^{-1/6}]$, is discussed extensively in [27 eq. A3.6] but is not relevant for this investigation.

2.2 *Combining together collections of integrals that form the sum (19) for $Z(t)$.*

The fact that the sum (19) effectively commences at $\alpha = INT_O(a) + 2$ rather than $\alpha = 1$, provides a means of stripping away the computational redundancy inherent in the formulation (11). For example, suppose we have two neighbouring terms at $\alpha_1 = \alpha_E - 1$ and $\alpha_2 = \alpha_E + 1 > a$, where $\alpha_E$ is the corresponding *even* integer between them. (In what follows the terminology *pivot* integer will be used for $\alpha_E$.) Then it straightforward to show that

$$\sum_{k=1}^2 \alpha_k \int_0^1 e^{i\pi\alpha_k^2 x/4} \frac{\exp[i(t/2)log\{(1-x)/x\}]}{[x(1-x)]^{1/4}} dx$$

$$= \sum_{k=1}^2 (\alpha_E + (-1)^k) \int_0^1 e^{i\pi(\alpha_E + (-1)^k)^2 x/4} \frac{\exp[i(t/2)log\{(1-x)/x\}]}{[x(1-x)]^{1/4}} dx$$

$$= 2\int_0^1 e^{i\pi\alpha_E^2 x/4} \frac{\exp[i(t/2)log\{(1-x)/x\}]}{[x(1-x)]^{1/4}} \left\{\alpha_E \cos\left(\frac{\pi\alpha_E x}{2}\right) + i\sin\left(\frac{\pi\alpha_E x}{2}\right)\right\} e^{i\pi x/4} dx. \quad (26)$$

Now because $\alpha_E > a$, so that $\alpha_E \geq O(\sqrt{t})$, integral (26) for the two combined terms has much the same *asymptotic properties* as that of (19) for a single term. The main phase is of $O(\alpha_E^2) \geq O(t)$, whilst the terms in { } involve phases of $O(\alpha_E)$ which change relatively slowly in comparison. Hence estimates for integral (26) are computable using very much the same methodology as described in section 2.1. And of course it is only the cosine term in (26) that *really matters*, since the amplitude of the sine term is at least a factor $O(\sqrt{t})$ smaller.

Pairing off a couple of neighbouring terms in this way is fine, but computationally this makes only a small difference to operational count necessary for the evaluation of (19). But suppose one is more ambitious and attempts to pair off a whole *collection* of such terms $\alpha = \alpha_E \pm 1, \alpha_E \pm 3 \ldots \alpha_E \pm M_t \geq a$ around a single pivot, where $M_t$ is some, as yet, undetermined odd integer fixing the collection size. The equivalent expression to (26) for such a collection is given by

$$2\sum_{\substack{j=1\\odd}}^{M_t}\int_0^1 e^{i\pi\alpha_E^2 x/4}\frac{\exp[i(t/2)\log\{(1-x)/x\}]}{[x(1-x)]^{1/4}}\left\{\alpha_E\cos\left(\frac{\pi j\alpha_E x}{2}\right)+ij\sin\left(\frac{\pi j\alpha_E x}{2}\right)\right\}e^{i\pi j^2 x/4}dx. \quad (27)$$

Now provided $M_t \ll \alpha_E \Rightarrow Max(\alpha_E M_t, M_t^2) \ll \alpha_E^2$ is satisfied, then all the remarks concerning (26) for the case $j = 1$ should still apply to (27) for $j = 3, 5, \ldots, M_t$ as well. To see if this supposition is true, one must examine a generic integral of the type (27) in more detail.

2.3 *Estimation of an upper bound for the collection size $M_t$.*

The first step is to rewrite integral (27), first by dividing the range of integration into two and then substituting $y = 1 - x$ for the integral over the range $x \in [1/2, 1]$, which gives

$$\int_0^1 e^{i\pi(\alpha_E^2+j^2)x/4}\frac{\exp[i(t/2)\log\{(1-x)/x\}]}{[x(1-x)]^{1/4}}\cos\left(\frac{\pi j\alpha_E x}{2}\right)dx = B_C(\alpha_E, t, j) + e^{i\pi(\alpha_E+j)^2/4}\overline{B_C(\alpha_E, t, j)},$$

(28a, b)

$$\int_0^1 e^{i\pi(\alpha_E^2+j^2)x/4}\frac{\exp\left[i(t/2)\log\left\{\frac{(1-x)}{x}\right\}\right]}{[x(1-x)]^{1/4}}\sin\left(\frac{\pi j\alpha_E x}{2}\right)dx = B_S(\alpha_E, t, j) - e^{i\pi(\alpha_E+j)^2/4}\overline{B_S(\alpha_E, t, j)}.$$

Making the further substitution $x = (w+1)^{-1} \in [0, 1/2]$ (cf. eqs. 19 & 20) one obtains

$$B_C(\alpha_E, t, j) = \int_1^\infty \frac{\exp[i\pi\alpha_E^2/4(w+1) + i(t/2)\log(w)]}{w^{1/4}(w+1)^{3/2}}\cos\left(\frac{\pi\alpha_E j}{2(w+1)}\right)e^{i\pi j^2/4(w+1)}dw,$$

$$B_S(\alpha_E, t, j) = \int_1^\infty \frac{\exp[i\pi\alpha_E^2/4(w+1) + i(t/2)\log(w)]}{w^{1/4}(w+1)^{3/2}}\sin\left(\frac{\pi\alpha_E j}{2(w+1)}\right)e^{i\pi j^2/4(w+1)}dw.$$

(29a, b)

These two integrals can be estimated using similar methods to those used to obtain (23a) for $B(\alpha, t)$. But as will be seen there are some *important* differences associated with the error terms. The integration contour is the same as Fig. 1b. The contribution to (28) originating from integrating around the unit circle remains zero, as explained under (22). Hence the only significant contribution to (28) from (29) arises from the integration through the saddle point at $w = pc(\alpha_E)$, associated with the most rapidly varying $O(\alpha_E^2)$ and $O(t)$ terms of the exponential phase. By contrast the phase terms in $\alpha_E j$ and $j^2$ vary relatively slowly near $pc(\alpha_E)$, since $M_t^2 \ll \alpha_E^2$. Integrating along the contour $w = pc(\alpha_E) + ue^{i\pi/4}$, where $u \in [-(pc+1)/\sqrt{2}, \infty]$, gives the following expression for $B_C(\alpha, t, j)$ (the case for $B_S(\alpha, t, j)$ is essentially identical)

$$B_C(\alpha_E, t, j) = e^{i\pi/4}\int_{-\frac{(pc+1)}{\sqrt{2}}}^\infty e^{d_2 u^2}\left[f(u)e^{i\pi j^2/4(pc+1+ue^{i\pi/4})}\right]\cos\left(\frac{\pi\alpha_E j}{2(pc+1+ue^{i\pi/4})}\right)du. \quad (30)$$

Setting $j = 0$ in (30) recovers $B(\alpha, t)$. In (30) the real coefficient $d_2 = -\frac{t(pc-1)}{4(pc)^2(pc+1)} < 0$,

$$f(u) = \frac{\exp\left[\sum_{i=3}^{\infty} d_i(pc)u^i + i\left\{g(u) - \frac{\arctan\left(u/(pc\sqrt{2}+u)\right)}{4} - \frac{3\arctan\left(u/(pc\sqrt{2}+\sqrt{2}+u)\right)}{2}\right\}\right]}{[(pc)^2 + u^2 + u\sqrt{2}pc]^{1/8}[(pc+1)^2 + u^2 + u\sqrt{2}(pc+1)]^{3/4}},$$

where

$$g(u) = \frac{(\pi\alpha_E^2/4)(pc+1+u/\sqrt{2})}{[(pc+1)^2 + u^2 + u\sqrt{2}(pc+1)]} + \frac{t\log\left((pc)^2 + u^2 + u\sqrt{2}pc\right)}{4}. \quad (31a,b)$$

These expressions (31) were derived in [27 eqs. A49, 50a, 50b], which also lists the values of $d_{i\geq 3}(pc)$, but these are irrelevant to the argument that follows. Now in general, if $\alpha = O(\sqrt{t}) \geq a$ and $pc = O(1)$, then $g(u) = O(t)$ and it is this term that dominates the behaviour of the derivatives of $f(u)$. In particular, one would expect $f'(0) \sim g'(0) = O(t)$ and $f''(0) \sim [g''(0) + (g'(0))^2] = O(t^2)$. However, since *both* $g'(0) = 0$ (as a consequence of the definition of $pc$ given by 12) and $g''(0) = 0$ (by differentiation), it means that both $f'(0)$ and $f''(0)$ are terms of $O(1)$, just like $f(0)$ itself. In fact it is easy to show from these results that $f'(0) = O(t/\alpha^2)f(0) \; \forall \alpha \geq a$. Expanding the term in square brackets in (30) around the point $u = 0$ gives the following

$$f(u)e^{i\pi j^2/4(pc+1+ue^{i\pi/4})} = f(0)e^{i\pi j^2/4(pc+1)} + f(0)u\left\{\frac{-i\pi j^2 e^{i\pi j^2/4(pc+1)}}{4(pc+1)^2} + O(t/\alpha^2)\right\} + O(u^2), \quad (32)$$

which puts (30) in the form

$$B_C(\alpha_E, t, j) \approx e^{i\pi/4}f(0)e^{i\pi j^2/4(pc+1)} \int_{-\frac{(pc+1)}{\sqrt{2}}}^{\infty} e^{d_2 u^2} \cos\left(\frac{\pi\alpha_E j}{2(pc+1+ue^{i\pi/4})}\right) du$$

$$-e^{i\pi/4}f(0)\left\{\frac{i\pi j^2 e^{i\pi j^2/4(pc+1)}}{4(pc+1)^2} + O\left(\frac{t}{\alpha^2}\right)\right\} \int_{-\frac{(pc+1)}{\sqrt{2}}}^{\infty} ue^{d_2 u^2} \cos\left(\frac{\pi\alpha_E j}{2(pc+1+ue^{i\pi/4})}\right) du. \quad (33)$$

The next step is to find an approximation of the first integral in (33) with a prescribed relative error bound. Consider the following generic integral defined by

$$I(A, B, C) = \int_{-X}^{X} e^{-Au^2} e^{\pm iJB/(C+ue^{i\pi/4})} du, \quad \text{with } A > 0 \text{ and } X > 0. \quad (34)$$

The following assignments makes the integrands in (33 & 34) identical:

$$A \equiv |d_2| = \frac{t(pc-1)}{4(pc)^2(pc+1)}, \quad B \equiv \frac{\pi\alpha_E}{2} = \sqrt{\frac{\pi t}{2pc}}(pc+1), \quad C \equiv pc+1. \quad (35)$$

Expansion of the term $(C + ue^{i\pi/4})^{-1} = C^{-1}\left(1 - ue^{i\pi/4}/C + i(u/C)^2 + O(u^3/C^3)\right)$ in (34) gives

$$I(A, B, C) = \int_{-X}^{X} e^{-Au^2} e^{\pm iJB\left(1 - ue^{i\pi/4}/C + i(u/C)^2 + O(u^3/C^3)\right)/C} du. \quad (36)$$

Now define $X = \pi^{-1}(pc + 1)(pc/tM_t^2)^{1/6} \ll C$, where $M_t = \max(j)$ in (27) and employ the assignments (35). Then $\forall |u| \leq X$ the term $jB|u|^3/C^4 \leq jBX^3/C^4 = \pi^{-3}\sqrt{\frac{\pi}{2}}\left(\frac{j}{M_t}\right) < 0.041$ and the exponential of this factor will be very close to unity. Now *suppose* one could fix $M_t$ so that when $|u| = X$, the $Re[-AX^2 \pm iM_tB/(C \pm Xe^{i\pi/4})] = -O(t^{n>0})$. Then the terms $O(u^3/C^3)$ in (36) could be ignored, because they would only become significantly different from unity in the region $u > X$, when the value of $Re[-AX^2 + iM_tB/(C \pm Xe^{i\pi/4})]$ would render the integrand a totally insignificant $e^{-O(t^n)}$ in magnitude. Subject to this supposition, $I(A, B, C)$ could then be approximated using the standard result

$$I(A,B,C) \approx e^{\pm ijB/C} \int_{-X}^{X} e^{-(A \pm jB/C^3)u^2} e^{\mp ijBe^{i\pi/4}u/C^2} du \approx e^{\pm ijB/C} \sqrt{\frac{\pi}{(A \pm jB/C^3)}} \exp\left[-i\frac{j^2B^2}{4C^4(A \pm jB/C^3)}\right],$$

(37)

utilising [15 eq. 3.323(2)], with the proviso $A > jB/C^3$. If this proviso were satisfied, then one can define $\delta = jB/AC^3 < 1$ as a small parameter. In which case (37) itself could be written in the form

$$I(A,B,C) \approx e^{\pm ijB/C} e^{-ij^2B^2/4AC^4} \sqrt{\frac{\pi}{A}} \left\{ \exp\left(\pm i\frac{j^2B^2}{4AC^4}[\delta + O(\delta^2)]\right) \times \left(1 \mp \frac{\delta}{2} + O(\delta^2)\right) \right\}. \quad (38)$$

Studying equation (38) carefully, one is struck by the fact that *if* all the terms involving $\delta$ are indeed small, the term inside { } could be written as $1 + some\ small\ error$. This would leave 'relatively' simple expressions for both $I(A, B, C)$ and the integral $B_C(\alpha, t, j)$ required for (28a, b). Together with (35), the term inside { } would indeed lie close to unity provided

$$\frac{j^2B^2}{4AC^4}\delta = \frac{j^3B^3}{4A^2C^7} \leq \frac{M_t^3 B^3}{4A^2C^7} = \frac{M_t^3}{(pc^2-1)^2}\sqrt{\frac{2\pi^3 pc^5}{t}} < 1,$$

$$\Rightarrow M_t(pc) < \frac{(pc^2 - 1)^{2/3} t^{1/6}}{(2\pi^3 pc^5)^{1/6}}, \quad (39)$$

is chosen as a suitable upper bound on the value of $M_t$.

2.4 *The upper bound on the collection size $M_t$, its implications and validity*

At first glance, the surprising feature of (39) is the *sheer scale* of this upper bound. Here are some simple observations. If $\alpha_E = [1 + O(1)]a \Rightarrow (pc - 1) = O(1)$, then $M_t = O(t^{1/6})$; whilst if $\alpha_E = O\left(t^\eta; \eta \in \left(\frac{1}{2}, 1\right]\right) \Rightarrow pc = O(t^{2\eta-1})$ and $M_t = O(\sqrt{pc}\,t^{1/6}) = O(t^{\eta-1/3})$. Hence (39) suggests that provided $j \leq M_t$, then pairs of integrals $\alpha = \alpha_E \pm j$, which make up the terms of (19), can be first be combined as in (27) and then estimated to some prescribed accuracy, say of $O(\varepsilon_t)$, using the very much simpler approximation given by the leading part of (38). Taken together, all these estimates form a collection, of at least $O(t^{1/6})$ in size, characterised by a

*common, fixed* pivot integer $\alpha_E$. The implications of this last point regarding the sum of such a collection (27), will only become clear in the next section. But given all the assumptions made in section 2.3 to reach this point, are such large collection sizes truly feasible?

Consider the following proposition. Let $M_t(pc)$ satisfy the following approximation (a *precise definition* of $M_t(pc)$ will be prescribed in section 4.2.2 when assessing the computational aspects of the proposition)

$$M_t(pc) \approx \varepsilon_t^{1/3} \frac{(pc^2 - 1)^{2/3} t^{1/6}}{(2\pi^3 pc^5)^{1/6}} \geq 1. \tag{40}$$

Here $\varepsilon_t$ is a small parameter, which for the moment will simply prescribed to satisfy $\varepsilon_t \xrightarrow[t \to \infty]{} 0$, but do so *more slowly* than any power of $t^{-\mu}$ with $\mu > 0$. With this prescription for $\varepsilon_t$, backtrack and consider the validity of all the various assumptions that were necessary to formulate (31-39). (In the remainder of this section, let $\kappa = O(1)$ constant.) The first necessary condition for (38) was that $\delta = jB/AC^3 < 1$. Using approximation (40) one has

$$\delta = jB/AC^3 < M_t B/AC^3 = \varepsilon_t^{1/3} \left(\frac{16 pc^2}{t(pc^2 - 1)}\right)^{\frac{1}{3}},$$

$$= \begin{cases} O\left(\dfrac{\varepsilon_t^{1/3}}{\kappa^{1/6} t^{1/4}}\right) \text{ when } \alpha_E = a + \kappa \Rightarrow pc \approx 1 + 2\sqrt{\kappa}\left(\dfrac{2\pi}{t}\right)^{1/4}, \\ \\ O\left(\left(\dfrac{\varepsilon_t}{t}\right)^{1/3}\right) \text{ otherwise} \end{cases} \tag{41}$$

and hence $\delta \ll 1$ across the entire range of possible pivots $\alpha_E$. This means that the proviso $A > jB/C^3$ is satisfied. Equation (40) also implies that $\delta < j^2 B^2 \delta / 4AC^4 \leq M_t^3 B^3 / 4A^2 C^7 = O(\varepsilon_t)$ from (39), which means in (38) the term $\{\ \} = 1 + O(\varepsilon_t)$. The variable $X = \pi^{-1}(pc+1)(pc/tM_t^2)^{1/6} \ll C = (pc+1)$, introduced after (36) satisfies

$$X = \frac{pc+1}{\varepsilon_t^{1/9}} \left(\frac{2^{1/4} pc^2}{\pi^{15/4}(pc^2-1)t}\right)^{2/9} = \begin{cases} O\left(\varepsilon_t^{-1/9} \kappa^{-2/9} t^{-1/6}\right) \text{ when } \alpha_E = a + \kappa \\ \\ O\left(\varepsilon_t^{-1/9} t^{2\eta - 11/9}\right) \text{ otherwise} \end{cases}$$

$$\Rightarrow Re\left[-AX^2 \pm iM_t B/(C \pm Xe^{i\pi/4})\right] = \begin{cases} -O\left(\varepsilon_t^{-2/9} \kappa^{5/18} t^{5/12}\right) \text{ when } \alpha_E = a + \kappa, \\ \\ -O\left(\varepsilon_t^{-2/9} t^{5/9}\right) \text{ otherwise.} \end{cases} \tag{42}$$

Since this real part is $\ll 0$, the requirement that integral (34) is negligible outside the range $u \in [-X, X]$ will also be met. Taken together these results demonstrate that the supposition regarding exclusion of the terms $O(u^3/C^3)$ from (36) is justified and one can write

$$I(A, B, C) = e^{\pm ijB/C} e^{-ij^2 B^2/4AC^4} \sqrt{\frac{\pi}{A}} (1 + O(\varepsilon_t)),$$

$$= e^{\pm ij\sqrt{\pi t/2pc}} e^{-i\pi j^2 pc/2(pc^2-1)} 2pc \sqrt{\frac{\pi(pc+1)}{t(pc-1)}} [1 + O(\varepsilon_t)], \tag{43}$$

with $A, B$ and $C$ replaced by (35) and $j \leq M_t$. Hence the *first* integral appearing in (33) can be approximated by

$$\int_{-\frac{(pc+1)}{\sqrt{2}}}^{\infty} e^{d_2 u^2} \cos\left(\frac{\pi \alpha_E j}{2(pc + 1 + ue^{i\pi/4})}\right) du = \frac{I(A, B, C) + I(A, -B, C)}{2} + e^{-O\left(\varepsilon_t^{-2/9} t^{5/9}\right)}$$

$$= e^{-i\pi j^2 pc/2(pc^2-1)} 2pc \sqrt{\frac{pc+1}{t(pc-1)}} \cos\left[j \frac{\sqrt{\pi t}}{\sqrt{2pc}}\right] [1 + O(\varepsilon_t)]. \tag{44}$$

The final step is to show that the estimate (44) for the first integral in (33) dominates the *second* integral appearing in (33). The latter is clearly identical to the former except for an extra factor of $u$ in the integrand and so one would expect it to be much smaller than (44). However, the second integral is also multiplied by a factor of magnitude $\pi j^2/4(pc+1)^2$ (the $O(t/\alpha^2)$ factor is smaller still and does not alter the conclusions), which has the *potential* to outweigh the reduction in the integral itself. All the analysis discussed in deriving approximation (44) remains applicable, so the resulting generic analogue of (37) satisfies

$$I_1(A, B, C) = e^{\pm ijB/C} \int_{-X}^{X} u e^{-(A \pm jB/C^3)u^2} e^{\mp ijBe^{i\pi/4}u/C^2} du \approx \left(\frac{\mp ijBe^{i\pi/4}}{2C^2(A \pm jB/C^3)}\right) I(A, B, C). \tag{45}$$

This is obtained from the exact result $\int_{-\infty}^{\infty} x e^{-p^2 x^2 \pm qx} dx = \sqrt{\pi}(q/2p^3)\exp(q^2/4p^2)$ which is the derivative w.r.t. $q$ of [15 eq. 3.323(2)]. So the term involving the second integral in (33) differs in magnitude from (44) by a factor

$$\frac{\pi j^2}{4(pc+1)^2} \times \frac{jB}{2AC^2} = \frac{j^3(\pi pc)^{3/2}}{2\sqrt{2t}(pc-1)(pc+1)^2} \leq \frac{M_t^3 (\pi pc)^{3/2}}{2\sqrt{2t}(pc-1)(pc+1)^2} \approx \frac{\varepsilon_t(pc-1)}{4pc}, \tag{46}$$

using (35) and (40). Hence the second integral in (33) is smaller than the first integral by at least a factor of $\varepsilon_t/4$ for all $j \leq M_t$. This is no more significant than the $O(\varepsilon_t)$ correction already accounted for in (44). One could continue the expansion of (32) to higher powers of $u^n$, but this only gives rise to leading order corrections to the integral $B_C(\alpha_E, t, j)$ of at most $(\varepsilon_t/4)^n/n!$, which are insignificant compared to the $O(\varepsilon_t)$ corrections already obtained. Hence, using (44), (46) and the value of $f(0)$ found from (31), one finds

$$B_C(\alpha_E, t, j) = e^{i\pi/4} \left\{\frac{e^{i\pi \alpha_E^2/4(pc+1)+it\log(pc)/2}}{pc^{1/4}(pc+1)^{3/2}}\right\} e^{i\pi j^2/4(pc+1) - i\pi j^2 pc/2(pc^2-1)}$$

$$\times 2pc \sqrt{\frac{\pi(pc+1)}{t(pc-1)}} \cos\left[j \frac{\sqrt{\pi t}}{\sqrt{2pc}}\right] [1 + O(\varepsilon_t)],$$

$$= \frac{2pc^{3/4}}{(pc+1)} \sqrt{\frac{\pi}{t(pc-1)}} \exp\left[\frac{i\pi}{4}\left\{\frac{\alpha_E^2}{(pc+1)} - \frac{j^2}{(pc-1)} + 1\right\} + \frac{it}{2}\log(pc)\right] \cos\left[j\frac{\sqrt{\pi t}}{\sqrt{2pc}}\right][1 + O(\varepsilon_t)],$$

(47)

provided $j \leq M_t$. This approximation for $B_C(\alpha_E, t, j)$ provides a similar estimate for integral $B_S(\alpha_E, t, j)$, defined by (29b), except the cosine term is replaced by a sine term. Substituting these results into (28) and noting that because $\alpha_E$ is even and $j$ is odd, the term $e^{i\pi(\alpha_E+j)^2/4} \equiv e^{i\pi/4}$, one obtains (upper terms ≡ cosine, lower terms ≡ sine)

$$\int_0^1 e^{i\pi(\alpha_E^2+j^2)x/4} \frac{\exp\left[i(t/2)\log\left\{\frac{(1-x)}{x}\right\}\right]}{[x(1-x)]^{1/4}} \frac{\cos}{\sin}\left(\frac{\pi j \alpha_E x}{2}\right) dx =$$

$$e^{i\pi/8} \frac{pc^{3/4}}{(pc+1)} \sqrt{\frac{\pi}{t(pc-1)}} \frac{4}{4i} \frac{\cos}{\sin}\left(\frac{\pi \alpha_E^2}{4(pc+1)} - \frac{\pi j^2}{4(pc-1)} + \frac{t}{2}\log(pc) + \frac{\pi}{8}\right) \frac{\cos}{\sin}\left[j\frac{\sqrt{\pi t}}{\sqrt{2pc}}\right][1 + O(\varepsilon_t)].$$

(48)

Note that when $j = 0$ one automatically recovers (23), although the relative error in (48) is (as highlighted after 29) *much* larger, since it arises from estimates at $j = M_t$.

### 2.5 $Z(t)$ *composed in terms of Generalized Quadratic Gauss Sums*

Using approximation (48) for the integral (28) one can now consider the sum of such integrals as prescribed in (27). Now there is one small caveat to the results established in 2.2-2.4. The definition of $M_t(pc)$ given by (40) stipulates that $M_t(pc) \geq 1$. Now if $\alpha_E \approx a$ then this may not happen. For instance, if $\alpha_E = a + \kappa \Rightarrow pc \approx 1 + 2\sqrt{\kappa}(2\pi/t)^{1/4}$, this means that one must have $\kappa \geq \pi\varepsilon_t^{-1}/16$ to ensure that $M_t(pc) \geq 1$. Depending upon how $\varepsilon_t$ is actually defined, there are a (arbitrarily) large number of $\alpha_E \in (a, a + \kappa]$ for which (40) is not satisfied. For $\alpha$ values lying in this range, the integrals making up the terms of $Z(t)$ in (11) cannot be collected together and estimated asymptotically in the manner described above, without producing relative errors greater than $O(\varepsilon_t)$. However, this is not a serious problem because the range $(a, a + \kappa]$ forms only an infinitesimal proportion of the total range of $\alpha \in (a, N_\alpha]$ in (11). Hence all one needs to recognize is that there exists a constant, to be denoted by $b \geq \pi\varepsilon_t^{-1}/16$ (as will be discussed later it is computationally expedient to set $b$ very much larger than this minimum value), for which the terms of (11) associated with $\alpha \in (a, a + b]$ can only be be added together one by one, rather than as amalgamative collections of length $M_t$.

With this caveat one can now substitute expression (48) into (27) (N.B. note the additional factor of 2 in front), which in turn can be utilised in (19). The sum in (19) is cut off beyond $\alpha > N_\alpha \geq INT_O(t/\pi) + 2$ and then continued analytically as in (10-11). This gives

$$Z(t) = T_{a+\varepsilon} + \mathcal{H}(t)2\sqrt{2}\left\{\sum_{\substack{\alpha>a\\odd}}^{a+b} \frac{\cos\left(\frac{t}{2}\{log(pc)+\frac{1}{pc}\}+\frac{t}{2}+\frac{\pi}{8}\right)}{(\alpha^2-a^2)^{\frac{1}{4}}}[1+O(\Lambda_\alpha^{-1})] + R_{EM}\right\} +$$

$$\mathcal{H}(t)\left(\frac{\pi}{2^5 t}\right)^{\frac{1}{4}} \sum_{\substack{\alpha_E>a+b\\step\ 2M_t+2}}^{N_\alpha} \frac{8\,pc^{3/4}}{(pc+1)}\sqrt{\frac{\pi}{t(pc-1)}} \sum_{\substack{j=1\\odd}}^{M_t} \alpha_E\left\{CC_{\alpha_E,j} - \frac{j}{\alpha_E}SS_{\alpha_E,j}\right\}[1+O(\varepsilon_t)]. \quad (49)$$

Here $CC_{\alpha_E,j} = \cos(\pi\alpha_E^2/4(pc+1) - \pi j^2/4(pc-1) + t\log(pc)/2 + \pi/8)\cos(j\sqrt{\pi t}/\sqrt{2pc})$ and $SS_{\alpha_E,j}$ the corresponding term with cosine replaced by sine. Now since $\frac{j}{\alpha_E} = \left(\frac{j}{M_t}\right)\frac{M_t}{\alpha_E} \leq \left(\frac{j}{M_t}\right) \times O(t^{-1/3}), \forall\, \alpha_E \in (a+b, N_\alpha]$ and $O(\sum_j CC_{\alpha_E,j}) = O(\sum_j jSS_{\alpha_E,j}/M_t)$, the relative error brought about by the $jSS_{\alpha_E,j}/\alpha_E$ term in (49) is no more than $O(t^{-1/3}) \ll O(\varepsilon_t)$. Substituting (12) for the amplitude $\alpha_E$ in the second sum and incorporating identity (24), one arrives at

$$Z(t) = T_{a+\varepsilon} + \mathcal{H}(t)2\sqrt{2}\left\{\sum_{\substack{\alpha>a\\odd}}^{a+b} \frac{\cos\left(\frac{t}{2}\{log(pc)+\frac{1}{pc}\}+\frac{t}{2}+\frac{\pi}{8}\right)}{(\alpha^2-a^2)^{\frac{1}{4}}}[1+O(\Lambda_\alpha^{-1})] + R_{EM}\right\} +$$

$$\sum_{\substack{\alpha_E>a+b\\step\ 2M_t+2}}^{N_\alpha} \frac{\mathcal{H}(t)4\sqrt{2}}{(\alpha_E^2-a^2)^{\frac{1}{4}}} \sum_{\substack{j=1\\odd}}^{M_t} \left\{\cos\left(\frac{t}{2}\{log(pc)+\frac{1}{pc}\}+\frac{t}{2}+\frac{\pi}{8}-\frac{\pi j^2}{4(pc-1)}\right)\cos\left(j\frac{\sqrt{\pi t}}{\sqrt{2pc}}\right)\right\}[1+O(\varepsilon_t)]. \quad (50)$$

Note the step-size for the first sum in (49-50) is $2M_t + 2$, so as the scale of $\alpha_E$ increases, $M_t$ increases too, and the gaps between the pivots widen. The number of pivots is simply $O(\alpha_E/M_t) = O(t^{1/3})$. Effectively this means the single sum (11) carried out from $\alpha = O(\sqrt{t})$ to $\alpha = O(t)$, has been rewritten as double sum (50) with total summands $O(t^{1/3}M_t) = O(\sqrt{t})$ to $O(t)$. So apart from significantly increasing the relative error from $\Lambda_\alpha^{-1}$ to $\varepsilon_t$, *seemingly nothing of any significance achieved*. *But this is not the case*. The sum over the collections $j = 1, 3, \ldots, M_t$ in (50) can rewritten in the form

$$\frac{1}{2}\mathrm{Re}\left\{e^{-i[t(log(pc)+1/pc+1)/2+\pi/8]}\sum_{\substack{j=1\\odd}}^{M_t} e^{i\pi j^2/4(pc-1)}\left(e^{i\pi aj/4\sqrt{pc}} + e^{-i\pi aj/4\sqrt{pc}}\right)\right\}. \quad (51)$$

At each pivot, define $\omega^\pm(pc) = (\pi/4)[1/(pc-1) \pm a/\sqrt{pc}] - t(\log(pc)+1/pc+1)/2 - \pi/8$ and fix $M_t^- = (M_t-1)/2$. Substituting $(2k+1)$ for $j$ into the above, with $k = 0, 1, 2, \ldots, M_t^-$, a little algebra transforms (51) to

$$\frac{1}{2}\mathrm{Re}\left\{e^{i\omega^+(pc)}S_{M_t^-}^*\left(\frac{1}{pc-1}, \frac{1}{2(pc-1)}+\frac{a}{4\sqrt{pc}}\right) + e^{i\omega^-(pc)}S_{M_t^-}^*\left(\frac{1}{pc-1}, \frac{1}{2(pc-1)}-\frac{a}{4\sqrt{pc}}\right)\right\}, \quad (52)$$

where

$$S_N^*(x,\theta) = S_N(x,\theta) + \frac{(1+e^{i\pi N(N+2\theta)})}{2} = \sum_{k=0}^{N}{}' e^{i\pi k[kx+2\theta]} + \frac{(1+e^{i\pi N(Nx+2\theta)})}{2}, \quad (53)$$

and $S_N(x, \theta)$ is the standard notation for the generalised quadratic Gauss sum. The prime $'$ on the summation sign indicates the first and last terms are halved ($\Rightarrow S_0^* = 1, S_0 = 1/2$), hence the two additional endpoint terms in (53). Substituting (53) into (50) gives

$$Z(t) = T_{a+\varepsilon} + \mathcal{H}(t)2\sqrt{2}\left\{\sum_{\substack{\alpha > a \\ odd}}^{a+b} \frac{\cos\left(\frac{t}{2}\left\{\log(pc) + \frac{1}{pc}\right\} + \frac{t}{2} + \frac{\pi}{8}\right)}{(\alpha^2 - a^2)^{\frac{1}{4}}}[1 + O(\Lambda_\alpha^{-1})] + R_{EM}\right\} +$$

$$\sum_{\substack{\alpha_E > a+b \\ step\ 4M_t^- + 4}}^{N_\alpha} \frac{\mathcal{H}(t)2\sqrt{2}}{(\alpha_E^2 - a^2)^{\frac{1}{4}}} \text{Re}\{e^{i\omega^+(pc)}S_{M_t^-}^*(x(pc), \theta^+(pc)) + e^{i\omega^-(pc)}S_{M_t^-}^*(x(pc), \theta^-(pc))\}[1 + O(\varepsilon_t)], \quad (54)$$

where $x(pc) = 1/(pc - 1)$ and $\theta^\pm(pc) = x/2 \pm a/4\sqrt{pc}$. It is now apparent that there is an asymptotic (if $\varepsilon_t \xrightarrow[t \to \infty]{} 0$ as prescribed) structural connection between the Hardy function and a whole series of sub-sequences based upon pairs of generalised quadratic Gauss sums, each of length $M_t^- = (M_t - 1)/2$ corresponding to half the collection size. A (long established) principle of such Gaussian sums is that their evaluation involves a computational process of only $O(\log(M_t^-))$ in complexity (see below), not the $O(M_t^-)$ as might be assumed from the number of summands. Given that the number of pivots in (54) is $O(t^{1/3})$, one is led to the tentative conclusion that a single arbitrary evaluation of $Z(t)$ is potentially an $O\left(t^{1/3} \times 2\log(M_t^-)\right)$ computational operation for $t \to \infty$. The problem now is to construct an algorithmic scheme which can be formulated into a suitable code for the estimation $Z(t)$, the performance of which comes as close as possible to realising this conclusion.

In his paper [18] Hiary, after re-formulating some original ideas published in [41], reached very much the same conclusion. However, in his work the connection between $Z(t)$ and quadratic Gaussian sums was established directly from (8) rather than (11), by performing a truncated Taylor expansion of the logarithmic part of the phase about some suitable (pivot) integer up to the quadratic terms. The scaling factors in the associated error term introduced at this level of truncation means that the resulting sub-sequences of Gauss sums can extend to a length $O(t^{1/6})$, comparable to the collection size $M_t$ utilised here. So in retrospect the scale on the upper bound predicted by (39) is not so surprising, as the formulations are analogous to one another.

### 3. Efficient evaluation of Generalized Quadratic Gauss Sums

The topic of quadratic and higher order Gauss sums has been the subject of study by many authors over a long period of time. The paper of [3] gives a comprehensive review of the many of the main results. Restricting the discussion to quadratic Gauss sums, the most famous result is that of Gauss himself, who first proved that for any positive integer $q$

$$S_{q-1}^*\left(x = \frac{2}{q}, \theta = 0\right) = \sum_{k=0}^{q-1} e^{2\pi i k^2/q} = \sqrt{q}\frac{(1+i)(1+i^{-q})}{2}. \quad (55)$$

Many alternative proofs were established later by Direchlet, Kronecker, Schur and others. This result was subsequently extended by [37] for rational $x = p/q$, where $p$ and $q$ are relatively prime, to give

$$S_q\left(\frac{p}{q}, 0\right) = e^{i\pi/4}\sqrt{\frac{q}{p}} S_p\left(-\frac{q}{p}, 0\right). \tag{56}$$

Relatively simple proofs of both results (55) and (56) are presented in [30]. This in turn led, by utilising the methods of contour integration, to the discovery of the remarkable *quadratic reciprocity formula* [40] for generalized quadratic Gauss sums

$$S^*_{|c|-1}\left(\frac{a}{c}, \frac{b}{2c}\right) = \sqrt{\left|\frac{c}{a}\right|} e^{i\pi(|ac|-b^2)/4ac} S^*_{|a|-1}\left(-\frac{c}{a}, -\frac{b}{2a}\right), \tag{57}$$

where $a$, $b$ and $c$ are integers with $ac \neq 0$ and $ac + b$ even. Expression (56) is essentially a special case of (57) with $b = 0$, which is known as Schaar's formula [3].

3.1 *A Quadratic Reciprocity approximation for a Generalized Quadratic Gauss Sum*

Expressions (55-57) apply for quadratic Gaussian sums $S^*_N(x, \theta)$ with $(x, \theta)$ restricted to certain rational values, so they cannot be utilised to compute the sums $S^*_{M_t^-}\left(x(pc), \theta^{\pm}(pc)\right)$ which appear in (54) because $(x, \theta)$ are (usually) irrational. However, a more recent paper by Paris [33] provides a means of computing generalised quadratic Gaussian sums for essentially arbitrary parameters. The method, based on the Abel-Plana form of the Euler-Maclaurin summation formula discussed in [32, section 2.10], produces an asymptotic expansion for $S_N(x, \theta)$ in terms of another (shorter) quadratic Gaussian sum (an approximate quadratic reciprocity relation) and a series of terms based around the complementary error function [32, Chap. 7] in $\sqrt{x}$. The crucial result is Theorem 2 of [33], which is reproduced in full below:

*Theorem [Paris].* Let $S_N(x, \theta)$ be quadratic Gaussian sum, $x \in (0, 1), \theta \in [-1/2, 1/2], \xi = Nx + \theta, \varepsilon = \xi - \lfloor\xi\rfloor = \xi - M, \omega = e^{-i\pi/4}, z^-_M = -\omega\varepsilon\sqrt{2\pi/x}$ and $f(N) = e^{i\pi Nx(N+2\theta)}$. Then $S_N(x, \theta)$ satisfies the following asymptotic approximation

$$S_N(x, \theta) = \frac{e^{-\frac{i\pi\theta^2}{x}}}{\omega\sqrt{x}} \sum_{k=1}^{\lfloor\xi\rfloor *} e^{\frac{i\pi k[-k+2\theta]}{x}} - \frac{i}{2\sqrt{\pi}} \sum_{r=0}^{P-1} \Gamma(r + 1/2)(-i\pi x)^r \{f(N)c'_r(\xi) - c_r(\theta)\}$$

$$+ \frac{e^{-\frac{i\pi\theta^2}{x}}}{2\omega\sqrt{x}}\left\{\text{erf}\left(\frac{\omega\xi\sqrt{\pi}}{\sqrt{x}}\right) - \text{erf}\left(\frac{\omega\theta\sqrt{\pi}}{\sqrt{x}}\right)\right\} + \frac{\text{sgn}(-\varepsilon)f(N)e^{(z^-_M)^2/2}}{2\omega\sqrt{x}} \text{erfc}\left(\text{sgn}(-\varepsilon)z^-_M/\sqrt{2}\right) + \tilde{R}'_P. \tag{58}$$

Here erf and erfc are the error and complimentary error functions respectively and sgn is the signum function. The $*$ on the first sum denotes the last term is halved if $\varepsilon = 0$. This sum vanishes if $\xi = Nx + \theta < 1$. If $\lfloor\xi\rfloor = M = 0$, then the erfc term also vanishes. The coefficients $c_r(\theta)$ in the second sum are defined by

$$c_r(\theta) = \sum_{\substack{k=-\infty \\ k \neq 0}}^{\infty} \frac{1}{[\pi(k+\theta)]^{2r+1}}, \qquad \theta \in (-1/2, 1/2). \qquad (59)$$

If $\lfloor \xi \rfloor = M = 0$, then $c'_r(\xi) = c_r(\xi)$. However, these coefficients are singular at all non-zero integer values. So for $M > 0$ the corresponding regularised coefficients $c'_r(\xi)$ have to be substituted into (58). These are defined by removing the singular values of $c_r(\xi)$ to give

$$c'_r(\xi) = c_r(\xi) - \frac{1}{[\pi\varepsilon]^{2r+1}} = c_r(\varepsilon) - \frac{1}{[\pi\xi]^{2r+1}}, \qquad \xi \geq 0. \qquad (60)$$

The second sum in (58) is not actually convergent because $\Gamma(r + 1/2)$ factor ensures that the terms grow without bound, so this sum must be truncated at some *small integer* $P - 1$, much as *Stirling's* series for the logarithm of the Gamma function must be truncated when the remainder terms stop decreasing. Finally the remainder term are comprises the following:

$$\tilde{R}'_P = R'_P(x, \xi) - R_P(x, \theta),$$

$$\text{with } |R_P(x,y)| \leq \frac{\Gamma(P+1/2)}{\pi^2\sqrt{2}} \left(\frac{x}{\pi}\right)^P \sqrt{\frac{\pi}{2}}^* A_P(y), \text{ where } A_P(y) = \sum_{\substack{k=-\infty \\ k \neq -NINT(y)}}^{\infty} \frac{1}{|k+y|^{2P+1}}. \quad (61)$$

(N.B. there is an additional $\sqrt{2/\pi}$ factor in [33 eq. 2.10-1] which appears to be inconsistent with the earlier eqs. 2.2, 2.4-5 of that paper. This is corrected for here by the inclusion of the factor $\sqrt{\pi/2}^*$ in the bound for $|R_P(x,y)|$. This point is addressed again in Section 4.4.2 where the remainder term is subject to greater scrutiny.) The function $A_P(y)$ has a maximum value when $y = 1/2$ and $\lfloor y \rfloor = 0$ ($\lfloor \ \rfloor \equiv NINT$ is rounded to nearest even integer for half integers values) given by

$$Max\{A_P(y)\} = 2^{2P+1} + \frac{2}{(2P)!} \int_0^\infty \frac{v^{2P} e^{-3v/2}}{(1-e^{-v})} dv, \qquad (62)$$

for all $P \geq 1$. (This result is readily established from the integral representation of the Lerch function, see [15 eq. 9.556.]) The integral in (62) dies off rapidly, giving successive maximum values $Max\{A_1(y)\} = 8.2879$, $Max\{A_2(y)\} = 32.2895$, $Max\{A_3(y)\} = 128.1207$ and $Max\{A_4(y)\} = 512.0525$, which approach $2^{2P+1}$ as $P$ gets large. Now if $x = 1/2$ in (58) then the $Max\{|R_P(x,y)|\}$ is at its smallest when $P = 2$. However, if on average $x = 1/4$ then the $Max\{|R_P(x,y)|\}$ is at its smallest when $P = 3$, providing a suitable cut-off value.

The technical details of the proof of this theorem are presented in [33], along with some illustrative calculations. From a computational point of view, the importance of (58) lies in the recursive structure of this (approximate) quadratic reciprocity formula. It means that any $S_N(x, \theta)$ can be written in terms of another, shorter, quadratic Gauss sum, plus a few correction terms. So as [33] points out (although postulated much earlier by [17]), repeated application of (58) together with the utilisation of simple symmetry properties inherent in all such sums, will potentially allow the evaluation of $S_N(x,\theta)$ in far fewer computational operations than the $O(N)$ incurred simply by the summation of $N$ separate terms. To bring this out, it is useful to rewrite equation (58) in the form

$$S_N(x,\theta) = \frac{e^{-\frac{i\pi\theta^2}{x}}}{\omega\sqrt{x}} S_{\lfloor\xi\rfloor}(-1/x, \theta/x) + C(N, x, \theta, P), \qquad (63)$$

where

$$C(N,x,\theta,P) = \frac{e^{-\frac{i\pi\theta^2}{x}}}{2\omega\sqrt{x}}\left[H(\lfloor\xi\rfloor - 1/2)e^{i\pi\lfloor\xi\rfloor\{-\lfloor\xi\rfloor+2\theta\}/x} - 1 + \mathrm{erf}\left(\frac{\omega\xi\sqrt{\pi}}{\sqrt{x}}\right) - \mathrm{erf}\left(\frac{\omega\theta\sqrt{\pi}}{\sqrt{x}}\right)\right]$$
$$+ \frac{\mathrm{sgn}(-\varepsilon)f(N)e^{(z_M^-)^2/2}}{2\omega\sqrt{x}}\mathrm{erfc}\left(\mathrm{sgn}(-\varepsilon)z_M^-/\sqrt{2}\right) - \frac{i}{2\sqrt{\pi}}\sum_{r=0}^{P-1}\Gamma\left(r+\frac{1}{2}\right)(-i\pi x)^r\{f(N)c_r'(\xi) - c_r(\theta)\} + \tilde{R}_P'.$$
(64)

N.B. The factor $e^{i\pi\lfloor\xi\rfloor\{-\lfloor\xi\rfloor+2\theta\}/x}$ vanishes if $\lfloor\xi\rfloor = 0$ or $-1$, hence the appearance of the Heaviside step function $H(\lfloor\xi\rfloor - 1/2)$. If $M = 0$ the alterations discussed in (59-60) still apply.

The aim now is to construct a recursive algorithm based on (63), which will reduce the length of the Gauss sum quickly from $L_1 = N$ to $L_{n+1} = \lfloor\xi_n = L_n x_n + \theta_n\rfloor, n = 1,2,...$ by formulating successive pairs of coefficients $(x_n, \theta_n)$ until $L_n$ reaches some 'small' $\sim O(1)$ termination value. The best way to achieve this is to exploit the symmetry properties inherent in $S_N(x,\theta)$ and shift the range of $x_n$ from $(0, 1)$, as specified in (58), to $(-1/2, 1/2)$ and then apply (63-64) for $|x_n| \in (0, 1/2)$. With $L_{n+1} = \lfloor\xi_n = L_n|x_n| + \theta_n\rfloor$, the length of the Gauss sum can be reduced from $L_1 = N$ to $L_{n+1} \in [1,2]$ in at most $n_{max} = \lfloor log(N)/log(2)\rfloor$ iterations. However, the cost of computing successive correction values of $C(L_n, |x_n|, \theta_n, P)$ down to a value of $L_{n+1} \in [1,2]$ for instance, would be rather inefficient. A much better strategy is to estimate the average computational cost of a single $C(L_n, |x_n|, \theta_n, P)$ valuation in terms of basic arithmetical operations, link this operational count to the equivalent cost of computing a Gaussian sum $S_K(x,\theta)$ of length $K$ and then terminate the iteration procedure when $L_{n+1} < K$. Termination would then occur after at most $n_K = \lfloor log(N/K)/log(2)\rfloor + 1$ iterations. The value of $K$ will be a constant (in practice $\sim O(10-100)$) because the cost of computing of $C(L_n, |x_n|, \theta_n, P)$ comes down to the cost of computing the various error function terms which are essentially independent $|L_n|$. The factors influencing the best choice of *termination constant K* are discussed in Section 3.6. Having computed the successive $(L_{n+1}, |x_{n+1}|, \theta_{n+1})$ parameters for $n = 1,..,n_K$, one initially computes $S_{L_{n_k+1}}(|x_{n_k+1}|, \theta_{n_k+1})$ exactly, and then repeatedly applies (63) to reach an estimate for the sum $S_{L_1}(|x_1|, \theta_1)$ actually desired. So to achieve a workable and reliable algorithm one needs to devise a specific iterative scheme for the calculation of $(L_{n+1}, |x_{n+1}|, \theta_{n+1})$ parameters, establish bounds both for the value of $K$ and the resulting error in the estimate for $S_{L_1}(|x_1|, \theta_1)$.

3.2 *A recursive scheme for the efficient computation of a Generalized Quadratic Gauss Sum*

In what follows let $S_N^{s=+1}(x,\theta) \equiv$ standard quadratic Gauss sum and $S_N^{s=-1}(x,\theta) \equiv$ its complex conjugate. Consider the calculation of $S_N(x,\theta)$, where initially $(x,\theta)$ are just an arbitrary pair of real numbers. Let $L_1 = N$ and define $x_1 = x - \lfloor x\rfloor \in (-1/2, 1/2]$. Now if $\lfloor x\rfloor$ is even, one can simply substitute $x_1$ for $x$ in the sum, but if $\lfloor x\rfloor$ is odd an extra quadratic factor $i\pi k^2$ is effectively introduced into the exponent of all the summands of $S_N(x,\theta)$. This can be dealt with by replacing the quadratic factor by a corresponding linear factor, since $e^{i\pi k^2} = (-1)^k = e^{i2\pi(k/2)}$, which can be achieved by adjusting the value of $\theta_1$ by an appropriate half factor. At the same time one wishes to be able to utilise *Theorem[Paris]*, which specifies that $x > 0$ and $-1/2 \leq \theta \leq 1/2$, which can be easily achieved by defining

$$\theta_{1\pm} = \text{sgn}(x_1)\{\theta - \lfloor\theta\rfloor \pm (1-(-1)^{\lfloor x \rfloor})/4\} \Rightarrow \theta_{1\pm} \in (-1/2,\ 1/2] \tag{65}$$

and the choice of $\pm$ is always made so that $|\theta_{1\pm}|$ lies inside the range $[0, 1/2]$. Then

$$S_N(x,\theta) = S_{L_1}^{s_1}(|x_1|,\theta_{1\pm}), \qquad x_1, \theta_{1\pm} \in (-1/2,\ 1/2]. \tag{66}$$

where $s_1 = \text{sgn}(x_1)$. Now from (63) one can see that the quadratic parameter $x$ of the primary Gauss sum is replaced by a factor $-1/x$ in the secondary Gauss sum. The change of sign can be easily handled by incorporating a switch from a standard $\leftrightarrow$ complex conjugate Gauss sum into the recursion procedure. This suggests computing a hierarchical chain of primary and secondary Gauss sums linked together formulaically

$$S_{L_n}^{s_n}(|x_n|,\theta_n) \mapsto \begin{cases} S_{L_{n+1}}^{s_{n+1}=+1}(|x_{n+1}|,\theta_{n+1\pm}) & \text{if } x_{n+1} > 0 \\ & \qquad n \geq 1. \\ S_{L_{n+1}}^{s_{n+1}=-1}(|x_{n+1}|,\theta_{n+1\pm}) & \text{if } x_{n+1} < 0 \end{cases} \tag{67}$$

To keep track of the sign changes introduced into the linear parameter $\theta \mapsto \theta/x$ by this procedure it is easiest to switch from $\theta \mapsto \pm\{\theta/|x|\}$ depending upon whether $x \gtrless 0$. So starting from (66) one can apply the following algorithm to estimate $S_N(x,\theta)$.

3.3 *An algorithm for the recursive computation of a Generalized Quadratic Gauss Sum using Theorem [Paris]*. (*algorithm QGS*)

*Replace $S_N(x,\theta)$ by $S_{L_1}^{s_1}(|x_1|,\theta_{1\pm})$, with $L_1 = N, x_1 = x - \lfloor x \rfloor$ and $\theta_{1\pm}$ defined by equation (65) above. Fix the parameters $(K > 10, P)$ with $P$ a small positive integer.*

*While $L_n > K$, compute for $n = 1, 2, \ldots, n_K$*

$$x_{n+1} = -\left\{\frac{1}{|x_n|} - \left\lfloor\frac{1}{|x_n|}\right\rfloor\right\}, \qquad \Rightarrow x_{n+1} \in (-1/2,\ 1/2],$$

$$L_{n+1} = \lfloor L_n|x_n| + \theta_{n\pm}\rfloor, \qquad s_{n+1} = \text{sgn}(x_{n+1}),$$

$$\theta_{n+1\pm} = s_{n+1}\left\{\frac{\theta_{n\pm}}{|x_n|} - \left\lfloor\frac{\theta_{n\pm}}{|x_n|}\right\rfloor \pm (1-(-1)^{\lfloor 1/|x_n|\rfloor})/4\right\},$$

*with the $\pm$ sign chosen to ensure $|\theta_{n+1\pm}| \in [0,\ 1/2]$.*

*Provided either i)* $L_{n_K} > K > L_{n_K+1} \geq 10$, *or ii)* $L_{n_K} > 3K > 10 > L_{n_K+1}$, *compute* $S_{L_{n_K+1}}^{S_{n_K+1}}(|x_{n_K+1}|, \theta_{n_K+1\pm})$ *exactly. If* $L_{n_K+1} = 0$ *or* $L_{n_K+1} = -1$, *then* $S_{L_{n_K+1}} = 1/2$.
*If iii)* $3K \geq L_{n_K} > K > 10 > L_{n_K+1}$ *compute* $S_{L_{n_K}}^{S_{n_K}}(|x_{n_K}|, \theta_{n_K\pm})$ *exactly.*
*From this starting point compute, for* $n = n_K, n_K - 1, \ldots, 1$, *the following iterations*

$$S_{L_n}^{S_n}(|x_n|, \theta_{n\pm}) = \left\{ \frac{e^{-\frac{i\pi(\theta_{n\pm})^2}{|x_n|}}}{\omega\sqrt{|x_n|}} S_{L_{n+1}}(|x_{n+1}|, \theta_{n+1\pm}) + C(L_n, |x_n|, \text{sgn}(x_n)\theta_{n\pm}, P) \right\}^{S_{n+1}}.$$

*The final iteration* $S_{L_1}^{S_1}(|x_1|, \theta_{1\pm})$ *gives an estimate to* $S_N(x, \theta) = S_{L_1}^{S_1}(|x_1|, \theta_{1\pm}) + \delta(K, N, x, \theta, P)$, *where* $\delta(K, N, x, \theta, P)$ *is the error term. The magnitude of this and the corresponding relative error, denoted by* $\varepsilon_{GS}(K)$, *are discussed in Section 3.5.*

Note that although each iterate $x_n \in (-1/2, \ 1/2]$, the algorithm invokes *Theorem [Paris]* at step *4* for $|x_n|$, which conforms with specified condition that $x \in (0, 1)$. Obviously if ever $x_n$ happened to be exactly zero, one would terminate the algorithm immediately and compute the corresponding linear Gauss sum as a geometric series. In step 3 an upper bound on the size of $L_{n_K+1} \geq 10$ is enforced, where realistically possible. Although *Theorem [Paris]* still gives a valid approximation for Gaussian sums with $L_{n_K+1} < 10$, it is in nature an *asymptotic* approximation most suitable 'large' summation values. Consequently, it is better to avoid commencing with an initial sum of length $L_{n_K+1} < 10$, as this can produce a relatively large error in the first iteration, which then percolates through into a relatively large error $\delta(K, N, x, \theta, P)$ in the final answer. Provided $L_{n_K}$ is not too large, say $3K \geq L_{n_K} > K$, then it is better to commence with an exact calculation of the initial sum with length $L_{n_K}$. Sometimes this bound is impossible to enforce, as when the final iteration jumps from a potentially very large $L_{n_K} > 3K$ to a very small $L_{n_K+1}$. Then commencement with a very small initial sum becomes unavoidable. The algorithm requires a storage capacity of $O(n_K \log(N)) = O((\log(N))^2)$ bits for the contents of the arrays $x_n, \theta_{n+1\pm}$ and $L_n$.

*3.4 Sample Computations using algorithm QGS.*

In order to get an idea for the workings of *algorithm QGS*, it is worth performing some sample calculations of quadratic Gauss sums. The issues highlighted will help to estimate its computational efficiency and accuracy. All the following calculations of $S_N(x, \theta)$ are for $N = N_0 = 129901233$. This is an unremarkable 9 digit composite integer, large enough to highlight the main features of *algorithm QGS*, but not so large to make an accurate term by term check computation of $S_{N_0}(x, \theta)$ a really time consuming computational exercise. Values for the termination constant $K = 20$ and cut-off integer $P = 3$ were also used throughout.

**Case A** $S_{N_0}(x, \theta)$, with $x = 1/\sqrt{45} = 0.149071198\ldots$ and $\theta = 1 - \sqrt{23/71} = 0.43083951\ldots$

The results of applying the computational algorithm for this sum are presented in Table 3.4A. At the top of each table are shown the respective calculations of $(x_n, \theta_{n\pm}, s_n)$ from steps 1 and 2. Underneath are shown the respective iteration calculations of the intermediate Gaussian sums, starting from the exact initial value for the sum over $L_{n_K+1}$ until a final estimate for $S_{N_0}(x, \theta)$ is reached, steps 3 to 5. Next to these values are the corresponding term by term calculations (accurate to the number of decimal places shown in the Table). The associated absolute and relative errors in these estimates are also shown.

In this first example $x_1$ is chosen to be a quadratic irrational $\{p + q\sqrt{d}, \ p, q \in \mathbb{Q}, \ d \in \mathbb{N}, \sqrt{d} \notin \mathbb{N}\}$, which means by the Euler-Lagrange theorem (see for example [8], and references therein] that it possesses a *positive* continued fraction representation that eventually becomes periodic. In this instance $1/\sqrt{45} = [0; 6, \overline{1,2,2,2,1,12}]$, where in the standard notation

$$[a_0; \ a_1, a_2, \overline{a_3, a_4}] = a_0 + \cfrac{1}{a_1 + \cfrac{1}{a_2 + \cfrac{1}{a_3 + \cfrac{1}{a_4 + \cfrac{1}{a_3 + \cfrac{1}{a_4 + 1 \ldots}}}}}}, \qquad (a_{1,2,3,\ldots} \geq 1) \qquad (68)$$

and $\overline{a_3, a_4}$ denotes the composition of the periodic sequence of positive partial quotients. This means that the quadratic sum coefficients $x_2, x_3, x_4 \ldots$ will also demonstrate periodic behaviour, with a period less than or equal in length to that of the continued fraction representation of $x_1$. The properties of the $x_2, x_3, x_4 \ldots$ can be deduced from their *nearest integer* continued fraction representations (because they themselves represent differences from a nearest integer). These nearest integer representations can be computed quite easily from any positive continued fraction representation of $x_1 \in (-1/2, 1/2]$, by employing the following procedure.

1. Move along the representation $[0; a_1 \neq 1, a_2, \ldots]$ (or $[-1; 1, a_2, \ldots]$ if $x_1 < 0$) from left to right and whenever a 1 is encountered replace it by a 0 and add one to the two adjacent digits either side. So $[0; 6, \overline{1,2,2,2,1,12}]$ becomes $[0; 7, \overline{0,3,2,3,0,14}]$.
2. From the initial zero again move along the representation left to right. On encountering the second zero, change the sign of all the subsequent non-zero digits bracketed by zeros two and three. Between zeros three and four leave the signs alone, but between zeros four and five change the signs again and repeat this pattern. Then concertina out all the zeros (except the first), which are irrelevant.
   So $[0; 7, \overline{0,3,2,3,0,14}]$ becomes $[0; +7, \overline{0, -3, -2, -3,0, +14}] = [0; 7, \overline{-3, -2, -3, 14}]$, which is the nearest integer continued fraction representation of $1/\sqrt{45}$. Notice that in any nearest integer continued fraction the partial quotients satisfy $|a_n| \geq 2$. The procedure can also be applied to non-periodic positive representations, e.g. $-1/\pi = [-1; 1,2,7,15,1,292,1,1,1,3,1,14,2,\ldots] \to [0; -3, -7, -16, 294, -3,5, -15, -3, \ldots]$, or $e = [2; (1,2k,1)_{k=1}^{\infty}] \to [3; -4, (-1)^k(2, 2k+1)_{k=2}^{\infty}]$.

The nearest integer continued fraction representations of the subsequent $x_{2,3,...}$ iterates can now be read off by multiplying by $-\text{sgn}(a_1)$, deleting $a_1$ and shifting the remaining $a_n$ one place to the left (so $a_2$ becomes $a_1$ etc). In this case $x_2 = [0; \overline{3,2,3,-14}]$, $x_3 = [0; \overline{-2,-3,14,-3}]$, $x_4 = [0; \overline{-3,14,-3,-2}]$ and $x_5 = [0; \overline{14,-3,-2,-3}]$, before the pattern repeats itself. The period length of four is a result of two zeros being concertinaed out of the original period six pattern exhibited by the positive continued fraction for $x_1$. The values of $x_{2-5}$ are shown in their quadratic irrational forms in Table 3.4A.

By contrast the $\theta_{n\pm}$ iterates display no such pattern and essentially form a sequence of pseudo random numbers (although as $\theta_1$ is also a quadratic irrational it seems likely a periodic pattern might, eventually, establish itself). The sequence terminates after $n_K = 11$ iterations, when $L_{12} = 22$. From the exact calculation of $S_{22}^{-1}(x_{12}, \theta_{12\pm})$ one achieves an estimate for $S_{N_0}(x, \theta)$ in a CPU time $\sim O(10^{-2}$ sec) compared to $\sim O(10^3$ sec) for the corresponding term by term computation. Note that computation of the latter is highly prone to round off error arising from the computation of the angular phase, so it essential that the values $(x_n, \theta_{n\pm})$ are stored to an accuracy of at least $3\log(N)/\log(10)$ digits. Special care must also be taken when naively applying the iteration scheme for the $x_n$, as it is numerically unstable and insignificant initial errors build up rapidly after a few repetitions. If $N$ were an integer consisting of millions of digits, the indirect algorithm of [38], devised for computing large scale continued fraction representations, could be utilised to calculate the $(x_n, \theta_{n\pm})$ values speedily and accurately. For the very much more modest calculations presented here it is sufficient to apply *algorithm QGS* as it stands, with $(x_n, \theta_{n\pm})$ computed to an accuracy of 30 digits.

Notice how the magnitude of the error builds up after each successive iteration, but only does so in conjunction with the overall size of the sum. This is illustrated by the behaviour of the relative error which, after some initial fluctuations, stabilises. The basis of this result is due to the relative sizes of the terms appearing in the approximate quadratic reciprocity formula (63). At the $m$th iteration the intermediate Gaussian sum is $O\left([|x_{n_K}||x_{n_K-1}|...|x_m|]^{-1/2}\right)$, whilst the additive factor (64) is no more than $O(|x_m|^{-1/2})$ and so is much smaller. Hence any error in the correction term quickly becomes negligible in comparison with the error in the intermediate Gaussian sum. As both the intermediate Gaussian sum and its error are subsequently scaled by factors of $\sqrt{|x_{m-1}|}$, $\sqrt{|x_{m-2}|}$,..., the relative error plateaus to a constant after just a few iterations (see next section). For example, multiplying 18.946, the final error in Table 3.4A, by $\sqrt{|x_1|} = 45^{-1/4}$ gives 7.3150, the error in the final intermediate Gaussian sum. Likewise $7.3150 \times \sqrt{|x_2|} = 3.9514$ etc, a pattern that only starts to break down when $n = 8$ is reached.

Case B $S_{N_0}(x, \theta)$, with $x = 1 - e/\pi = 0.13474402...$ and $\theta = 1/e = 0.36787944...$

The results for this calculation are shown in Table. 3.4B. In this instance the values of $(x, \theta)$ are not quadratic irrationals and so their corresponding positive continued fractions representations are non periodic. Hence the corresponding iterates $(x_n, \theta_{n\pm})$ form a series of pseudo random numbers, never conforming to any pattern. Again as in Table 3.4A the sequence terminates when $n_K = 11$, with $L_{12} = 60$ on this occasion. This is an example where the next iteration would give rise to a relatively steep drop from $L_{12} = 60$ to $L_{13} = 1$ as $|x_{12}| = 0.01506 \sim O(1/L_{12})$. Since the computation of the exact initial sum of length $L_{12} = 60 = 3K$ with no error, can be achieved almost as quickly as an estimate derived by iteration from a sum of length $L_{13} = 1$, it makes sense to terminate step 2 of the algorithm at this point.

The average number of iterations required by the algorithm is related to the growth of the product of partial quotients $|a_1 a_2 \ldots a_n|$ of the nearest integer continued fraction of $x_1$. For the positive continued fraction the behaviour of this product for large $n$ has been known since [24], see also [36] and [1], who proved that for generic irrationals $x = [0; a_1 a_2 \ldots a_n]$, (more specifically those irrationals that have unbounded partial quotients which follow no specific pattern, so excluding countable subsets of $\mathbb{R}$ such as the quadratic irrationals and numbers like $e$)

$$[a_1 a_2 \ldots a_n]^{1/n} \xrightarrow[n \to \infty]{} K = \prod_{m=1}^{\infty} \left(1 + \frac{1}{m(m+2)}\right)^{\log(m)/\log(2)} = 2.685452001 \ldots \quad (69)$$

where K is Khintchine's constant. This result is based on the fact that in the limit of large $n$ the number of occurrences of the positive integer $m$ in the positive continued fraction representation of $x$ follows the Gauss-Kuzmin distribution ([25], [26], see also [8])

$$P(a_n = m) = \frac{\log(1 + 1/m(m+2))}{\log(2)}, m = 1,2,\ldots \quad (70)$$

The above formula can be used to approximate the corresponding distribution for nearest integer continued fractions. The conversion procedure from positive to nearest continued fractions outlined above, removes all the 1 digits and converts most of them to zeros (not all because a small proportion become $\pm 3$ and an even smaller proportion convert to $\pm 2$). Assuming that all 1 digits were to convert to zeros entirely, the partial quotients of a nearest integer continued fraction would satisfy

$$P(a_n = m) \approx \frac{\log(1 + 1/m(m+2))}{\log(3/2)}, m = 2,3,\ldots$$
$$|a_1 a_2 \ldots a_n|^{1/n} \xrightarrow[n \to \infty]{} \approx K^{\log(2)/\log(3/2)} = 5.41265167 \ldots \quad (71)$$

Hence the average number of iterations $\langle n_K \rangle$ needed by the *algorithm QGS* to reduce a Gaussian sum of length $N$ to a length $K$ starting from a generic irrational $x$ should be about $y \log(N/K)$, where $y = \frac{1}{\log(5.41265..)} = 0.592 \ldots$, to a first order approximation. (A more precise calculation of $P(a_n = m)$ for nearest integer continued fractions gives a figure of $5.3081 \ldots$, when $y = 0.59907 \ldots$) For the sum of length $N_0 = 129901233$ and $K = 20$ computed here, this gives a termination value of $\langle n_K \rangle \approx 9$, which agrees well with the value of $n_K = 11$ found for $x = 1 - e/\pi$.

The behaviour of the error estimates presented in Table 3.4B are also very similar to those shown in Table 3.4A. The error grows with each iteration, but as before the relative error after some small initial fluctuations, settles down to a constant value. The scale of the final relative error differs little from that of the initial value.

Case C $S_{N_0}(x, \theta)$, with $x = \sqrt{2}/10 = 0.141213562\ldots$ and $\theta = \sqrt{10/71} = 0.375293312\ldots$

The two previous cases illustrate that the final relative error in the estimate of a Gaussian sum of this length is typically of $O(10^{-3})$. However, it is quite easy to find $(x, \theta)$ values for which the algorithm gives more precise results. If $x$ is chosen to be to the quadratic irrational $\sqrt{2}/10 = [0; 7, \overline{14}]$, all the following coefficients $x_{2,3,..K} = [0; \overline{14}]$ and the length of each intermediate sum is reduced rapidly by a factor close to 14, bringing about termination after just $n_K = 6$ iterations. Here the initial intermediate sum computation produces an estimate which has a relative error of only about $10^{-6}$ rather than the $10^{-3}$ found in the previous examples. The size of this initial error is due to the fact that the value of $|x_{n_K}| = 0.071..$ is smaller than the corresponding $|x_{n_K}| = 0.427\ldots$ & $0.496\ldots$ values found earlier, which naturally reduces the error term (61). This comparatively small initial error percolates through the rest of the calculation, to produce a final estimate still accurate to six significant figures. This example illustrates that the accuracy of the Gaussian sum estimate will be improved if it turns out that the algorithm terminates immediately after a relatively small $|x_{n_K}|$ value. In one off calculations, it could be worthwhile to allow small changes in $K$ to achieve this.

Case D $S_{N_0}(x, \theta)$, with $x = 0.332613390928725685017 4\ldots$ and $\theta = \frac{1}{2e} = 0.183939720\ldots$

Case E $S_{N_0}(x, \theta)$, with $x = 1/2 - \sqrt{\pi}/N_0^2$ and $\theta = 1/N_0\pi$

In both these cases the initial value of $x$ is a 'non-generic' irrational, specially chosen so that its nearest integer continued fraction representation has a large partial quotient relatively early in its expansion. In Case D, $x = [0; 3, 154, -275596610848, 8, 3, -2, \ldots]$ and in Case E, $x = [0; 2, 2380080351076487, -3, \ldots]$. More specifically, a quotient satisfying $|a_n| > L_n \gg K$, occurring within the first $n \sim 0.59 \log(N/K)$ terms or so of the expansion. Statistically, the chances of encountering such a large partial quotient so early in the expansion are very small, since the probability that $|a_n| > L_n \gg K$ over such a small range is only $O(L_n^{-1})$. (When $n \sim O(\log(N/K))$ and $L_n \to O(K)$, such occurrences become more likely, e.g. Case B above.) Nevertheless it is a possibility, and the behaviour of the algorithm under such scenarios is interesting. In Case D the presence of the very large quotient $|a_3| = 275596610848 > N_0$ causes the algorithm to immediately terminate (Table 3.4D) after just $n_K = 3$ iterations. Notice too that the first intermediate sum does not scale up as $1/\sqrt{x_{n_K}}$ as in the previous cases. This is because the additive factor $C(L_n, |x_n|, \text{sgn}(x_n)\theta_{n\pm}, P)$ given by (64) approaches $(-1/2)$ times the multiplicative factor in (63). (Since $L_3 x_3 \ll |\theta_3|$, $\xi \approx \theta_3$, the two erf terms offset to $O(1)$, whilst the term involving erfc vanishes when $\lfloor \xi \rfloor = M = 0$.) So the usual $1/\sqrt{x_{n_K}}$ scaling cancels out completely. Hence the first intermediate sum is very much smaller in magnitude than would be normally be expected given $L_3 = 280564$. This small intermediate value percolates through the remaining iterations to give an unusually small value for the full sum, much less than the $O(\sqrt{N})$ behaviour seen normally. Case E is slightly different, although as in Case D the algorithm terminates very quickly (Table 3.4E) after just $n_K = 2$ iterations. But on this occasion the $1/\sqrt{x_{n_K}}$ scaling *is present*, which means that the final sum is $O(N)$, very much larger than the standard $O(\sqrt{N})$ behaviour. (The reason lies in the fact that although $L_2 x_2$ is very small, $\theta_2$ is also close to zero and which means that $\xi \not\approx \theta_2$. Hence the cancellation effect seen in Case D does not come about, because the two erf terms no longer offset.) The final sum is $O(N)$ as a consequence of the original $(x, \theta) \approx (1/2, 0)$, and so its magnitude lies close to the exact value $S_{N_0}(1/2, 0) = N_0 e^{i\pi/4}/\sqrt{2}$.

In both these examples the regular behaviour of the relative error seen in Cases A-C is disrupted, because of the abrupt termination of the algorithm. In each case the relative error commences at very small value (see section 3.5) and does not get the chance to stabilise. In Case D the magnitude of the intermediate sums has not built up sufficiently to dominate the additive factors $C(L_n, |x_n|, \text{sgn}(x_n)\theta_{n\pm}, P)$ and so stabilise the error. If the same initial $x$ was used in conjunction with a larger $N$ value, the regular behaviour of the relative error and the $O(\sqrt{N})$ estimate for the sum size would gradually reassert itself. In Case E the effect of increasing $N$ would be offset by tighter approximation of $x \approx 1/2$, so the $O(N)$ scaling and error disruption would remain features of the results.

*3.5 Error Bounds*

Application of *algorithm QGS* gives an approximation to the original quadratic Gaussian sum $S_N(x, \theta) \approx S_{L_1}^{S_1}(|x_1|, \theta_{1\pm})$, with an error $\delta(K, N, x, \theta, P)$ which must be bounded. Initially consider the simplest case when $x \equiv x_1$ is a 'generic' irrational, characterised by a lack of any large partial quotients *early on* in its nearest integer continued fraction representation. Specifically *early on* means a partial quotient $|a_n| > L_n \gg \log(N/K)$ within the first $n \in [1, \langle n_K \rangle \sim 0.59 \log(N/K)]$ iterations needed to reduce $N$ to a number less than $K$. So the initial sum, computed at either steps *3i)* or *3iii)* of the algorithm, is either of length $L_{n_K+1} = O(K)$ or $L_{n_K} = O(K)$ respectively. In practice that means most irrationals with expansions following Gauss-Kuzmin statistics, most quadratic irrationals (unless they have a very large $|a_n|$ early in their periodic pattern) and most irrationals with regular, but non-periodic expansions, such as $e$, where the really large partial quotients are only encountered at positions $n \gg n_K$. But irrationals similar to the example discussed in Cases D & E above would be excluded.

Both 'generic' and 'non-generic' irrationals form uncountable sets. However, members of the former will tend to be much more common than the latter. One can make a rough estimate of the probability that an irrational is 'generic', in the specific sense set out above, using the fact that almost all irrationals obey Gauss-Kuzmin statistics. Choose values for $N$ and $K$, with $N \gg K$. Construct an irrational $x = [0; a_1, a_2, \ldots, a_n, \ldots]$, by means of a series of $n_K$ independent random selections of partial quotients, satisfying distribution (71). Then apply the constructed $x$ to reduce $N = L_1 > L_2 > \cdots > L_{n_K}$ in the manner prescribed in Step 2 of the algorithm and cease reduction if $|a_n| > L_n$. The probability that this will occur, $P(|a_n| > L_n) \approx 1/\{L_n \log(3/2)\}$, can be estimated from (71). If one further assumes that each iteration has - up to this point - reduced $N$ by an average factor $p \sim e^{-1/0.59} = 0.183$

(commensurate with the estimate that it takes $\langle n_K \rangle \sim 0.59 log(N/K)$ iterations to go from $N$ to $L_{n_K} \sim K$), then $L_n \sim p^{n-1} N$. In which case the probability that the constructed $x$ will be 'generic' and the reduction process reaches $n = \langle n_K \rangle$, satisfies

$$\sim \left[ \prod_{n=1}^{0.59 log(N/K)} \left( 1 - \frac{1}{p^{n-1} N \left\{ log\left(\frac{3}{2}\right) \right\}} \right) \right] \xrightarrow[N \to \infty]{} e^{-1/[(1-p)K log(3/2)]} \approx e^{-3.02/K}. \quad (72)$$

So if the termination constant $K \sim O(100)$, then this 'generic' probability would be close to one.

If $x \equiv x_1$ is a 'generic' irrational, then the initial sum in the iteration process computed exactly at step 3, will be of length $O(\sqrt{K})$. From this exact initial sum, step 4 of the algorithm computes successive approximations to the various intermediate sums. Let $\varepsilon_{n_K - m}$ represent the respective error term associated with the $(n_K - m - 1)th$ iteration, defined by

$$\varepsilon_{n_K - m} = S_{L_m}^{S_m}(|x_m|, \theta_{m\pm}) - \left\{ \frac{e^{-\frac{i\pi(\theta_{m\pm})^2}{|x_m|}}}{\omega \sqrt{|x_m|}} S_{L_{m+1}}(|x_{m+1}|, \theta_{m+1\pm}) + C(L_m, |x_m|, sgn(x_m)\theta_{m\pm}, P) \right\}^{S_{m+1}} \quad (73)$$

for $m = n_K, n_K - 1, \ldots, 2, 1$ (if the initial sum starts at step 3 iii) then replace $n_K$ by $n_K - 1$ in the above and what follows). The first iteration, when $m = n_K$, approximates the sum $S_{L_{n_K}}^{S_{n_K}}(|x_{n_K}|, \theta_{n_K \pm})$ to within an error $\varepsilon_0$. The second iteration uses this approximation to estimate the next intermediate sum $S_{L_{n_K - 1}}^{S_{n_K - 1}}(|x_{n_K - 1}|, \theta_{n_K - 1 \pm})$ to within an error $\left[ |\varepsilon_0| / |x_{n_K - 1}| \right]^{1/2} + |\varepsilon_1| \right]$. Continuing in this vein, one finds that the final error $\delta(K, N, x, \theta, P)$ satisfies

$$|\delta(K, N, x, \theta, P)| \leq \frac{|\varepsilon_0|}{\sqrt{|x_{n_K - 1} x_{n_K - 2} \ldots x_1|}} + \frac{|\varepsilon_1|}{\sqrt{|x_{n_K - 2} \ldots x_1|}} + \cdots + \frac{|\varepsilon_{n_K - 2}|}{\sqrt{|x_1|}} + |\varepsilon_{n_K - 1}|,$$

$$\leq \frac{max(|\varepsilon_{n_K - m}|)}{\sqrt{|x_{n_K - 1} x_{n_K - 2} \ldots x_1|}} \left\{ 1 + \sqrt{|x_{n_K - 1}|} + \sqrt{|x_{n_K - 1} x_{n_K - 2}|} + \cdots + \sqrt{|x_{n_K - 1} x_{n_K - 2} \ldots x_1|} \right\}$$

$$(74)$$

Since $|x_i| \leq 1/2$ the sum in the curly brackets of (74) will always be less than $\sqrt{2}/(\sqrt{2} - 1) \approx 3.41$, and one deduces that

$$|\delta(K, N, x, \theta, P)| < \frac{3.41 max(|\varepsilon_{n_K - m}|)}{\sqrt{|x_1 x_2 \ldots x_{n_K - 1}|}} \approx 3.41 max(|\varepsilon_{n_K - m}|) \left| \frac{S_{L_1}^{S_1}(|x_1|, \theta_{1\pm})}{S_{L_{n_K}}^{S_{n_K}}(|x_{n_K}|, \theta_{n_K \pm})} \right|$$

$$\lesssim \frac{3.41 |max(|\varepsilon_{n_K - m}|) S_{L_1}^{S_1}(|x_1|, \theta_{1\pm})|}{\sqrt{K}}. \quad (75)$$

The last result comes about because $L_{n_K} > K$ and hence $S^{s_{n_K}}_{L_{n_K}}(|x_{n_K}|, \theta_{n_K \pm}) \gtrsim \sqrt{K}$. So one concludes that the modulus of the relative error $\varepsilon_{GS}(K)$ will satisfy

$$|\varepsilon_{GS}(K)| = \frac{|\delta(K, N, x, \theta, P)|}{|S_N(x, \theta)|} \approx \frac{|\delta(K, N, x, \theta, P)|}{\left|S^{s_1}_{L_1}(|x_1|, \theta_{1\pm})\right|} < \frac{3.41 \max(|\varepsilon_{n_K - m}|)}{\sqrt{K}}. \tag{76}$$

Each iteration error $\varepsilon_{n_K - m}$ is given by (61), and the $max(|\varepsilon_{n_K - m}|)$ will never exceed the maximum modulus $|R_3(x = 1/2, y = 1/2)| \approx 15/\pi^4$ if a cut-off value of $P = 3$ is employed (cf. discussion after eq. 61). Hence one concludes that the final relative error cannot exceed and indeed must be very much less than

$$|\varepsilon_{GS}(K)| \ll \frac{3.41}{\sqrt{K}} \times \frac{15}{\pi^4} \approx \frac{1}{2\sqrt{K}}. \tag{77}$$

All the relative errors for the calculations shown in Table 3.4 with $K = 30$ are orders of magnitude smaller than this rather crude upper bound. What is clear is that the size of both the actual $\delta(K, N, x, \theta, P)$ and relative $\varepsilon_{GS}(K)$ errors are regulated by the choice of $K$.

The case when $x \equiv x_1$ is a 'non-generic' irrational is only slightly different. In such an instance the algorithm terminates immediately the small quadratic parameter $x_n = x_{n_K} < 1/L_{n_K}$ associated with the large partial quotient $|a_{n_K}|$ is encountered. However, because $x_{n_K}$ is so small, the first source of error from (61) with $P = 3$ is given by $|\varepsilon_0| \sim O(x^3_{n_K}) \lesssim O(L^{-3}_{n_K})$, which is completely negligible compared to what comes next. For example the error in the first intermediate sum computed in Table 3.4D† is beyond the 9 significant digits shown in the display, whilst in Table 3.4E the relative error only $\sim 10^{-10}$. These relatively high precision initial estimates mean that the main source of error in the quadratic Gaussian sum arises from the error $|\varepsilon_1|$ in the *second* intermediate sum, percolating through the rest of the calculation. An upper bound for this error can be found in an analogous manner to the upper bound established in (74-77). The magnitude of the error satisfies

$$|\delta(K, N, x, \theta, P)| < \frac{3.41 |\varepsilon_1|}{\sqrt{|x_{n_K - 2} \dots x_2 x_1|}} \approx 3.41 |\varepsilon_1| \left| \frac{S^{s_1}_{L_1}(|x_1|, \theta_{1\pm})}{S^{s_{n_K}}_{L_{n_K - 1}}(|x_{n_K - 1}|, \theta_{n_K - 1 \pm})} \right|, \tag{78}$$

whilst the relative error is less than

$$|\varepsilon_{GS}(K)| = \frac{|\delta(K, N, x, \theta, P)|}{|S_N(x, \theta)|} \approx \frac{|\delta(K, N, x, \theta, P)|}{\left|S^{s_1}_{L_1}(|x_1|, \theta_{1\pm})\right|} < \frac{1}{2\sqrt{L_{n_K - 1}}} < \frac{1}{2\sqrt{L_{n_K}}} < \frac{1}{2\sqrt{K}}. \tag{79}$$

So the actual and relative errors in the computation of a general quadratic Gaussian sum $S_N(x, \theta)$ are *always smaller* when $x \equiv x_1$ is a 'non-generic' irrational as opposed to a 'generic' irrational, by at least a factor of $\sqrt{L_{n_K - 1}/L_{n_K}}$.

## 3.6 The choice of termination constant K

The most appropriate choice of $K$ comes down to a trade off between greater computational accuracy (large $K$) and computational speed (small $K$). However, if one is motivated purely to minimize the computational speed of the algorithm, one can obtain a useful estimate for a lower bound on $K$ that can be employed for all $S_N(x,\theta)$. Here $N$ is any arbitrarily large integer and $x \equiv x_1$ is any irrational ('generic' or 'non-generic'). Computing $S_N(x,\theta)$ directly from (53), requires 4 basic multiplication and addition operations to find each phase $\pi k[kx + 2\theta]$ (the value of $2\theta$ would be computed once and then stored) and a sine and cosine evaluation to find each term. There are $N$ terms and these must all be added together, which required a further $2N$ (complex) addition operations. So this gives a total of $\approx N\{2ops_{s;c} + 6ops_{ao}\}$ operations, where $ops_{s;c}$ is the operational count to compute the sine/cosine function for an arbitrary phase and $ops_{ao}$ is the operational count for any of the basic arithmetic operations $+;-;\times;\div$. However, for any computing platform the operational count for a single basic arithmetic operation is *at least* (or should be!) fifty times faster than for a single intrinsic sine/cosine evaluation $\Rightarrow ops_{s;c} \leq 50ops_{ao}$. Hence to compute $S_N(x,\theta)$ directly requires an operational count $\leq 2.12Nops_{s;c}$.

When applying *algorithm QGS*, the main computational problem concerns the evaluation of the appropriate $C(N,x,\theta,P)$ term (64) at each iteration step. This term contains three error function calls within it, which are key to assessing its computational speed. Now each of these calls is of the form $\text{erf}(\omega \times real)$ (since $\text{erfc}(\text{sgn}(-\varepsilon)z_M^-/\sqrt{2}) = 1 - \text{erf}(\omega\sqrt{\pi/x}\,|\varepsilon|)$) and hence $\arg(z) = \pm\pi/4$ for every $\text{erf}(z)$ evaluation. This being the case it is relatively simple to construct a purpose built erf routine which takes advantage of inherent symmetry properties that apply when $|\arg(z)| = \pi/4$. Since $|z|$ is not restricted, the routine would have to employ both the standard power series [32 eq. 7.6.1] and asymptotic expansions [32 eq. 7.12.1] valid for small and large $|z|$ respectively, as well as adopting rapidly converging power series local to the range $3/4 < |z| < 9/4$, when neither of the standard representations are particularly efficient. All these series methods involve the evaluation of a $2ops_{s;c}$ $e^{iz^2}$ term, plus the equivalent of a further $ops_{s;c}$ operation to compute the various power/asymptotic series that arise. Such a purpose built routine can easily be constructed so that the relative error lies below 0.5% for all $|z|$. However, most modern mathematical packages (such as MATLAB or MAPLE) contain an intrinsic routine built into their software designed to compute $\text{erf}(z)$ for arbitrary $\arg(z)$. These routines are somewhat slower, but much more precise, generally requiring the equivalent of $9ops_{s;c}$ to achieve an accuracy to within 10 significant figures. The operational count for remainder of the $C(N,x,\theta,P)$ term is easier to estimate. The short sum over $r$ takes about $P\{4ops_{s;c} + 40ops_{ao}\}$ depending upon whether exact expressions (cf. definitions 59-60)

$$c_0(x) = \cot(\pi x) - \frac{1}{\pi x}, \quad c_r(x) = \frac{1}{(2r)!}\frac{d^{2r}c_0(x)}{d(\pi x)^{2r}}, \quad r \geq 1$$

$$\Rightarrow c_1(x) = \frac{\cot(\pi x)}{\sin^2(\pi x)} - \frac{1}{(\pi x)^3} \text{ etc,} \qquad (80)$$

or power series expansions (the exact expressions very prone to round-off error if $|\pi x| < 0.2$)

$$c_0(x) = -\frac{\pi x}{3}\left(1 + \frac{(\pi x)^2}{15} + \frac{2(\pi x)^4}{315} + \cdots\right), c_1(x) = -\frac{\pi x}{3}\left(\frac{1}{5} + \frac{4(\pi x)^2}{63} + \frac{(\pi x)^4}{75} + \cdots\right), \text{etc } (81)$$

are used to compute the coefficients. The remaining computation of $C(N, x, \theta, P)$ and its amalgamation into iteration step (63) requires about $\{10ops_{S;c} + 80ops_{ao}\}$. So with $P = 3$ the total operational count comes to approximately $\{22ops_{S;c} + 200ops_{ao}\} \leq 26ops_{S;c}$ plus the three necessary erf(z) evaluations. Hence one obtains an estimate between $35 - 53ops_{S;c}$ for the evaluation of $C(N, x, \theta, P)$, depending upon whether a purpose built or an intrinsic routine is used to calculate erf(z).

Suppose the algorithm were to be employed with termination constant $K = 1$. Then the number of iterations required to truncate the initial sum size $N$ down to a value $L_{n_1+1} < 1$ would vary from an average value $\langle n_1 \rangle = 0.59 log(N)$ up to an maximum value $n_1 = log(N)/log(2) + 1 \approx 1.44 log(N)$. In what follows it is useful to define $\Omega^* = O(10)$ constant, representing the mean value of $y \times$ operational count of $C(N, x, \theta, P)$, where $y$ is the iteration count factor derived from (71). The absolute maximum of this quantity will be denoted by $\Omega^*_{max}$, and representative values (using an intrinsic erf routine) of $\Omega^* = 0.59 \times 53$ and $\Omega^*_{max} = 1.44 \times 53$ will be adopted for calculation purposes. So in the worst case scenario the algorithm would require $\Omega^*_{max} log(N) ops_{S;c}$ to compute $S_N(x, \theta)$ compared to the $2.12 N ops_{S;c}$ needed for the direct term by term calculation. The value of $N$ at which these two operational counts become comparable provides a suitable means of choosing $K$. In the worst case scenario, solving $\Omega^*_{max} log(K) = 2.12 K$ gives a value of $K \approx 189$. However, if the average number of iterations is substituted, then $K$ would satisfy $\Omega^* log(K) = 2.12 K \Rightarrow K \approx 61$. So in order to maximise the speed of the algorithm, a prescribed value for $K \in [61, 189]$ would be appropriate. The round figure of $K = 100$ suggests itself, which is equivalent to a value $\Omega^* = 0.87 \times 53$. Values for $K \in [61, 189]$ are motivated purely by the desire to maximize computational speed of the algorithm. If more precise estimates of $S_N(x, \theta)$ are required, then $K$ would need to be made larger. The corresponding $O(\{2K + \Omega^* log(N/K)\} ops_{S;c})$ operational count would then be increased, slowing the algorithm down.

## 4. The Computational Complexity of the Hardy Function $Z(t)$ in the limit $t \to \infty$

At the conclusion of Section 2, equation (54) expresses $Z(t)$ as a sum of $O(\alpha_E/M_t) = O(t^{1/3})$ subsidiary pairs of quadratic Gaussian sums $S^*_{M_t^-}(x(pc), \theta^\pm(pc))$, each of (ever increasing) length $M_t^- = (M_t(pc) - 1)/2$ given by (40). The algorithm developed in Section 3 from the *Theorem [Paris]* demonstrates that a general quadratic Gaussian sum $S_N(x, \theta)$ of arbitrary length can be evaluated, using an operational count no larger than $O(\{2K + \Omega^*_{max} log(N/K)\}) ops_{S;c}$, to within a relative error $|\varepsilon_{GS}(K)| \ll 1/2\sqrt{K}$ for a suitable choice of $K = 10 - 100$. In combination, these two results allow one to reduce the computational complexity of $Z(t)$ significantly below the $O(\sqrt{t}) ops_{S;c}$ of (8) or its variants. This assertion is framed in the following *Theorems 2 & 3*, the proofs of which are the subject of this section.

## 4.1 Theorem 2

Let $K, K_0, L \in \mathbb{N}$ and set $K_0 = 50$. For any $t > (5\sqrt{\pi}K_0)^6 > 0$, define $a = \sqrt{8t/\pi}$, $L = \lfloor \log(t) \rfloor - 28$ and choose $\varepsilon_t$ as a fixed parameter, prescribed so that $\varepsilon_t \xrightarrow[t \to \infty]{} 0$ and $\varepsilon_t^{-1} = o(t^\mu)$ for any $\mu > 0$. Define $h_t = \log(\varepsilon_t)/\log(t)$ and choose $K$ to lie in the range $K_0 \leq K \leq MAX\{K_0, \ell(\varepsilon_t^2 t)^{1/6}\}$, where $\ell = \exp(55[4h_t - 1]/48)/\sqrt{2\pi} \leq 0.127$. Fix $b > 0$ ($b \ll a$), and let $\alpha \in (a, a+b)$ denote the odd integers in this range. Fix $\alpha_E^c = INT_E\left[(\varepsilon_t t^2)^{1/3}/\left(\sqrt{\pi}\log(t)\right)\right]$, and let $\alpha_E \in [a+b, \alpha_E^c]$ denote elements of a subset of the even (pivot) integers within the given interval. Define $ZP(t)$ by the following sum:

$$ZP(t) = 2\sqrt{2}\left\{\sum_{\substack{\alpha > a \\ odd}}^{a+b} \frac{\cos\left(\frac{t}{2}\{\log(pc) + \frac{1}{pc}\} + \frac{t}{2} + \frac{\pi}{8}\right)}{(\alpha^2 - a^2)^{1/4}} + \sum_{\substack{\alpha_E \geq a+b \\ step\ 4M_t^- + 4}}^{\alpha_E^c} \frac{\text{Re}\{e^{i\omega^+} S^*_{M_t^-}(x, \theta^+) + e^{i\omega^-} S^*_{M_t^-}(x, \theta^-)\}}{(\alpha_E^2 - a^2)^{1/4}}\right\}.$$

(82)

*The various parameters of the sub-sums that make up $ZP(t)$ are given by*

$$pc(\alpha \text{ or } \alpha_E) = 2(\alpha/a)^2\left[1 + (1 - (a/\alpha)^2)^{\frac{1}{2}}\right] - 1;\quad x = 1/(pc - 1);\quad \theta^\pm = x/2 \pm a/4\sqrt{pc};$$

$$\omega^\pm = \pi[\theta^\pm - x/4] - t(\log(pc) + 1/pc + 1)/2 - \pi/8;\quad M_t^- = (M_t - 1)/2 \geq 0;$$

$$M_t(\alpha_E) = \begin{cases} NINT_O\left[\left\{\frac{\varepsilon_t a}{\pi}\left(\frac{\alpha_E}{a} - 1\right)\right\}^{1/3}\right] \geq 1, & a+b < \alpha_E \leq Ya \\ NINT_O\left[\frac{(\varepsilon_t^2 t)^{1/6}}{\sqrt{2\pi}}\left(\frac{2\alpha_E}{a} - \frac{a}{2\alpha_E}\right)\right], & Ya < \alpha_E \leq \alpha_E^c \end{cases}$$

(83a, b)

*where $Y = 1 + O(1), \forall t$. Suppose the quadratic Gaussian sums $S^*_{M_t^-}(x, \theta^\pm)$ are estimated using algorithm QGS with termination constant $K$.*

*i) If $K_0 \leq K < \sqrt{3}[\varepsilon_t^4 t/\pi^4]^{1/12}/2$, the calculation of $ZP(t)$ requires only*

$$\left(\frac{t}{\varepsilon_t}\right)^{\frac{1}{3}} \log(t)\{O(\log(t)) + O(K)\}ops_{s;c}$$

*computational operations in the limit $t \to \infty$.*

*ii) Otherwise, if $\sqrt{3}[\varepsilon_t^4 t/\pi^4]^{1/12}/2 < K \ll \ell(\varepsilon_t^2 t)^{1/6}$, then the calculation of $ZP(t \to \infty)$ requires $\ll O(\sqrt{t})ops_{s;c}$ computational operations.*

## 4.2 The cut-off and operational conversions

Before embarking on the analysis necessary to prove *Theorem 2*, a couple of points are worth elucidating with regards to the implications for estimating the computational complexity of $Z(t)$ itself. Equation (54) gives an expression for the Hardy function of which $ZP(t)$ forms only a part, because the cut-off set in (82), namely $\alpha_E^c \approx (\varepsilon_t t^2)^{1/3}/\left(\sqrt{\pi}\log(t)\right)$ is very much less than the cut-off $N_\alpha \geq INT_O(t/\pi) + 2$ in (54). The idea behind this choice of $\alpha_E^c$ is to use (82) to estimate a portion of the terms that appear in the main sum of (8) and then compute the remainder directly. Structurally this will give rise to a second hybrid formula for $Z(t)$, analogous to (17). In that case the last $N \in \left[\frac{1}{2}(t/\pi)^{1/2}, N_t\right]$ terms of (8) were replaced by the first $\alpha \in \left[INT_O(a) + 2, INT_O\left(3\sqrt{t/\pi}\right)\right]$ odd terms of (11), reducing the total summands by 17.16%. But this proposed second hybrid formula is *more ambitious,* since $ZP(t)$ is representative of *all but* the first $N \in [1, t^{1/3} \log(t)/\sqrt{\pi}\varepsilon_t^{1/3}]$ terms of (8). (The size of the relative error, encapsulated by the term "representative" in the previous sentence, is the subject *Theorem 3*.) The implication of *Theorem 2,* is that the sum of the overwhelming majority of the $N \sim O(\sqrt{t})$ terms in the RSF can be computed, via $ZP(t)$, in something close to $O(t^{1/3})$ operations. This is what lies at the heart of reducing the computational complexity of the Hardy function.

The other point is a small technical issue concerning some conversions linking the standard sine/cosine operational count $ops_{s;c}$, with $ops_{RS}$ and $ops_{hybrid}$ the operational counts required to compute a typical RSF term and a typical principal term of the new series (11) respectively. A typical RSF term, including the amplitude, consists of *sqrt, log* and *cosine* evaluations, so $ops_{RS} \equiv ops_{s;c} + ops_{log} + ops_{sqrt}$. Estimates devised from the MAPLE mathematical package show that $ops_{log} \approx 0.62 ops_{s;c}$ and $ops_{sqrt} \approx 0.23 ops_{s;c}$, so $ops_{RS} \approx 1.85 ops_{s;c}$, which for simplicity will be rounded up to define a conversion relation $ops_{RS} \equiv 2 ops_{s;c}$. In [27 eq. 72] an optimal method for computing $ops_{hybrid} \equiv ops_{s;c} + ops_{log} + 2 ops_{sqrt}$ is presented. Based on the estimates above, this implies that $ops_{hybrid} = \Omega ops_{RS}$ where $\Omega \approx 1.12$. Large scale sample computations of the two respective sums in (17) suggest this estimate is a little too low, with $\Omega \approx 1.17 - 1.30$ depending upon the programming scheme employed. So for the purpose of analysing *Theorem 2,* let $ops_{hybrid} \equiv 2\Omega ops_{s;c}$ define a second conversion relation, with $\Omega$ lying inside the range given above.

*4.3 Proof of Theorem 2*
*4.3.1 Operational count of first sum making up $ZP(t)$.*

The first part of the proof is very straight forward. Assuming all the bounds concerning $t$ set out in *Theorem 2* are met, then define $b = INT_E(t^{1/4}) > 0$. This particular choice of $b$ is motivated by the desire that the number $\alpha_E > a + b$, which commences the second sum making up $ZP(t)$, is large enough to ensure that $M_t \geq 1$. As discussed at the start of section 2.5, this is a necessary pre-requisite which allows the principal terms of the new series (11) for $Z(t)$ to be collected together into the quadratic Gaussian sums of (54), in such a way that the relative error introduced is always $\leq O(\varepsilon_t)$. With this choice of $b$, one finds that at $\alpha = a + b$ the corresponding value of $M_t$, defined by (83a), will have reached

$$M_t = \left(\frac{\varepsilon_t}{\pi}\right)^{1/3} t^{1/12} + O\left(\frac{\varepsilon_t^{1/3}}{t^{1/24}}\right). \tag{84}$$

The bounds imposed on the definition of $\varepsilon_t$ means that the $t^{1/12}$ factor dominates (84) and $M_t \geq 1 \ \forall \ t \geq (5\sqrt{\pi}K_0)^6 = O(10^{15})$. Computation of the first sum in (82) is then equal to $b/2$, the number of typical principal terms, times $ops_{hybrid}$, the operational count of a typical term. Using the conversions discussed in 4.2, this is equivalent to $INT_E(t^{1/4})/2 \times 2\Omega ops_{s;c} = O(\Omega t^{1/4}) ops_{s;c}$ operations.

*4.3.2 The precise definition of $M_t$ for pivots $\alpha_E \in (a + b, Ya]$.*

The next part of the proof is to establish the computations needed to evaluate the second sum starting from the first pivot integer $\alpha_E > a + b$, where $pc(\alpha_E) \approx 1$, to some 'suitable' pivot $\alpha_E = Ya$ where $Y = 1 + O(1)$ constant (termed the 'step-up' parameter), at which point $pc(\alpha_E) = 1 + O(1)$. Specifically the value of $Y$ is deemed 'suitable' provided it is large enough so that *at some point* within the pivots $\alpha_E \in (a + b, Ya]$, the length of the quadratic Gaussian sums $M_t^-$ exceeds the set value of $K$ in the range specified in *Theorem 2*. When $M_t^- > K$, substitution of the *algorithm QGS* for their computation becomes applicable.

Now equation (40) prescribes an upper bound on the collection size $M_t$ which ensures the process of combining the principal terms of the new series (11) into quadratic Gaussian sums always yields relative errors $\leq O(\varepsilon_t)$. However, the upper bound (40) is a little awkward to deal with as the $(pc - 1)$ factor causes $M_t$ to move rapidly through a number of different scale regimes as $pc$ increases from unity. These rapid changes in the size of $M_t$ must be carefully be accounted for, in order to estimate the operational count. Suppose $\alpha_E = Ya$ and $1 < Y \leq 3/\sqrt{8}$, so that $pc = 1 + x$ with $x = \sqrt{8(Y - 1)} + O(Y - 1) \leq 1$. If $M_t$ were to be defined using (40) then

$$M_t \approx \sqrt{\frac{2}{\pi}} \varepsilon_t^{1/3} x^{2/3} t^{1/6} \left(1 - \frac{x}{2} + \frac{11x^2}{24} \ldots\right) = 2\left(\frac{\varepsilon_t a}{\pi}\left(\frac{\alpha_E}{a} - 1\right)\right)^{\frac{1}{3}} (1 - small\ term). \qquad (85)$$

It is tricky to keep track of rapidly changing scale of the *small term* in (85), so ideally one would like to forget about it by setting it to zero. However, because the *small term* is negative, the resulting value of $M_t$ would then exceed the upper bound. To avoid this problem, $M_t$ is defined to equal $(\varepsilon_t(\alpha_E - a)/\pi)^{1/3}$ for $\alpha_E \in [a + b, Ya]$ in *Theorem 2*, i.e. approximately one half of the upper bound (40, 85). Fixing $M_t$ in this way, also ensures it never exceeds its upper bound even if $pc$ lies considerably above unity. For example, if $Y = 7/4 \Rightarrow pc \approx 10.15$, then (83a) gives $M_t \approx 0.72\varepsilon_t^{1/3}t^{1/6} < 1.59\varepsilon_t^{1/3}t^{1/6}$, the corresponding upper bound value (40). This means that the length of the quadratic Gaussian sums is continually constrained to keep the relative errors $\leq O(\varepsilon_t)$ for all $\alpha_E \in [a + b, Ya]$.

In practical calculations (see Section 5), one prescribes a 'suitable' fixed value of $Y$ and then steps up the size of the $M_t$ values (utilising the second definition given in *Theorem 2*), which are computationally more efficient when $\alpha_E$ is no longer close to $a$. However, the 'suitability' of $Y$ is difficult to pin down precisely, since it is determined by the respective values of $t, K$ and $\varepsilon_t$ in *Theorem 2*. So instead for the purposes of proving the theorem, the parameter $L = \lfloor log(t) \rfloor - 28$ is chosen and utilised in such a way as to guarantee $Y$ is 'suitable' for all possible choices of $t, K$ and $\varepsilon_t$, whilst at the same time ensuring the 'stepping up' *always occurs* before $Y$ can grow beyond the $1 + O(1)$ scale regime. This is the reason behind the appearance of the $-28$ factor in $L$, as will become clear shortly.

*4.3.3 Operational count of second sum of $ZP(t)$ for pivots $\alpha_E \in (a+b, Ya]$.*

Define the constant $X = exp[1/12 - h_t/3] > 1.08690 ...$, with $h_t = log(\varepsilon_t)/log(t)$. Let the interval $\alpha_E \in [a+b, Ya]$ be divided up into a series of slowly increasing subintervals or pivot 'Blocks', numbered $p = 1, 2, .., L$. Each Block consists of $\lfloor X^p t^{1/4} \rfloor$ different pivots $\alpha_E$, and associated with each Block is a fixed step parameter $M_{t,p}$. This prescribes the lengths of the Gaussian sums $S^*_{M_t^-}(x(\alpha_E), \theta^\pm(\alpha_E))$ defined by (53 & 83), and means that within the Block each pivot is separated from its two nearest neighbours by a factor $2M_{t,p} + 2 = 4M_{t,p}^- + 4$. Let $M_{t,1} = 1$ and $\alpha_{E,0} = INT_O(a+b) + 1$. The choice of $M_{t,1} = 1$ is conservative, but serves to smooth out the transition from the first to the second sum of $ZP(t)$. Starting at $\alpha_{E,0} + 2$, compute the second sum (82) over Block 1, which consists of the $\lfloor Xt^{1/4} \rfloor$ pivots $\alpha_E = \alpha_{E,0} + 1 \times 2, \alpha_{E,0} + 3 \times 2, \alpha_{E,0} + 5 \times 2, ..., \alpha_{E,0} + (2\lfloor Xt^{1/4} \rfloor - 1) \times 2$. The first even integer $\alpha_{E,1}$ bounding the pivots residing in Block 1 is given by

$$\alpha_{E,1} = \alpha_{E,0} + 2\lfloor Xt^{1/4} \rfloor (M_{t,1} + 1) = \alpha_{E,0} + 4\lfloor Xt^{1/4} \rfloor. \tag{86}$$

Now define the next step parameter for Block 2 to be $M_{t,2} = M_t(\alpha_{E,1})$. Block 2 will consist of the $\lfloor X^2 t^{1/4} \rfloor$ pivots $\alpha_E = \alpha_{E,1} + 1 \times (M_{t,2} + 1), ..., \alpha_{E,1} + (2\lfloor X^2 t^{1/4} \rfloor - 1) \times (M_{t,2} + 1)$. The lowest even integer $\alpha_{E,2}$ bounding the pivots residing in Block 2 is given by

$$\alpha_{E,2} = \alpha_{E,1} + 2\lfloor X^2 t^{1/4} \rfloor (M_{t,2} + 1). \tag{87}$$

The step parameter for Block 3 is then $M_{t,3} = M_t(\alpha_{E,2})$ and the process repeated for all the subsequent Blocks. So for Block $p$ one has

$$\begin{cases} \text{Step parameter } M_{t,p} = M_t(\alpha_{E,p-1}), \text{ where} \\ \alpha_{E,p-1} = INT_O(a+b) + 1 + 2 \sum_{i=1}^{p-1} \lfloor X^i t^{1/4} \rfloor (M_{t,i} + 1), \\ \text{with pivots } \alpha_E = \alpha_{E,p-1} + 1 \times (M_{t,p} + 1), ..., \alpha_{E,p-1} + \left(2\lfloor X^p t^{\frac{1}{4}} \rfloor - 1\right) \times (M_{t,p} + 1). \end{cases} \tag{88}$$

Next one needs to estimate some lower bounds on $M_{t,p}$ as $p = 1, 2, .., L$. Since $\alpha_{E,0} \approx a + b \approx a + t^{1/4}$, substituting (86) into (83a) gives

$$M_{t,2} = M_t(\alpha_{E,1}) \approx \left(\frac{\varepsilon_t}{\pi}\right)^{1/3} 4^{1/3} (X)^{1/3} t^{1/12} \left[1 + \frac{1}{4X}\right]^{1/3} = (4Xg)^{1/3} \left(\frac{\varepsilon_t}{\pi}\right)^{1/3} t^{1/12}, \tag{89}$$

where $g = [1 + 1/4X]$. Substituting (89) into (87) gives

$$\alpha_{E,2} \approx a + t^{\frac{1}{4}} + 4\left[Xt^{\frac{1}{4}}\right] + 2\left[X^2 t^{\frac{1}{4}}\right] \times (2\mathcal{g})^{\frac{1}{3}}(X)^{\frac{1}{3}}\left[\left(\frac{2\varepsilon_t}{\pi}\right)^{\frac{1}{3}} t^{\frac{1}{12}}\right]^1,$$

$$= a + 2(2\mathcal{g})^{\frac{1}{3}}\left(\frac{2\varepsilon_t}{\pi}\right)^{\frac{1}{3}}(t^{1/4})^{1+1/3}(X)^{\frac{1}{3}}(X^2)^1 + O\left(t^{\frac{1}{4}}\right). \tag{90}$$

Dropping the $O(t^{1/4})$ terms underestimates $\alpha_{E,2}$. So on substituting (90) into (83a), one obtains a lower bound for $M_{t,3} = M_t(\alpha_{E,2}) > (2\mathcal{g})^{1/9}(X)^{1/9}(X^2)^{1/3}\left[\left(\frac{2\varepsilon_t}{\pi}\right)^{1/3} t^{1/12}\right]^{1+1/3}$. This leads to an estimate $\alpha_{E,3} > a + 2(2\mathcal{g})^{\frac{1}{9}}\left(\frac{2\varepsilon_t}{\pi}\right)^{1/3+1/9}(t^{1/4})^{1+1/3+1/9}(X)^{1/9}(X^2)^{1/3}(X^3)^1 + O(t^{1/3})$, which (to first order) can be used to find lower bound for $M_{t,4}$ etc. The pattern for $p = 1, 2, \dots, L-1$ that emerges is given by

$$M_{t,p+1} > (2\mathcal{g})^{\frac{1}{3^p}}\left[\left(\frac{2\varepsilon_t}{\pi}\right)^{\frac{1}{3}} t^{\frac{1}{12}}\right]^{1+\frac{1}{3}+\frac{1}{9}+\dots+\frac{1}{3^{p-1}}} \prod_{i=1}^{p}(X^i)^{1/3^{p+1-i}}. \tag{91}$$

This bound has important implications. From the definition of $X$ prescribed at the start of this section, one can compute

$$\log\left\{\prod_{i=1}^{p}(X^i)^{1/3^{p+1-i}}\right\} = \frac{\log(X)}{3^{p+1}}\sum_{i=1}^{p} i \times 3^i = \frac{\log(X)}{3^{p+1}} \times \frac{3^{p+1}(2p-1)+3}{4},$$

$$\Rightarrow \prod_{i=1}^{p}(X^i)^{1/3^{p+1-i}} = X^{p/2-1/4+1/(4\times 3^p)} = \exp\left\{\frac{2p-1+3^{-p}}{48}\left[1 - \frac{4\log(\varepsilon_t)}{\log(t)}\right]\right\},$$

$$\Rightarrow M_{t,p+1} > (2\mathcal{g})^{\frac{1}{3^p}}\left[\left(\frac{2\varepsilon_t}{\pi}\right)^{\frac{1}{2}} t^{\frac{1}{8}}\right]^{1-\frac{1}{3^p}} X^{p/2-1/4+1/(4\times 3^p)}. \tag{92}$$

Now consider the final Block $L = \lfloor \log(t) \rfloor - 28$ which has step parameter $M_{t,L}$. If $p+1 = L$ in (92), the factor $3^{-L+1}$ factor is negligible in comparison with the other terms (e.g. if $t > (5\sqrt{\pi}K_0)^6 \sim O(10^{16})$, then $L \geq 9, \Rightarrow 3^{-L+1} < 0.00015$ and $(t^{1/8})^{-3^{-L+1}} \approx 1$). Hence

$$M_{t,L} > \frac{\sqrt{2}\varepsilon_t^{1/3} t^{1/6}}{\sqrt{\pi} e^{55(1-4h_t)/48}} = 2\ell\varepsilon_t^{1/3} t^{1/6} \xrightarrow[t\to\infty]{} 0.254\varepsilon_t^{1/3} t^{1/6}. \tag{93}$$

*Theorem 2* specifies $\varepsilon_t^{-1} = o(t^\mu)$ for any $\mu > 0$, and so the factor $\exp\left\{\frac{55h_t}{12}\right\} \to 1$, although potentially rather slowly. *Theorem 2* also specifies that $K < \ell\varepsilon_t^{1/3} t^{1/6}$, so the implication of (93) is that there *must exist* a $p^* \in [1, L-1]$ such that $M_{t,p^*}^- < K < M_{t,p^*+1}^-$. Hence for those pivots $\alpha_E$ lying inside Blocks $p > p^*$, computation of their respective quadratic Gaussian sums can be achieved using *algorithm QGS*, each to within a relative error $< 1/2\sqrt{K}$, by means of a significantly reduced operational count.

For example, suppose $t \geq O((2K_0)^9)$ is chosen to be somewhat larger than the minimum value $(5\sqrt{\pi}K_0)^6$ specified in *Theorem 2*. Then substituting $p = 2$ into (92) gives

$$M_{t,3} > \left(\frac{5}{2}\right)^{1/9} \left(\frac{2\varepsilon_t}{\pi}\right)^{4/9} t^{1/9} \exp\{7(1-4h_t)/108\} \approx 0.97 e^{-7h_t/27} \varepsilon_t^{4/9} t^{1/9}. \tag{94}$$

Hence one could select $K$ so that it satisfies $K_0 \leq K < e^{-7h_t/27}(\varepsilon_t^4 t)^{1/9}/2 \approx M_{t,3}^-$, in which case the value of $p^* = 2$. For still larger values of $t \geq O((2K_0)^{12})$, it would be possible to select $K$ so that it lies between $K_0 < K < M_{t,2}^-$ (from the bound on $M_{t,2}$ given by eq. 89), when $p^*$ would fall to unity. In summary, $p^* \leq L$ for any $K$ lying in the range specified in *Theorem 2*, and if $K$ is restricted to values near the lower end of its prescribed range, e.g. either $K < M_{t,2}^-$ or $M_{t,2}^- < K < M_{t,3}^-$, then $p^*$ approaches 1 or 2 respectively as $t \to \infty$.

If the lower bound (92) is substituted into (88), by the end of Block $L$, one has reached a pivot value

$$\alpha_{E,L} > a + b + 1 + 4Xt^{1/4} + 2\sum_{p=2}^{L} \lfloor X^i t^{1/4} \rfloor (2\mathcal{g})^{\frac{1}{3^{p-1}}} \left[\left(\frac{2\varepsilon_t}{\pi}\right)^{\frac{1}{2}} t^{\frac{1}{8}}\right]^{1+\frac{1}{3^{p-1}}} X^{(p-1)/2 - 1/4 + 1/(4\times 3^{p-1})},$$

$$\approx a + 1 + 4X\mathcal{g}t^{1/4} + 2X^{-3/4}\left(\frac{2\varepsilon_t}{\pi}\right)^{\frac{1}{2}} t^{3/8} \sum_{p=2}^{L}\left[2\mathcal{g}\left(\frac{2\varepsilon_t}{\pi}\right)^{\frac{1}{2}} X^{1/4} t^{1/8}\right]^{\frac{1}{3^{p-1}}} \left(X^{3/2}\right)^p. \tag{95}$$

Now the vast majority of the $p$ values in (95) lie in the range $\log(\log(t)/8) < p < L$, which means that the factor in square brackets $[O(\sqrt{\varepsilon_t})t^{1/8}]^{-1/3^{-p+1}} \approx 1$. For the smallest $p$ values this does not hold, but these terms themselves are completely dominated by those associated with large $p$ values close to $L$. Hence an excellent lower bound estimate on the sum in (95) can be found from

$$\sum_{p=2}^{L}\left(X^{3/2}\right)^p = \frac{X^{3/2}\left(X^{3L/2} - X^{3/2}\right)}{(X^{3/2} - 1)} \approx \frac{X^{3/2}}{(X^{3/2}-1)} X^{3\log(t)/2 - 42} \approx \frac{e^{14h_t - 7/2} X^{3/2} t^{1/8}}{\sqrt{\varepsilon_t}(X^{3/2} - 1)}. \tag{96}$$

So by the end of Block $L$ one *must have reached at least*

$$\alpha_{E,L} > a + 1 + 4X\mathcal{g}t^{1/4} + \frac{e^{14h_t - 7/2} X^{3/4} a}{(X^{3/2} - 1)} \approx Ya,$$

$$\text{with}\quad Y = 1 + \frac{e^{14h_t - 7/2} X^{3/4}}{(X^{3/2} - 1)} \xrightarrow[t\to\infty]{} 1.2414\dots. \tag{97}$$

The factor $-28$ in the definition for $L = \lfloor \log(t) \rfloor - 28$ introduces the $e^{-7/2}$ term into this lower bound estimate for $Y$, ensuring that the stepping up of the $M_t$ values occurs at a reasonable small $Y = 1 + O(1)$ value. Numerical calculations (for $t$ up to $10^{1500}$ and appropriate choices of $\varepsilon_t$) show that the value $Y = 1.2401$ is usually achieved by the time $p \approx L - \log(\log(t))$ and reaches $Y \approx 1.6$ at $p = L$. Since (97) gives only a lower bound estimate for $Y$, this higher value is to be expected. But it too remains well within the $1 + O(1)$ regime.

The operational cost of all these calculations is as follows. For Block 1, $M_{t,1} = 1 \Rightarrow M^-_{t,1} = 0$, so one effectively has two $ops_{hybrid}$ terms to compute for each of the $\lfloor Xt^{1/4} \rfloor$ pivots. At every pivot within Blocks $p = 2, 3, \ldots, p^*$ one has to compute two quadratic Gaussian sums of length $M^-_{t,p}$ longhand. That is equivalent to $2M^-_{t,p} = (M_{t,p} - 1) ops_{hybrid}$ terms for each of the $\lfloor X^p t^{1/4} \rfloor$ pivots. For Blocks $p = p^* + 1, \ldots, L$ computation of the two quadratic Gaussian sums of length $M^-_{t,p}$ can be achieved to within a relative error $1/2\sqrt{K}$ using *algorithm QGS*. As discussed in Section 3.6 this will take about $2\{2K + \Omega^* \log((M_{t,p} - 1)/2K)\} ops_{S;c}$; certainly no more than the maximum possible number of $2\{2K + \Omega^*_{max} \log((M_{t,p} - 1)/2K)\} ops_{S;c}$. One must also add in a further $2ops_{hybrid}$ to account for the calculation of the $e^{-i\omega^\pm}/(\alpha_E^2 - a^2)^{\frac{1}{4}}$ amplitude factors that appear in (82). So across all the Blocks the average number of computations required will amount to

$$2\lfloor Xt^{1/4} \rfloor ops_{hybrid} + \sum_{p=2}^{p^*} \lfloor X^p t^{1/4} \rfloor \times (M_{t,p} - 1) ops_{hybrid} + 2 \sum_{p=p^*+1}^{L} \lfloor X^p t^{1/4} \rfloor$$
$$\times [\{2K + \Omega^* \log((M_{t,p} - 1)/2K)\} ops_{S;c} + ops_{hybrid}],$$

$$\equiv \left[ 4\Omega \lfloor Xt^{1/4} \rfloor + 2\Omega \sum_{p=2}^{p^*} \lfloor X^p t^{1/4} \rfloor (M_{t,p} - 1) + (4\Omega + 4K - 2\Omega^* \log(2K)) \sum_{p=p^*+1}^{L} \lfloor X^p t^{1/4} \rfloor \right.$$
$$\left. + 2\Omega^* \sum_{p=p^*+1}^{L} \lfloor X^p t^{1/4} \rfloor \log(M_{t,p} - 1) \right] ops_{S;c}, \tag{98}$$

using the $ops_{hybrid} \equiv 2\Omega ops_{S;c}$ conversion factor. Now it is just a matter of estimating the three sums in (98).

The second sum is trivial to estimate since

$$\sum_{p=p^*+1}^{L} X^p = \frac{X^{L+1}}{(X-1)} \left(1 - \frac{X^{p^*}}{X^L}\right) \approx \frac{e^{9h_t - 9/4} t^{1/12}}{\varepsilon_t^{1/3} (X-1)} B_0, \tag{99}$$

where $B_0 = (1 - X^{p^*}/X^L)$. The third sum is a bit more difficult, because it requires an upper bound estimate for $M_{t,p}$, not yet established. One can achieve this by demonstrating that the value of $M_{t,p}$ will certainly be bounded above by the lower bound (92) established for $M_{t,p+1}$. To show this, consider how the difference $(\alpha_{E,p} - a)$ grows with $p$. From (88) one can see that to *second* order this difference is given by $2\lfloor X^{p-1} t^{1/4} \rfloor (XM_{t,p} + M_{t,p-1} + 2)$. Hence, to second order (q.v. eqs. 91-92)

$$M_{t,p+1} \geq (2g)^{\frac{1}{3^p}} \left[\left(\frac{2\varepsilon_t}{\pi}\right)^{\frac{1}{2}} t^{\frac{1}{8}}\right]^{1-\frac{1}{3^p}} X^{p/2-1/4+1/(4\times 3^p)} \left[1 + \frac{d_p M_{t,p-1}}{XM_{t,p}}\right]^{1/3}, \quad p = 1 \ldots L-1 \quad (100)$$

where $M_{t,0} = 0, M_{t,1} = 1, d_2 = (2g)^{2/3}$ and $d_p = 1$ otherwise. This means that for $p \geq 2$ the ratio

$$\frac{M_{t,p}}{M_{t,p+1}} \approx X^{-(1+3^{-p})/2} \left[\frac{2\pi g^2}{\varepsilon_t t^{1/4}}\right]^{1/3^p} \frac{\left(1 + M_{t,p-2}/XM_{t,p-1}\right)^{1/3}}{\left(1 + M_{t,p-1}/XM_{t,p}\right)^{1/3}}. \quad (101)$$

The ratio (101) is initially very small and then gradually increases, approaching, but never exceeding the value $1/\sqrt{X}$ as $p \to \infty$. Hence $M_{t,p} \leq M_{t,p+1}/\sqrt{1.08690} < M_{t,p+1} \, \forall p$, as postulated. Using this result, one can now estimate an upper bound on the third sum in (98) as follows:

$$\sum_{p=p^*+1}^{L} X^p \log(M_{t,p} - 1) < \sum_{p=p^*+1}^{L} X^p \log(\text{lower bound on } M_{t,p+1})$$

$$\approx \log\left(\left(\frac{2\varepsilon_t}{\pi}\right)^{1/2} t^{1/8}\right) \sum_{p=p^*+1}^{L} X^p + \log\left(g\sqrt{2\pi/\varepsilon_t t^{1/4}}\right) \sum_{p=p^*+1}^{L} \left(\frac{X}{3}\right)^p$$

$$+ \frac{\log(X)}{4} \left\{ 2 \sum_{p=p^*+1}^{L} pX^p - \sum_{p=p^*+1}^{L} X^p + \sum_{p=p^*+1}^{L} \left(\frac{X}{3}\right)^p \right\}. \quad (102)$$

Since $X < 3$, the sums involving $(X/3)^p$ are negligibly small. The sums involving just $X^p$ are given by (99), which leaves

$$\sum_{p=p^*+1}^{L} pX^p = \frac{LX^{L+1}}{(X-1)} \left(1 - \frac{p^* X^{p^*}}{LX^L} - \frac{1}{L(X-1)} \left\{1 - \frac{X^{p^*}}{X^L}\right\}\right)$$

$$\approx \frac{\{\log(t) - 28\} e^{9h_t - 9/4} t^{1/12}}{\varepsilon_t^{1/3}(X-1)} B_1, \quad (103)$$

where $B_1 = \left(1 - p^* X^{p^*}/LX^L - B_0/L(X-1)\right)$. Hence

$$\sum_{p=p^*+1}^{L} X^p \log(M_{t,p} - 1)$$

$$< \frac{e^{9h_t - 9/4} t^{\frac{1}{12}}}{2\varepsilon_t^{\frac{1}{3}}(X-1)} \left[\log(X)\left[\{\log(t) - 28\}B_1 - \frac{B_0}{2}\right] + \left\{\frac{\log(t)}{4} + \log(2\varepsilon_t/\pi)\right\} B_0\right]. \quad (104)$$

Finally, one must estimate the first sum in (98), running from $= 2, \ldots, p^*$. Since $p^* \ll L$ for large $t$ (indeed $p^* = 1$ for $t \geq O((2K_0)^{12})$, in which case this sum will vanish), procuring an upper bound by simply replacing $M_{t,p}$ by $M_{t,p+1}$ leads to a rather crude approximation, which will vastly overestimate (98) if $p^*$ is small, e.g. if $p^* < \log(\log(t)/8)$. A much better estimate can be obtained by first modifying expression (100) for $M_{t,p+1}$ to include corrections of *all* order of magnitude and then employing $M_{t,p}/M_{t,p+1} < 1/\sqrt{X}$. From (88), one obtains for $p \geq 0$

$$M_{t,p+1} \approx (2g)^{\frac{1}{3^p}} \left[\left(\frac{2\varepsilon_t}{\pi}\right)^{\frac{1}{2}} t^{\frac{1}{8}}\right]^{1-\frac{1}{3^p}} X^{p/2 - 1/4 + 1/(4 \times 3^p)} \left[1 + \sum_{j=1}^{p-1} \frac{d_{p+1-j} M_{t,p-j}}{X^j M_{t,p}}\right]^{1/3},$$

$$< (2g)^{\frac{1}{3^p}} \left[\left(\frac{2\varepsilon_t}{\pi}\right)^{\frac{1}{2}} t^{\frac{1}{8}}\right]^{1-\frac{1}{3^p}} X^{p/2 - 1/4 + 1/(4 \times 3^p)} \left[\sum_{j=0}^{p-1} \frac{1}{X^{\frac{3j}{2}}}\right]^{\frac{1}{3}},$$

$$< (2g)^{\frac{1}{3^p}} \left[\left(\frac{2\varepsilon_t}{\pi}\right)^{\frac{1}{2}} t^{\frac{1}{8}}\right]^{1-\frac{1}{3^p}} X^{p/2 - 1/4 + 1/(4 \times 3^p)} \left[\frac{X^{\frac{3}{2}}}{X^{\frac{3}{2}} - 1}\right]^{\frac{1}{3}}. \quad (105)$$

Hence the first sum in (98) is bounded by

$$\sum_{p=2}^{p^*} \lfloor X^p \rfloor (M_{t,p} - 1) < \left[\frac{X^{\frac{3}{2}}}{X^{\frac{3}{2}} - 1}\right]^{\frac{1}{3}} \left[X^{-1/4} \left(\frac{2\varepsilon_t}{\pi}\right)^{\frac{1}{2}} t^{\frac{1}{8}}\right]^{p^*} \sum_{p=2}^{p^*} X^{3p/2} \left[2gX^{\frac{1}{4}} \left(\frac{2\varepsilon_t}{\pi}\right)^{-\frac{1}{2}} t^{-\frac{1}{8}}\right]^{\frac{1}{3^p}} - \sum_{p=2}^{p^*} X^p,$$

$$< \left[\frac{X^{\frac{3}{2}}}{X^{\frac{3}{2}} - 1}\right]^{\frac{1}{3}} (2g)^{1/3^{p^*}} \left[X^{-1/4} \left(\frac{2\varepsilon_t}{\pi}\right)^{\frac{1}{2}} t^{\frac{1}{8}}\right]^{1-\frac{1}{3^{p^*}}} \sum_{p=2}^{p^*} X^{3p/2} - \sum_{p=2}^{p^*} X^p,$$

$$= (2g)^{1/3^{p^*}} \left[X^{-1/4} \left(\frac{2\varepsilon_t}{\pi}\right)^{\frac{1}{2}} t^{\frac{1}{8}}\right]^{1-\frac{1}{3^{p^*}}} \frac{X^2 (X^{3p^*/2} - X^{3/2})}{(X^{3/2} - 1)^{4/3}} - \frac{X(X^{p^*} - X)}{(X - 1)}. \quad (106)$$

Substituting bounds (99), (104) and (106) for the three sums into (98), one obtains an upper bound to the operational count necessary to compute the second sum of $ZP(t)$ for pivots $\alpha_E \in (a + b, Ya]$. One finds that this computation requires no more than

$$< 2\Omega t^{1/4} \left[ 2X + (2g)^{1/3p^*} \left[ X^{-1/4} \left(\frac{2\varepsilon_t}{\pi}\right)^{\frac{1}{2}} t^{\frac{1}{8}} \right]^{1-\frac{1}{3p^*}} \frac{X^2(X^{3p^*/2} - X^{3/2})}{(X^{3/2} - 1)^{4/3}} - \frac{X(X^{p^*} - X)}{(X - 1)} \right] ops_{S;c}$$

$$+ \frac{\Omega^* e^{9h_t - 9/4} t^{1/3} log(t)}{\varepsilon_t^{1/3}(X - 1)} \left\{ \left\{ B_1 log(X) + \frac{B_0}{4} \right\} + \frac{1}{log(t)} \left[ B_0 log\left(\frac{\varepsilon_t}{2\pi K^2 \sqrt{X}}\right) - 28 B_1 log(X) \right] \right.$$
$$\left. + \frac{4 B_0 (\Omega + K)}{\Omega^*} \right\} ops_{S;c}. \qquad (107)$$

*4.3.4 The precise definition of $M_t$ for pivots $\alpha_E \in (Ya, \alpha_E^c]$.*

In section 4.3.2, the precise definition (83a) for the step parameter $M_t$ based around (40) and (85), was motivated by considerations of the first pivot regime in which both $Y$ and $pc$ are close to unity. From the estimates leading up to (97) in the previous section, all the pivots in this second regime satisfy $\alpha_E \geq Ya$ where $Y > 1.2401$. This equates to a value of $pc > 3.894$, which is sufficiently above unity to approximate the upper bound (40) for $M_t$ by the following:

$$M_t(pc) \approx \left(\frac{\varepsilon_t^2 t}{2\pi^3}\right)^{1/6} \frac{(pc^2 - 1)^{2/3}}{(pc^5)^{1/6}} \approx \left(\frac{\varepsilon_t^2 t}{2\pi^3}\right)^{1/6} \left[\sqrt{pc} - \frac{2}{3 pc^{3/2}} + O(pc^{-5/2})\right]. \qquad (108)$$

The simpler definition (83b) proposed for $M_t(\alpha_E)$ in *Theorem 2* is asymptotically similar to (108), since

$$\frac{(\varepsilon_t^2 t)^{1/6}}{\sqrt{2\pi}} \left(\frac{2\alpha_E}{a} - \frac{a}{2\alpha_E}\right) = \frac{(\varepsilon_t^2 t)^{1/6}}{\sqrt{2\pi}} \left[\frac{(pc+1)}{\sqrt{pc}} - \frac{\sqrt{pc}}{(pc+1)}\right] \approx \frac{(\varepsilon_t^2 t)^{1/6}}{\sqrt{2\pi}} \left[\sqrt{pc} + \frac{1}{pc^{3/2}} + O(pc^{-5/2})\right]. \qquad (109)$$

However, because the term $O(pc^{-3/2})$ in (109) is positive as opposed to negative in (108), an extra factor of $2^{-1/3}$ is included in (109) to ensure it never exceeds (108). This condition is guaranteed provided $pc > 2.463 \equiv Y > 1.103$, a fact already established by the analysis leading to (97).

*4.3.5 Operational count of second sum of $ZP(t)$ for pivots $\alpha_E \in (Ya, \alpha_E^c]$.*

Computation of this element of the sum proceeds in an analogous manner to that described in section 4.2.3. The interval $\alpha_E \in (Ya, \alpha_E^c]$ is subdivided into more pivot 'Blocks' numbered $p = L + 1, L + 2, .., L + L_c$, only in this computational region the number of pivots within each block is fixed at $\lfloor X^L t^{1/4} \rfloor$. At the end of each Block the step parameter is increased, this time regulated by equation (109), much as before. The computation ceases when $\alpha_{E, L + L_c} > \alpha_E^c = INT_E[(\varepsilon_t t^2)^{1/3}/ log(t)]$ the designated maximum specified in *Theorem 2*. The first question to consider is the actual size of $L_c \in \mathbb{N}$ necessary to ensure $\alpha_{E, L + L_c} > \alpha_E^c$. Continuation of equation (88) gives

$$\alpha_{E,L+L_c} = \alpha_{E,L} + 2\lfloor X^L t^{1/4} \rfloor \sum_{i=1}^{L_c}(M_{t,L+i} + 1) \approx Ya + \frac{2e^{(28h_t-7)/3}t^{1/3}}{\varepsilon_t^{1/3}} \sum_{i=1}^{L_c}(M_{t,L+i} + 1). \quad (110)$$

The next step is to establish a lower bound estimate for the various step parameters $M_{t,L+1,\dots,L+L_c}$. Initially

$$M_{t,L+1} = M_t(\alpha_{E,L}) = \frac{(\varepsilon_t^2 t)^{1/6}}{\sqrt{2\pi}} \left( \frac{2\alpha_{E,L}}{a} - \frac{a}{2\alpha_{E,L}} \right) = \frac{(\varepsilon_t^2 t)^{1/6}}{\sqrt{2\pi}} 2Y\left(1 - \frac{1}{4Y^2}\right),$$

$$\Rightarrow \alpha_{E,L+1} = Ya + e^{(28h_t-7)/3} aY\left(1 - \frac{1}{4Y^2}\right) = YF_1 a, \quad (111)$$

where $F_1 = [1 + e^{(28h_t-7)/3}(1 - 1/4Y^2)]$. For the next Block,

$$M_{t,L+2} = M_t(\alpha_{E,L+1}) = \frac{(\varepsilon_t^2 t)^{1/6}}{\sqrt{2\pi}} \left( \frac{2\alpha_{E,L+1}}{a} - \frac{a}{2\alpha_{E,L+1}} \right) = \frac{(\varepsilon_t^2 t)^{1/6}}{\sqrt{2\pi}} 2YF_1 \left(1 - \frac{1}{4Y^2 F_1^2}\right),$$

$$\Rightarrow \alpha_{E,L+2} = YF_1 a + e^{(28h_t-7)/3} aYF_1 \left(1 - \frac{1}{4Y^2 F_1^2}\right) = YF_1 F_2 a, \quad (112)$$

where $F_2 = [1 + e^{(28h_t-7)/3}(1 - 1/4Y^2 F_1^2)]$. This pattern continues and one finds that

$$M_{t,L+j} = M_t(\alpha_{E,L+j-1}) = \frac{(\varepsilon_t^2 t)^{1/6}}{\sqrt{2\pi}} 2YF_1 F_2 \dots F_{j-1} \left(1 - \frac{1}{4Y^2 F_1^2 F_2^2 \dots F_{j-1}^2}\right),$$

$$\alpha_{E,L+j} = YF_1 F_2 \dots F_j a, \quad j = 1,2,\dots,L_c, \quad (113)$$

where $F_j = [1 + e^{(28h_t-7)/3}\{1 - 1/(4Y^2 F_1^2 F_2^2 \dots F_{j-1}^2)\}]$. Since $F_1 < F_2 < \dots < F_{L_c}$, by the end of the $p = (L + L_c)$th Block, one is guaranteed to have reached *at least* summand $\alpha_{E,L+L_c} = YF_1^{L_c} a$. This value of $YF_1^{L_c} a$, is certain to surpass the value $\alpha_E^c$ prescribed in *Theorem 2* for some $L_c$ satisfying

$$L_c < \frac{\log(t)}{6\log(F_1)} \left\{ 1 + \frac{6\log[\varepsilon_t^{1/3}/Y\sqrt{8}\log(t)]}{\log(t)} \right\} = Q\log(t). \quad (114)$$

Now the exact value of the parameter $Q$ is dependent upon the precise definition of $\varepsilon_t$. However, it is possible to get some idea of its value in the limit as $t \to \infty$, since the term in curly brackets will tend to unity. In this limit, $Y > 1.2414$ (see eq. 97) and $e^{28h_t/3} \to 1$, $\Rightarrow \log(F_1) > \log\left(1 + e^{-7/3}(1 - 1/4Y^2)\right) = \log(1.0812\dots) = 0.0781$, which in turn $\Rightarrow Q < 2.133\dots$ The numerical results mentioned earlier show that $Y \approx 1.6 \Rightarrow Q \approx 2$.

The values of $M_{t,L+j}$ associated with these Blocks are all sufficiently large for the associated quadratic Gaussian sums to be calculated using *algorithm QGS*. So the operational count of $ZP(t)$ with $\alpha_E \in [Ya, \alpha_E^c]$ will be similar to the operational count for Blocks $p = p^* + 1, \dots, L$ when *algorithm QGS* was first utilised, only on this occasion the number of pivots is fixed at

$\lfloor X^L t^{1/4} \rfloor$. Incorporating this change into the final two terms of (98), one obtains an operational count

$$\equiv \lfloor X^L t^{1/4} \rfloor \left[ \left(4\Omega + 4K - 2\Omega^* \log(2K)\right) \sum_{j=1}^{L_c} 1 + 2\Omega^* \sum_{j=1}^{L_c} \log(M_{t,L+j} - 1) \right] ops_{s;c}. \tag{115}$$

Since $F_j < (1 + e^{(28h_t - 7)/3}) \forall j$, $M_{t,L+j} < (\varepsilon_t^2 t)^{1/6} \sqrt{2} Y (1 + e^{(28h_t - 7)/3})^{j-1}/\sqrt{\pi}$, an upper bound on the logarithmic sum is given by

$$\sum_{j=1}^{L_c} \log(M_{t,L+j} - 1) < \log((\varepsilon_t^2 t)^{1/6} \sqrt{2} Y / \sqrt{\pi}) \sum_{j=1}^{L_c} 1 + \log(1 + e^{(28h_t - 7)/3}) \sum_{j=1}^{L_c} (j-1)$$

$$= L_c \left[ \log((\varepsilon_t^2 t)^{1/6} \sqrt{2} Y / \sqrt{\pi}) + \left(\frac{L_c - 1}{2}\right) \log(1 + e^{(28h_t - 7)/3}) \right]. \tag{116}$$

Substituting (116) into (115) gives an upper bound on the operational count as follows

$$\equiv \frac{e^{(28h_t - 7)/3} t^{1/3}}{\varepsilon_t^{1/3}} L_c \left[ \left(4\Omega + 4K - 2\Omega^* \log(2K)\right) + 2\Omega^* \log((\varepsilon_t^2 t)^{1/6} \sqrt{2} Y / \sqrt{\pi}) \right.$$
$$\left. + 2\Omega^* \left(\frac{L_c - 1}{2}\right) \log(1 + e^{(28h_t - 7)/3}) \right] ops_{s;c},$$

$$\equiv \frac{2\Omega^* Q t^{1/3} [\log(t)]^2}{\varepsilon_t^{1/3} e^{(7 - 28h_t)/3}} \left[ \left\{ \frac{1}{6} + \frac{Q}{2} \log(1 + e^{(28h_t - 7)/3}) \right\} + \frac{1}{\log(t)} \left\{ \left( \frac{2(\Omega + K)}{\Omega^*} \right) \right. \right.$$
$$\left. \left. + \log\left(\frac{\varepsilon_t^{1/3} \sqrt{2} Y}{\sqrt{\pi} K^2}\right) - \frac{1}{2} \log(1 + e^{(28h_t - 7)/3}) \right\} \right] ops_{s;c}. \tag{117}$$

*4.3.6 Total operational count for the evaluation of the sum $ZP(t)$, as defined in Theorem 2*

The total operational count for the evaluation of sum $ZP(t)$ of *Theorem 2*, using the *algorithm QGS* to estimate the quadratic Gaussian sums in the manner described in Sections 4.3.2-5, cannot exceed the combination of results (107) and (117). In full the evaluation of $ZP(t)$ requires no more than

$$< 2\Omega t^{1/4} \left[ 2X + (2\wp)^{1/3p^*} \left[ X^{-1/4} \left(\frac{2\varepsilon_t}{\pi}\right)^{\frac{1}{2}} t^{\frac{1}{8}} \right]^{1 - \frac{1}{3p^*}} \frac{X^2 (X^{3p^*/2} - X^{3/2})}{(X^{3/2} - 1)^{4/3}} - \frac{X(X^{p^*} - X)}{(X - 1)} \right] +$$

$$\frac{\Omega^* Q t^{\frac{1}{3}} [\log(t)]^2}{\varepsilon_t^{1/3} e^{(7 - 28h_t)/3}} \left[ \left\{ \frac{1}{3} + Q \log(1 + e^{(28h_t - 7)/3}) \right\} + \frac{1}{\log(t)} \left\{ \frac{4(\Omega + K)}{\Omega^*} + \log\left(\frac{2\varepsilon_t^{2/3} Y^2}{\pi K^4 [1 + e^{(28h_t - 7)/3}]}\right) \right\} \right] +$$

$$\frac{\Omega^* \varepsilon_t^{-1/3} t^{1/3} \log(t)}{e^{9/4 - 9h_t} (X - 1)} \left[ \left\{ B_1 \log(X) + \frac{B_0}{4} \right\} + \frac{1}{\log(t)} \left\{ B_0 \log\left(\frac{\varepsilon_t}{2\pi K^2 \sqrt{X}}\right) - 28 B_1 \log(X) + \frac{4(\Omega + K) B_1}{\Omega^*} \right\} \right] \tag{118}$$

standard sine/cosine function $ops_{s;c}$. Suppose that the termination parameter $K$ is chosen lie below $M_{t,2}^- \approx (4X\mathcal{G})^{1/3}[\varepsilon_t^4 t/\pi^4]^{1/12}/2 \approx \sqrt{3}[\varepsilon_t^4 t/\pi^4]^{1/12}/2$ as specified by postulate *i)* of *Theorem 2* and consider the limit of (118) as $t \to \infty$. As $K < M_{t,2}^-$ this implies $p^* = 1$ (see discussion after eq. 94) and hence the first term in (118) collapses to a negligible $4\Omega X t^{1/4}$ value. In the remaining two terms the all the factors $e^{h_t}, B_0\ B_1 \xrightarrow[t \to \infty]{} 1, \log(X) \xrightarrow[t \to \infty]{} 1/12$, whilst the terms of the form $\log(\varepsilon_t/K^2)/\log(t)$, whilst not necessarily tending to zero in this limit, would be dominated by their corresponding $K/\log(t)$ terms. Hence the evaluation of $ZP(t \to \infty)$ with $K < \sqrt{3}[\varepsilon_t^4 t/\pi^4]^{1/12}/2$, requires no more than

$$< \frac{\Omega^* Q t^{\frac{1}{3}}[\log(t)]^2}{\varepsilon_t^{1/3} e^{7/3}} \left[ \left\{\frac{1}{3} + Q\log(1 + e^{-7/3})\right\} + \frac{1}{\log(t)}\left\{\frac{4(\Omega + K)}{\Omega^*} + \frac{e^{1/12}}{3Q(X-1)}\right\} + O\left(\frac{1}{[\log(t)]^2}\right) \right] ops_{s;c}$$

$$\approx 3.14 \left(\frac{t}{\varepsilon_t}\right)^{1/3} [[\log(t)]^2 + [\log(t)]\{0.247(\Omega + K) + 4.02\} + O(1)] ops_{s;c}$$

$$< 8.38 \left(\frac{t}{\varepsilon_t}\right)^{1/3} [[\log(t)]^2 + [\log(t)]\{0.099(\Omega + K) + 3.68\} + O(1)] ops_{s;c}. \quad (119)$$

The approximate result in (119) is obtained by substituting in the average values for $\Omega^* = 0.59 \times 53$ and $Q \approx 2$ established at the end of section 3 and after (114) respectively. The upper bound for (118) is obtained using the respective maximum values $\Omega_{max}^* = 1.44 \times 53$ and $Q < 2.133$. The value of $\Omega \approx 1.17 - 1.30$.

*These operational counts are consistent with and so prove postulate i) of Theorem 2.*

The numerical factors that appear here obviously depend upon the (conservative) conversion factors $ops_{RS} \equiv 2ops_{s;c}$ and $ops_{hybrid} \equiv 2\Omega ops_{s;c}$ prescribed in Section 4.2. In reality these will vary slightly depending upon the operating system being used, but cannot change the overall conclusion. Section 5.2 compares (118) to actual operational counts derived from sample computations.

For larger values of $K$ things are a little more complicated. If $K > \sqrt{3}[\varepsilon_t^4 t/\pi^4]^{1/12}/2$ then $p^* > 1$ which means that the first term in (118) can contribute significantly to the overall operational count. In this instance $K\log(t) \gg [\log(t)]^2$ in (119), and so the most significant terms in the overall operational count (118) would comprise

$$\left\{ O\left(t^{1/4}[t^{1/8}\sqrt{\varepsilon_t}]^{1-1/3^{p^*}} X^{3p^*/2}\right) + O\left(K(t/\varepsilon_t)^{1/3}\log(t)\right) \right\} ops_{s;c}. \quad (120)$$

When $\sqrt{3}[\varepsilon_t^4 t/\pi^4]^{1/12}/2 < K \ll \ell(\varepsilon_t^2 t)^{1/6}$ the value of $p^*$ is relatively small and the second term is dominates the overall operational count, which remains $\ll O(\sqrt{t}) ops_{s;c}$ as $t \to \infty$. *This is consistent with postulate ii) of Theorem 2*. It is only when $K$ starts to approach the value $\ell(\varepsilon_t^2 t)^{1/6}$ at which point $p^* \to L - 1$, do the two terms in (120) start to coalesce to around $\left\{ O(\sqrt{t}) + O(\sqrt{t}\log(t)) \right\} ops_{s;c}$. Of course fixing $K$ at such a large value effectively means that all the quadratic Gaussian sums making up $ZP(t)$ are being evaluated long hand (since $K \sim M_{t,L+j}$; c.f. eq. 112), with no computational benefits accruing from the use of *algorithm QGS*. This completes the proof of *Theorem 2*.

The implications of this result for the operational complexity of $Z(t)$ in the limit $t \to \infty$ are summarised in the following theorem.

### 4.4 Theorem 3

*Let $K, K_0, P \in \mathbb{N}$ and set $K_0 = 50$. For any $t > (5\sqrt{\pi}K_0)^6 > 0$, define $a = \sqrt{8t/\pi}$ and choose $\varepsilon_t$ as a fixed parameter prescribed so that $\varepsilon_t \xrightarrow[t \to \infty]{} 0$ and $\varepsilon_t^{-1} = o(t^\mu)$ for any $\mu > 0$. Define $h_t = \log(\varepsilon_t)/\log(t)$ and choose $K$ to lie in the range $K_0 \leq K \leq MAX\{K_0, \ell(\varepsilon_t^2 t)^{1/6}\}$, where $\ell = \exp(55[4h_t - 1]/48)/\sqrt{2\pi}$. Fix $b > 0$ ($b \ll a$), and let $\alpha \in (a, a+b)$ denote the odd integers in this range. Fix $\alpha_E^c = INT_E\left[(\varepsilon_t t^2)^{1/3}/(\sqrt{\pi}\log(t))\right], N_c = \lfloor \alpha_E^c\{1 - \sqrt{1 - (a/\alpha_E^c)^2}\}/4 \rfloor \approx t/\pi\alpha_E^c$ and let $\alpha_E \in [a+b, \alpha_E^c]$ denote elements of a subset of the even (pivot) integers within the given interval. Define the sum $ZP(t)$ as given by equations (82-3) of Theorem 2. Then utilising a* **refined version** *(see 4.5.2b) of algorithm QGS with parameters $(K, P)$ to evaluate the sub-sums of $ZP(t)$, as discussed in sub-sections 4.3.1-4.3.5, one can obtain the following approximation for the majority of terms making up the Hardy function $Z(t)$*

$$\sum_{N=N_c}^{N_t} \frac{2\cos\{\theta_C(t) - t\log(N)\}}{\sqrt{N}} = ZP(t)[1 + O(\varepsilon_t)], \qquad (121)$$

*for an additional operational cost of $O(\log(t/\varepsilon_t^7 K^3))$ times the bounds specified in Theorem 2 (in particular section 4.3.6).*

*Corollary*

  *Suppose $K$ is chosen so that $K_0 \leq K \leq MAX\{K_0, \sqrt{3}[\varepsilon_t^4 t/\pi^4]^{1/12}/2\}$. Then the calculation of the Hardy function $Z(t)$ to a precision of $|O(\varepsilon_t) \times ZP(t)|$, requires only*

$$O\left(\left(\frac{t}{\varepsilon_t}\right)^{1/3} \log(t)[\log(t) + K] \times \log(t/\varepsilon_t^7 K^3)\right) ops_{S;c} \text{ as } t \to \infty.$$

### 4.5 Proof of Theorem 3

The statement of *Theorem 3* expresses two distinct points. Firstly it links the sum $ZP(t)$ defined in *Theorem 2*, to a partial sum of the terms in the RSF (8). Secondly it makes a statement about the size of the relative errors that can accrue if $ZP(t)$ is used as an approximation for this partial sum. These points will be discussed separately.

### 4.5.1 Connection of the sum $ZP(t)$ to the Riemann-Siegel formula.

This is relatively straightforward since it follows almost automatically from *Theorem 1*, quoted in the Introduction (q.v. eq. 14), links the sum of terms arising from the first order approximations of the various confluent hypergeometric functions in (10, 18) to the terms of (8), to within a guaranteed tolerance of $O(t^{-1/12})$. Since $ZP(t)$ forms a proportion of the full sum (54), which itself arises when one approximates the same confluent hypergeometric functions (18) into $M_t$ sized collections (Section 2, albeit with larger *relative* $O(\varepsilon_t)$ errors than the absolute error of *Theorem 1*), then $ZP(t)$ too is linked to a proportion of the terms making up the main sum of (8). The actual number depends upon the cut-off value $\alpha_E^c$, or more precisely the corresponding value of $\alpha_E^c + M_t(\alpha_E^c)$, which is the exact number of terms from the full sum represented by $ZP(t)$. Any choice cut-off satisfying $a < \alpha_E^c + M_t(\alpha_E^c) \leq N_\alpha \approx (t/\pi)$ is linked to a corresponding cut-off $N_c$ in (8) by means of equation (15). (This can be easily established using the exactly same methodology employed by [27 eq. 8] to proof *Theorem 1* for the full sum, but over a smaller circle $C(R)$ centred at $z = (\alpha_E^c + a)/2$ with radius $R \approx (\alpha_E^c - a)/2$.) Since $M_t(\alpha_E^c) \ll \alpha_E^c$ it can be ignored when inverting (15). This gives a first order approximation for $N_c$

$$N_c \approx \frac{\alpha_E^c[1 - \sqrt{1 - (a/\alpha_E^c)^2}]}{4} \approx \frac{a^2}{8\alpha_E^c} = \frac{t}{\pi \alpha_E^c} \approx \frac{t^{1/3} \log(t)}{\sqrt{\pi} \varepsilon_t^{1/3}}, \qquad (122)$$

accurate to $O(\{log(t)\}^3/\varepsilon_t)$. Hence the sum $ZP(t)$ corresponds to an approximation of the sum of $N \in [N_c, N_t]$ terms of the Riemann-Siegel formula, which represents the overwhelming majority of those needed to compute $Z(t)$.

*4.5.2 Error Analysis*

This raises some more difficult issues. There are two sources of error that arise in the approximation (121). These must be examined in much more detail than hitherto.

*4.5.2a) Errors arising from the derivation of the $ZP(t)$ approximation.*

The first source of error concerns those terms, loosely denoted by $O(\varepsilon_t)$ in Section 2, arising from the *derivation* of the connection between $Z(t)$ and the series of Gaussian quadratic sums, expressed by (54) of which $ZP(t)$ forms a part. The first such error arises from the estimation of the integral $I(A, B, C)$, defined by (34-35), using approximation (38) with the exponential factor $exp\left(\pm i \frac{j^2 B^2}{4AC^4} \delta\right)$ set equal to unity. This is justified by fixing $max(j) = M_t$ to lie below the upper bound given by (40), which ensures that $\delta = jB/AC^3$ satisfies $\delta < j^2 B^2 \delta/4AC^4 \leq M_t^3 B^3/4A^2 C^7 = O(\varepsilon_t)$. Suppose instead that the exponential factor is written as $1 \pm i \frac{j^2 B^2}{4AC^4} \delta$. The effect of including the $\pm i \frac{j^2 B^2}{4AC^4} \delta$ term explicitly in the analysis spelt in Section 2 (cf. eqs. 27-30 and 47-54) is to produce a first order correction to equation (54) which is *structurally similar*, except for the quadratic Gaussian sums $e^{i\omega^\pm(pc)} S_{M_t^-}^*(x(pc), \theta^\pm(pc))$ (also the substitution of the Imaginary part as opposed to the Real part), which are replaced by

$$\pm \frac{iM_t^3}{(pc^2 - 1)^2} \sqrt{\frac{2\pi^3 pc^5}{t}} e^{i\omega^\pm(pc)} \sum_{k=0}^{M_t^-} \left(\frac{2k+1}{M_t}\right)^3 e^{i\pi k[kx(pc) + 2\theta^\pm(pc)]}. \qquad (123)$$

The multiplicative factor in front arises from (39) and is $O(\varepsilon_t)$ from (40), whilst $2k+1$ has been substituted for $j$ in the summation part. The *structural similarity* stated above pertains to the exponential phases in (123), which are identical to the corresponding phases of the $e^{i\omega^{\pm}(pc)}S^*_{M_t^-}(x(pc),\theta^{\pm}(pc))$ terms in (54). This is important for the following reason. One would expect the sum (123) to be somewhat smaller than the corresponding $S^*_{M_t^-}(x(pc),\theta^{\pm}(pc))$ which occurs in the main part of (54), since all the amplitudes $(2k+1)/M_t < 1$. However, if the phase structure of (123) were to differ from (54) then when collections of these smaller sub-sums (each characterised by a different pivot value) are themselves added together, they could give rise (by some extremely fortuitous correlation in the signs of the sub-sums) to a first order correction to (54) and hence $ZP(t)$, very much larger than the predicted $O(\varepsilon_t)$. So it is necessary to check this cannot occur.

The correction term (123) can, to $O(M_t^{-1})$, be written in terms of triple derivative of $S^*_{M_t^-}(x(pc),\theta^{\pm}(pc))$ with respect to $\theta^{\pm}(pc)$

$$O(\varepsilon_t) \times e^{i\omega^{\pm}(pc)} \left(\frac{1}{i\pi M_t}\right)^3 \frac{d^3}{d\theta^{\pm 3}}\left[S^*_{M_t^-}(x(pc),\theta^{\pm}(pc))\right], \qquad (124)$$

If the pivot is such that $pc > 2$, then $0 < x \equiv x(pc) = (pc-1)^{-1} < 1$. Alternatively if $pc \in (1,2]$ replace $x \equiv frac[x] \in [0,1]$ and adjust $\theta^{\pm} \equiv \theta^{\pm}(pc)$ appropriately so that $-1/2 < \theta^{\pm} < 1/2$. In each case *Theorem[Paris]* is then applicable. Suppose $x(pc)$ is a 'generic' irrational in the sense elucidated in section 3.5, in which case the quadratic Gaussian sum in (124) is dominated by the first term of (58, 63). Substituting this term into (124) allows one to evaluate the derivatives explicitly giving, to $O(M_t^{-1})$

$$O(\varepsilon_t)\, e^{i\omega^{\pm}(pc)} \frac{e^{-\frac{i\pi\theta^{\pm 2}}{x}}}{\omega\sqrt{x}} \sum_{l=0}^{3} \binom{3}{l} \left(\frac{-2\theta^{\pm}}{M_t x}\right)^{3-l} \sum_{k=1}^{\lfloor\xi\rfloor} \left(\frac{2k}{M_t x}\right)^l e^{\frac{i\pi k[-k+2\theta^{\pm}]}{x}}, \qquad (125)$$

with $\xi = M_t^- x + \theta^{\pm} > 1$. The $l=0$ term in (125) corresponds to the first in the hierarchal chain of (shorter) quadratic Gaussian sums $S_{\lfloor\xi\rfloor}(-1/x,\theta^{\pm}/x)$ generated from $S^*_{M_t^-}(x(pc),\theta^{\pm}(pc))$ by repeated application of (58). The $l > 0$ terms correspond to derivatives of this first shorter sum. Just as (125) follows from (123), these derivatives can themselves be written in terms of derivatives of the *next* shorter sum, simply by applying *Theorem[Paris]* a second time. This process can be repeated recursively, just as in *algorithm QGS*, until termination is reached when $\xi = M_t^- x + \theta^{\pm} < 1$. Hence the hierarchal chain of sums originating from $S^*_{M_t^-}(x(pc),\theta^{\pm}(pc))$ in (58) would also be a feature of the derivatives (124-5) and the phase structure of the original sum would be repeated in its derivatives. Complicating things are the whole series of extra amplitudes, consisting of products of factors, commencing with powers of $(2\theta^{\pm}/M_t x)$, subsequently followed by powers of $(\theta_{n\pm}/L_n|x_n|)$ down the hierarchal chain. But these extra amplitudes can only increase the size of the error term if their various constituent factors have moduli greater than unity. Suppose the pivot is such that $pc \in (1,2]$ and $x \equiv frac[x]$ is generic, then the factor $|2\theta^\pm/M_t x| \approx |\theta^\pm/M_t^- x| \leq 1$, since $\xi \geq 1$. Alternatively if $pc > 2$, then the factor $|2\theta^\pm/M_t x| = |2\theta^\pm(pc-1)/M_t|$. Using the general definition (40) for $M_t$, this is bounded above by

$$\left|\frac{2\theta^\pm(pc-1)}{M_t}\right| \approx \left|\frac{2\theta^\pm(pc-1) \times (2\pi^3 pc^5)^{\frac{1}{6}}}{(pc^2-1)^{2/3}(\varepsilon_t \sqrt{t})^{1/3}}\right| < 2^{\frac{1}{6}}\sqrt{\pi}\left|\frac{2\theta^\pm\sqrt{pc}}{(\varepsilon_t\sqrt{t})^{1/3}}\right| < \sqrt{\frac{\pi}{2}}\frac{|2\theta^\pm|\alpha_E}{(\varepsilon_t t^2)^{1/3}}$$

$$\leq \sqrt{\frac{\pi}{2}}\frac{|2\theta^\pm|\alpha_{Emax}}{(\varepsilon_t t^2)^{1/3}} \approx \sqrt{\frac{\pi}{2}}\frac{|2\theta^\pm|}{\log(t)} < 1. \tag{126}$$

Hence all the initial factors in the chain satisfy $|2\theta^\pm/M_t x| < 1$. (It is an interesting 'coincidence' to note that the cut off $\alpha_{Emax} = INT_E[(\varepsilon_t t^2)^{1/3}/\log(t)]$, originally proposed in *Theorem 2* with a view to optimizing the computation of $Z(t)$, also bounds the error by guaranteeing the initial amplitude factors in (125) never exceed unity.) Moving on down the hierarchal chain, all the subsequent factors satisfy $|\theta_{n\pm}/L_n|x_n|| < 1$, since $(L_n|x_n| + \theta_{n\pm}) > 1$ and hence $|\theta_{n\pm}| < L_n|x_n|$, to the termination point of the chain. Hence (124) consists of an $O(\varepsilon_t)$ term multiplied by a term which cannot be any larger in magnitude than the size of $S^*_{M_t^-}\left(x(pc), \theta^\pm(pc)\right)$ itself. This fact, combined with the observation that the sub-sum (123-4) shares the same phase structure as $S^*_{M_t^-}\left(x(pc), \theta^\pm(pc)\right)$, means that when collections of such sub-sums are added together as prescribed by (82), the final total cannot give anything larger than the postulated $O(\varepsilon_t)$ correction to $ZP(t)$.

When $x(pc)$ happens to be a 'non-generic' irrational things are slightly more awkward. In such instances at some point in the hierarchal chain a large partial quotient is encountered, resulting in a very small $|x_n|$ value. This abruptly terminates the chain of Gaussian sums since $L_n|x_n| < 1$. For simplicity, suppose this termination happens immediately (if the termination occurs further down the hierarchal chain, the analysis is entirely analogous). Then the behaviour of the sum $S^*_{M_t^-}\left(x(pc), \theta^\pm(pc)\right)$ in (124) will no longer be governed the leading order term in (63), as this will be cancelled out by the first term in square brackets of (64). Instead it is the second set of terms, $\left[\text{erf}\left(\frac{\omega(M_t^- x + \theta^\pm)\sqrt{\pi}}{\sqrt{x}}\right) - \text{erf}\left(\frac{\omega\theta^\pm\sqrt{\pi}}{\sqrt{x}}\right)\right]e^{-\frac{i\pi\theta^{\pm 2}}{x}}/2\omega\sqrt{x}$ in (64), which will be relevant. The behaviour of these terms is complicated, but there are two main branches valid for small $x$.

$$\begin{cases} \dfrac{i}{2\pi}\left(\dfrac{1}{\theta^\pm} - \dfrac{e^{i\pi M_t^-(M_t^- x + 2\theta^\pm)}}{(M_t^-|x| + \theta^\pm)}\right), & \text{for } M_t^- x \ll |\theta^\pm|, \left|\dfrac{\theta^\pm\sqrt{\pi}}{\sqrt{x}}\right| \gg 1 \text{ and } \left|\dfrac{(M_t^- x + \theta^\pm)\sqrt{\pi}}{\sqrt{x}}\right| \gg 1, \\ \\ \dfrac{e^{-\frac{i\pi\theta^{\pm 2}}{x}}}{2\omega\sqrt{x}}, & \text{for } |\theta^\pm| < M_t^- x, \left|\dfrac{\theta^\pm\sqrt{\pi}}{\sqrt{x}}\right| < 1 \text{ and } \left|\dfrac{(M_t^- x + \theta^\pm)\sqrt{\pi}}{\sqrt{x}}\right| \gg 1. \end{cases}$$

(127a, b)

The most common case (127a), pertains when $|\theta^\pm|$ is not itself close to zero. An example of this is the termination exhibited in 3.4 Case D. Substituting (127a) into (124) and differentiating gives rise to (to first order) extra amplitudes $\equiv 4(M_t^-)^3/\pi M_t^3 = 1/2\pi < 1$. The less common case (127b) occurs when $|\theta^\pm|$ is also close to zero, small enough to ensure $|\theta^\pm| < M_t^- x$. An example of this kind of termination is exhibited in 3.4 Case E. Substituting (127b) and differentiating gives rise (to first order) extra amplitudes $\left|2\theta^\pm/M_t x\right|^3 \approx \left|\theta^\pm/M_t^- x\right|^3 < 1$, similar to the 'generic' case. So in both cases the extra amplitudes have modulus less than one, which means that (just as in the 'generic' case) the correction terms (123-4) are no larger in magnitude than $S^*_{M_t^-}\left(x(pc), \theta^\pm(pc)\right)$ itself. Since the associated phase structures are also preserved, one arrives at the same conclusion as for the 'generic' case, namely that when the various sub-sums (123) are added together their total will not exceed the postulated $O(\varepsilon_t)$ correction to $ZP(t)$.

Another $O(\varepsilon_t)$ relative error term arising from the *derivation* analysis of Section 2, pertains to the secondary integral appearing in (33). But as was demonstrated in (44), this secondary integral is no more than the main integral of (33) multiplied by a factor $\pm iO(\varepsilon_t)$. So the inclusion of such secondary integrals in the computation of $ZP(t)$ still only produces an $O(\varepsilon_t)$ correction to the estimate for the partial sum of (8), in concordance with *Theorem 3*.

4.5.2b) *Errors associated with the computation of $ZP(t)$*.

The second source of error arises in the *computation* of the various quadratic Gaussian sums by means of *algorithm QGS*. In section 3.5 it was shown that for a *single* general quadratic Gaussian sum the relative error $\varepsilon_{GS}(K)$ was always less than the crude bound $1/2\sqrt{K}$, for a termination constant $K$ (q.v. eqs. 76-77). In itself this is fine, but what is envisaged in *Theorem 3* is the computation of $ZP(t)$, which is a *sum* of multiple combinations of Gaussian sums, all estimated using *algorithm QGS*. Again this means there is the potential to obtain, through some (extremely) fortuitous correlation in the signs of the errors of the various sub-sums, a much larger overall error in the calculation of $ZP(t)$ than one would envisage from (77). The problem is that, unlike for the *derivation* errors discussed above, *structural similarity* is not preserved here. The phase of the overall error $\delta(K, N, x, \theta, P)$ in *algorithm QGS* is determined by combinations of the phases of the various $\varepsilon_{0-n_k}$ terms arising from each iteration step (73-74), and these phases will be quasi-random. Hence one cannot guarantee that some fortuitous correlation in signs will not occur. To overcome this problem one needs a methodology to reduce the overall error in a systematic way, should the need arise. The value of $\delta(K, N, x, \theta, P)$ is determined from the remainder terms (61) which arise from asymptotic expansions of the complimentary error function erfc, and as is pointed out in [32 Sections 2.11(iii) & 7.12(i)] such expansions can be *"exponentially improved"*. By employing similar methods to those of [4] it is possible to derive very much more accurate estimates for the remainder terms than those presented in [33]. Unfortunately these methods, on their own, are not enough for the purposes of proving *Theorem 3*, because although yielding a superior estimate for the remainder term, it still suffers from algebraic *not* exponentially small, relative errors. And the presence of these algebraic errors means that the resulting reduction in $\varepsilon_{GS}(K)$ is still insufficient to guarantee that the sum of all the relative errors, arising from the $O\left(t^{1/4}[L_{max}X^L + \sum_{p^*+1}^{L} X^p] \approx Qe^{-7/3}log(t)[t/\varepsilon_t]^{1/3}\right)$ Gaussian subsums making up $ZP(t)$, cannot combine to give an unacceptably large value. Nevertheless, refining the estimates of the remainder terms provides the key, because the resulting analysis highlights the means of ensuring that each individual $\varepsilon_{GS}(K)$ associated with a particular subsum can be set to any specified precision (without at the same time drastically increasing the order of the operational count estimate in *Theorem 2*). In which case the relative error in the computation of $ZP(t)$ can be reduced to the same $O(\varepsilon_t)$ scale as the *derivation* errors discussed previously, irrespective of the details of the constituent phase structure in the various $\varepsilon_{GS}(K)$. The first step towards attaining this desired level of accuracy is to make a more detailed examination of the remainder terms (61) of *Theorem[Paris]*.

The remainder terms $\tilde{R}'_P$ of (58 & 64) come about through the method employed by [33 eq. 2.5] to estimate the following integral

$$E(z_k^\pm) = e^{z_k^{\pm 2}/2}\text{erfc}\{z_k^\pm/\sqrt{2}\} = \frac{e^{z_k^{\pm 2}/2}z_k^\pm}{\sqrt{2\pi}}E_{1/2}(z_k^{\pm 2}/2). \tag{128}$$

The second expression follows from the identity (7.11.3) of [32] which links $\text{erfc}\{z\}$ to $E_{1/2}(z)$ the generalized exponential integral [32 section 8.19]. It applies for all $z$. In this instance [33 eq. 2.3]

$$z_k^\pm = \sqrt{\frac{2\pi}{x}}(k \pm \xi)e^{-i\pi/4}, \qquad x > 0, \qquad k = 1,2,\dots \tag{129}$$

Since $\arg(z) = -\pi/4$ and hence $\arg(z^2) = -\pi/2$ for most $z_k^\pm$ values, it will prove useful to establish the asymptotic expansion for $E_{1/2}(z)$ in the particular case when $z = -i\tau$ with $\tau > 0$. This is easily done by means of results given in [32], reproduced below:

$$E_{1/2}(z) = z^{-1/2}\Gamma(1/2, z), \qquad (8.19.2, \quad z \in \mathbb{C})$$

$$= \frac{e^{-z}}{\sqrt{\pi}}\int_0^\infty \frac{t^{-1/2}e^{-t}}{t+z}dt, \qquad (8.6.4, \quad \arg(z) < \pi)$$

$$= \frac{e^{-z}}{\sqrt{\pi}}\left[\sum_{r=0}^{P-1}\frac{(-1)^r}{z^{r+1}}\int_0^\infty t^{r-1/2}e^{-t}dt + \frac{(-1)^P}{z^P}\int_0^\infty \frac{t^{P-1/2}e^{-t}}{t+z}dt\right],$$

$$= \frac{e^{-z}}{\sqrt{\pi}}\left[\sum_{r=0}^{P-1}\frac{(-1)^r\Gamma(r+1/2)}{z^{r+1}} + \frac{(-1)^P}{z^P}\int_0^\infty \frac{t^{P-1/2}e^{-t}}{t+z}dt\right]. \tag{130}$$

The last expression follows from $(t + z)^{-1} = z^{-1}(1 + t/z)^{-1}$ and the application of [32 eq. 2.11.9]. Now in the special case $\arg(z = -i\tau) = -\pi/2$, one can use the substitution $t = izu^2 = \tau u^2$ to transform the integral in (130) into a more convenient form, giving

$$E_{1/2}(z = -i\tau) = \frac{e^{-z}}{\sqrt{\pi}} \left[ \sum_{r=0}^{P-1} \frac{(-1)^r \Gamma(r+1/2)}{z^{r+1}} + \frac{2(-i)^P e^{i\pi/4}}{z^{1/2}} \int_0^\infty \frac{u^{2P} e^{-i\tau u^2}}{1+iu^2} du \right]. \quad (131)$$

Now (131) is suitable for $z_k^+$ and $z_{k>\xi}^-$ as it stands, since both satisfy $\arg(z^2) = -\pi/2$. However, since $z_{k<\xi}^-$ satisfies $\arg(z^2) = 3\pi/2$, the use of equation (8.6.4) in the derivation of (130) is invalid in this instance. To overcome this problem one must invoke the analytical continuation formula [32 eq. 8.19.18] for the general function $E_{1/2}(z)$ with $z \in \mathbb{C}$. This gives

$$E_{1/2}\left(e^{i(2\pi - \pi/2)}|z_{k<\xi}^-|^2/2\right) = \left[\frac{2\pi i e^{i\pi/2}}{\sqrt{\pi}} \frac{\sqrt{2}}{|z_{k<\xi}^-| e^{-i\pi/4}} + E_{1/2}\left(e^{-i\pi/2}|z_{k<\xi}^-|^2/2\right)\right]. \quad (132)$$

Since the argument is now appropriate, one can use (131) to replace the exponential function on the right of (132). Using (128), (131) & (132) one obtains

$$E(z_{k<\xi}^-) = \frac{e^{-i|z_{k<\xi}^-|^2/2}}{\sqrt{2\pi}} e^{i3\pi/4} |z_{k<\xi}^-| \left[ \frac{2\sqrt{2\pi} e^{i5\pi/4}}{|z_{k<\xi}^-|} + E_{1/2}\left(e^{-i\pi/2}|z_{k<\xi}^-|^2/2\right) \right]$$

$$= e^{-i|z_{k<\xi}^-|^2/2} \left[ 2 + \frac{e^{i3\pi/4}|z_{k<\xi}^-|}{\sqrt{2\pi}} E_{1/2}\left(e^{-i\pi/2}|z_{k<\xi}^-|^2/2\right) \right]$$

$$= \left[ 2e^{-i|z_{k<\xi}^-|^2/2} + \frac{e^{i3\pi/4}|z_{k<\xi}^-|}{\pi\sqrt{2}} \left\{ \sum_{r=0}^{P-1} \frac{(-1)^r \Gamma(r+1/2)}{\left[-i|z_{k<\xi}^-|^2/2\right]^{r+1}} + \frac{2(-i)^P e^{i\pi/4}}{\left[-i|z_{k<\xi}^-|^2/2\right]^{1/2}} \int_0^\infty \frac{u^{2P} e^{-|z_{k<\xi}^-|^2 u^2/2}}{1+iu^2} du \right\} \right]$$

$$= \left[ 2e^{-i|z_{k<\xi}^-|^2/2} - \frac{1}{\pi}\sum_{r=0}^{P-1} \frac{(-1)^r \Gamma(r+1/2)}{\left[-i|z_{k<\xi}^-|^2/2\right]^{r+1/2}} - \frac{2(-i)^P}{\pi e^{-i\pi/4}} \int_0^\infty \frac{u^{2P} e^{-|z_{k<\xi}^-|^2 u^2/2}}{1+iu^2} du \right]. \quad (133)$$

Equation (133) is actually equation (2.6b) of [33], but with the remainder term now written explicitly in terms of a specific integral. The case for $z_k^+$ and $z_{k>\xi}^-$ is much easier, since (131) can be substituted into (128) directly, giving [33 q.v. eq. 2.6a]

$$E(z_{k>\xi}^\pm) = \frac{1}{\pi} \sum_{r=0}^{P-1} \frac{(-1)^r \Gamma(r+1/2)}{\left[-i|z_{k>\xi}^\pm|^2/2\right]^{r+1/2}} + \frac{2(-i)^P}{\pi e^{-i\pi/4}} \int_0^\infty \frac{u^{2P} e^{-|z_{k>\xi}^\pm|^2 u^2/2}}{1+iu^2} du. \quad (134)$$

The remainder terms in (133-134), denoted by $-R_P(-z_{k<\xi}^-)$ and $R_P(z_{k>\xi}^\pm)$ respectively in [33 eq. 2.6] but not actually defined there, satisfy

$$R_P(z_{k>\xi}^\pm) = \frac{2(-i)^P}{\pi e^{-i\pi/4}} \int_0^\infty \frac{u^{2P} e^{-|z_{k>\xi}^\pm|^2 u^2/2}}{1+iu^2} du, \quad -R_P(-z_{k<\xi}^-) = -\frac{2(-i)^P}{\pi e^{-i\pi/4}} \int_0^\infty \frac{u^{2P} e^{-|z_{k<\xi}^-|^2 u^2/2}}{1+iu^2} du$$

(135)

The relatively crude estimates for the integral in (135) used in [33] are what underlie the error bounds calculated from (61-62). For the general application of *algorithm QGS* these bounds are more than adequate, but for computation of $ZP(t)$ envisaged in *Theorem 3* something more precise is needed.

Both integrals in (135) can be written in the generic form

$$I = \int_0^\infty \frac{\exp[2P\log(u) - \tau u^2]}{1+iu^2} du, \quad \text{with } \tau > 0. \tag{136}$$

Simple calculus reveals that the exponential phase $f(u) = 2P\log(u) - \tau u^2$ possesses a saddle at $u = \sqrt{P/\tau} = \sqrt{\varrho}$ and $f''(u = \sqrt{\varrho}) = -4\tau$. Substituting $v = u - \sqrt{\varrho}$ into (136) gives

$$I = \left(\frac{P}{\tau}\right)^P e^{-P} \int_{-\sqrt{\varrho}}^{\sqrt{\varrho}} \frac{\exp[-2\tau v^2 + \sum_{s \geq 3} O(v^s)]}{1+i(v+\sqrt{\varrho})^2} dv + \int_{2\sqrt{\varrho}}^\infty \frac{u^{2P} e^{-\tau u^2}}{1+iu^2} du. \tag{137}$$

Since $f^s(u) = (-1)^{s+1}(s-1)! \, 2P/u^s$, the series of terms denoted by $O(v^s) = (-1)^{s+1} 2Pv^s/s\varrho^{s/2}$ only converges provided $|v| < \sqrt{\varrho}$, restricting the limits of the first integral. Now around $v = 0$, the denominator in (137) is given by

$$\frac{1}{1+i(v+\sqrt{\varrho})^2} = \frac{1}{1+i\varrho}\exp\left[\frac{-2\sqrt{\varrho}}{(\varrho-i)}v + \frac{(\varrho^2-3)\varrho + i(3\varrho^2-1)}{(1+\varrho^2)^2}v^2 + \sum_{s \geq 3} O(v^s)\right], \tag{138}$$

where the terms denoted by $O(v^s) = (-1)^s 2v^s/s\varrho^{s/2}$ (real part to first order in $\varrho$). This means

$$I = \frac{\varrho^P e^{-P}}{1+i\varrho}\int_{-\sqrt{\varrho}}^{\sqrt{\varrho}} \exp\left[\frac{-2\sqrt{\varrho}}{(\varrho-i)}v - \left\{2\tau - \frac{(\varrho^2-3)\varrho + i(3\varrho^2-1)}{(1+\varrho^2)^2}\right\}v^2 + \sum_{s \geq 3} O(v^s)\right] dv + \int_{2\sqrt{\varrho}}^\infty \frac{u^{2P} e^{-\tau u^2}}{1+iu^2} du \tag{139}$$

with $O(v^s) = (-1)^s 2(1-P)v^s/s\varrho^{s/2}$ (real part to first order $\varrho$, $|v| < \sqrt{\varrho}$). Now consider the behaviour of the real part of the exponential phase of integrand (139), as $v \to \pm\sqrt{\varrho}$, in both the limits $\varrho \to 0$ and $\varrho \to \infty$.

$$\text{Re}\{\text{phase at } v \to \pm\sqrt{\varrho}\} \to \begin{cases} -2P - 2(P-1)\sum_{s \geq 3}(\mp 1)^s s^{-1} & \text{as } \varrho \to 0 \\[1em] \mp 2 - (2P-1) - 2(P-1)\sum_{s \geq 3}(\mp 1)^s s^{-1} & \text{as } \varrho \to \infty \end{cases} \tag{140}$$

Equation (140) shows that provided $P \geq 2$ the real part of the phase is already significantly negative as $v \to \pm\sqrt{\varrho}$, irrespective of the actual value of $\varrho = P/\tau$. Consequently the magnitude of second integral in (139) will be exponentially smaller than the first integral. Since $|1+iu^2|^{-1} < 1 \, \forall \, u > 2\sqrt{\varrho}$, it is relatively easy to show that the second integral in (139) is no more than $O\left((4\varrho)^{P-1/2} e^{-4P}/2\tau\right)$. This means that integral $I$ satisfies

$$I = \frac{\varrho^P e^{-P}}{1+\mathrm{i}\varrho}\left[\int_{-\sqrt{\varrho}}^{\sqrt{\varrho}} exp\left[\frac{-2\sqrt{\varrho}}{(\varrho-\mathrm{i})}v - A_2(\tau,\varrho)v^2 + \sum_{s\geq 3} O(v^s)\right]dv\right] + O\left(\frac{(4\varrho)^{P-1/2}e^{-4P}}{2\tau}\right),$$

$$\text{where } A_2(\tau,\varrho) = 2\tau - \frac{(\varrho^2-3)\varrho + \mathrm{i}(3\varrho^2-1)}{(1+\varrho^2)^2}. \tag{141}$$

In this form $I$ can be approximated by elementary functions using the tabulated standard integral [15 eq. 3.462.3]. There is one slight drawback because the standard integral [15 eq. 3.462.3] covers the infinite range $v \in (-\infty, \infty)$ rather than the finite range $v \in (-\sqrt{\varrho}, \sqrt{\varrho})$ in (141). So employing [15 eq. 3.462.3] introduces a further error of the form

$$\frac{\varrho^P e^{-P}}{1+\mathrm{i}\varrho}\left[\int_{\sqrt{\varrho}}^{\infty} v^s exp\left[\frac{-2\sqrt{\varrho}}{(\varrho-\mathrm{i})}v - A_2(\tau,\varrho)v^2\right]dv\right] \sim O\left(\frac{\varrho^P e^{-3P}\Gamma(s)}{(1+\mathrm{i}\varrho)(2\sqrt{\tau P})^s}\right). \tag{142}$$

The presence of the $\Gamma(s)$ term in (142) restricts the use of [15 eq. 3.462.3] to the range $s \leq 2\sqrt{\tau P}$. With this caveat one can establish estimates for $I$ to increasing degrees of precision. Terminating the computation on reaching the terms of $O(v^3)$ in (141) gives

$$I = \frac{\varrho^P e^{-P}}{1+\mathrm{i}\varrho}\sqrt{\frac{\pi}{A_2(\tau,\varrho)}}[C(\tau,\varrho)] + O\left(\frac{(4\varrho)^{P-1/2}e^{-4P}}{2\tau}\right) + O\left(\frac{\varrho^P e^{-3P}\Gamma(3)}{(1+\mathrm{i}\varrho)(2\sqrt{\tau P})^3}\right),$$

with $\quad C(\tau,\varrho) = exp(A_1(\tau,\varrho))[1 + A_3(\tau,\varrho)], \quad A_1(\tau,\varrho) = \varrho/[A_2(\tau,\varrho)(\varrho-\mathrm{i})^2],$

$$\text{and } A_3(\tau,\varrho) = -\frac{(2A_1(\tau,\varrho)/3+1)}{(A_2(\tau,\varrho))^2(\varrho-\mathrm{i})}\left\{\tau - \frac{\varrho[\varrho^4 - 12\varrho^2 + 3 + \mathrm{i}2\varrho(3\varrho^2 - 5)]}{(1+\varrho^2)^3}\right\}. \tag{143}$$

Further explicit expressions for the various $A_{s\geq 4}(\tau,\varrho)$, arising from the contributions of the $O(v^{s\geq 4})$ in (141), can be found by repeated application of [15 eq. 3.462.3]. However, these expressions become increasingly complicated, whilst producing progressively smaller increases in precision. Since $A_1(\tau,\varrho), A_{s\geq 3}(\tau,\varrho) \ll A_2(\tau,\varrho)$ irrespective of the limit $\varrho \to 0$ or $\varrho \to \infty$, a concise and accurate estimate for the generic integral $I$ is obtained by setting the constant $C(\tau,\varrho) = 1$. For this application $\tau \geq \pi/2$, and direct numerical computation of $I$ shows that estimate (143) is accurate up to a maximum relative error of $< 5\%$. Inclusion of both the $A_1(\tau,\varrho)$ and $A_3(\tau,\varrho)$ terms reduces this to $< 1\%$.

Using (143) it is now possible to obtain more refined estimates for the remainder terms (135). In these specific instances $\tau = |z_k^\pm|^2/2 = \pi|k \pm \xi|^2/x$ and $\varrho = Px/\pi|k \pm \xi|^2$. Since the value $k = |\xi| = M$ is excluded from these remainders (that erfc function is found separately in 64), $\tau \geq \pi/2$ and $\varrho \leq 2P/\pi$. With $C(\tau,\varrho) = 1$, this gives an improved estimate for the remainder terms of the form

$$R_P(z_k^\pm) \approx \frac{\text{sgn}(k \pm \xi)e^{i\pi/4}(-i)^P P^P e^{-P} x^{P+1/2}}{\pi^{P+1}|k \pm \xi|^{2P+1}(1 + iPx/\pi|k \pm \xi|^2)} \left\{ \left(1 - \frac{(\varrho^2 - 3)\varrho + i(3\varrho^2 - 1)}{2\tau(1 + \varrho^2)^2}\right)^{-\frac{1}{2}} \right\}. \quad (144)$$

To first $O(|k \pm \xi|^{-1})$ equation (144) recovers (as it should) the same upper bound on $|R_P(z_k^\pm)|$ as that given by [33 eq. 2.7], which is the basis of the error estimates (61-62). Now the *exact* expression for the total remainder term $\tilde{R}_P'$ in (58, 64) is given by (q.v. 61)

$$\tilde{R}_P' = R_P'(x, \xi) - R_P(x, \theta),$$

$$R_P'(x, \xi) = \frac{f(N)e^{i\pi/4}}{\sqrt{2\pi x}} \left[ \sum_{k=M+1}^{\infty} R_P(z_k^-) + \sum_{k=1}^{M-1} R_P(z_k^-) - \sum_{k=1}^{\infty} R_P(z_k^+) \right],$$

$$R_P(x, \theta) = \frac{e^{i\pi/4}}{\sqrt{2\pi x}} \left[ \sum_{k=1}^{\infty} R_P(z_k^-) - \sum_{k=1}^{\infty} R_P(z_k^+) \right]. \quad (145a, b, c)$$

Equation (145b) is the same as equation 2.10 of [33] but with the deletion of the $k = M$ term [33 c.f. discussion pp 585-6] and an extra $e^{i\pi/4}$ included for consistency with [33 eq. 2.4]. In equation (145c) the $R_P(z_k^\pm)$ terms are given by (144) but with $\theta \leftrightarrow \xi$. The only other feature of note is the inclusion of the extra unaccounted for $\sqrt{2/\pi}$ factor, a point previously raised after equation (61), which seems inconsistent with early equations of [33]. It is included in (145) for pragmatic reasons discussed in a moment.

The computation of $\tilde{R}_P'$ requires a means of estimating the infinite sums arising from (145) following the substitution of (144). For $P \geq 2$ these sums converge very rapidly and are dominated by the terms at which the factor $|k \pm \xi \leftrightarrow \theta|^{-2P-1}$ reaches its maximum. This maximum could become relatively large when $1/2 < |k \pm \xi \leftrightarrow \theta| < 3/2$, but this is only possible for the $R_P(z_{M\pm 1}^-)$ terms in (145b) and/or the $R_P(z_1^\pm)$ terms in (145c). (The $R_P(z_k^+)$ term in $R_P'(x, \xi)$ is always much smaller than its companions since $|k + \xi| \gtrsim O(K) \, \forall k$.) These (potentially) relatively large terms, given explicitly by

$$R_P'(x, \xi, M \pm 1) - R_P(x, \theta, 1) = \frac{e^{i\pi/4}}{\sqrt{2\pi x}} [f(N)\{R_P(z_{M+1}^-) + R_P(z_{M-1}^-)\} - \{R_P(z_1^-) - R_P(z_1^+)\}],$$
$$(146)$$

are much the most significant of all the various factors making up the remainder $\tilde{R}_P'$ not specifically accounted for in *algorithm QGS*. So incorporation of (146) directly into (64), using estimates for the various $R_P(z)$ computed from (144), should produce a refined version of *algorithm QGS* with significantly improved precision capabilities. The performance of this *refined algorithm QGS* applied to the five quadratic Gaussian sums (Cases 3.4A-E) previously examined in Section 3.4, is illustrated in Table 4.5.2. Only the final answers are shown in Table 4.5.2 (with $P = 3$ and $K = 20$ as before), along with the absolute error and relative error values. The corresponding values for these errors obtained previously using the basic version of *algorithm QGS* (see Tables 3.4A-E) are shown in parenthesises. As anticipated, inclusion of (146) gives a much higher degree of precision, with a more than five-fold reduction in the error terms in all cases (even examples 3.4D&E where the errors were already extremely small). The extra operational count needed to compute (146) using (144) is comparatively small (many of the parameters, such as $|k \pm \xi \leftrightarrow \theta|$, have already computed previously, whilst factors such as $e^{i\pi/4}(-i)^P P^P e^{-P}/\pi^{P+1}$ need only be computed once and then stored), roughly equivalent to adding an extra term in the short sum evaluation to $(P+1)\{4ops_{s;c} + 40ops_{ao}\}$ (see discussion before eq. 80). This means that for $P = 3$ the maximum operational count for one evaluation of $C(N, x, \theta, P)$ goes up from $53ops_{s;c}$ to about $60ops_{s;c}$. The resulting changes in the values of $\Omega^* = 0.59 \times 60$ and $\Omega^*_{max} = 1.44 \times 60ops_{s;c}$ are much too small to effect the conclusions of *Theorem 2*. The 'stray' $\sqrt{2/\pi}$ factor is incorporated into (146) simply because without it the precision performance of *algorithm QGS* is unchanged. So it should appear, although its origin remains a small mystery.

Following the incorporation of the most significant correction terms given by (146) directly into *algorithm QGS*, the general remainder term can now be written as

$$|\tilde{R}'_P| = \left| \frac{f(N)e^{i\pi/4}}{\sqrt{2\pi x}} \left[ \sum_{k=M+l+1}^{\infty} R_P(z_k^-) + \sum_{k=1}^{M-l-1} R_P(z_k^-) - \sum_{k=1}^{\infty} R_P(z_k^+) \right] - \frac{e^{i\pi/4}}{\sqrt{2\pi x}} \sum_{k=l}^{\infty} [R_P(z_k^-) - R_P(z_k^+)] \right|$$

$$\lesssim 4 \times \frac{P^P e^{-P} x^P}{\sqrt{2}\pi^{P+3/2} |l \pm \theta \leftrightarrow M + l \pm 1 \pm \xi|^{2P+1}} < \left(\frac{2}{\pi}\right)^{3/2} \frac{2}{|2l-1|} \left(\frac{2P}{\pi|2l-1|^2}\right)^P e^{-P}, \quad (147)$$

for $l \in \mathbb{N}$, with $l = 2$ in this specific case. This comes about because $|x| < 1/2$ and the $MIN|l \pm \theta \leftrightarrow M + l \pm 1 \pm \xi| < (2l-1)/2$, irrespective of the values of $\xi \leftrightarrow \theta$. This means that the relative error in each individual Gaussian sum calculated from *algorithm QGS* using termination constant $K$ should satisfy (see 75-76)

$$|\varepsilon_{GS}(K)| \ll \frac{3.41|\tilde{R}'_P|}{\sqrt{K}} < \frac{2\sqrt{3}}{\sqrt{K}|2l-1|} \left(\frac{2P}{\pi|2l-1|^2}\right)^P e^{-P}. \quad (148)$$

For the calculations shown in Table 4.5.2 with $K = 20$, $P = 3$ and $l = 2$, this gives a value of $\varepsilon_{GS}(K) = 1.2 \times 10^{-4}$. Comparison with the results shown in Table 4.5.2, reveal this bound to be comparable to, but somewhat *below*, the actual values $4.5 \times 10^{-4}$ and $2.3 \times 10^{-4}$ found for the cases 3.4A-B respectively. The reason for these larger than expected values lies in the various $O(\varrho^{P-1} e^{-P}/2\tau)$ terms disregarded if the approximation $C(\tau, \varrho) = 1$ is adopted to estimate (146). Explicit computation of the terms $A_1(\tau, \varrho)$ and $A_3(\tau, \varrho)$ that make up $C(\tau, \varrho)$ (see 143) provides a complicated (and operationally expensive) means of obtaining extra, but incremental, reductions in $\varepsilon_{GS}(K)$. Rather than continuing any further with this incremental approach, a better strategy would be to modify equation (64) to include more and more of the various erfc terms explicitly; each computed using the intrinsic erf(z) routines common to modern software packages. This provides a much more efficient means of reducing (145), and hence $\varepsilon_{GS}(K)$, to the desired precision necessary to satisfy the constraints of *Theorem 3*. However, implementation of such a strategy would come at the expense of increasing the operational count needed to compute $ZP(t)$. This increase is estimated below.

Currently (64) already contains the $\text{erfc}(\text{sgn}(-\varepsilon)\, z_M^-/\sqrt{2})$ term associated with $k = M$ explicitly, since the corresponding asymptotic expansions (132 & 133) are not very accurate for $|\varepsilon| = |\xi - M| < 1/2$ [33 pp. 585-6]. The proposed strategy is to gradually add in further such $\text{erfc}(z_k^-/\sqrt{2})$ terms (q.v. eq. 129) into (64), starting with $k = M + \text{sgn}(\varepsilon)$, then $k = M - \text{sgn}(\varepsilon)$, followed by $k = M + 2\text{sgn}(\varepsilon)$, $M - 2\text{sgn}(\varepsilon)$, etc. In conjunction with this process one would also augment (64) with the $\text{erfc}(z_k^\pm/\sqrt{2})$ terms associated with the variables $z_k^\pm = (k \pm \theta)\omega\sqrt{2\pi/x}$ for $k = \text{sgn}(\theta), -\text{sgn}(\theta), 2\text{sgn}(\theta), \ldots$ Whenever a $\text{erfc}(z_k^\pm/\sqrt{2})$ term is added into (64), one would need to modify the expressions (59-60) for $c_r(\theta)$ and $c_r'(\xi)$ by removing the corresponding factor of $k$ from the overall sum. These are relatively minor modifications to *algorithm QGS* and for the purposes of proving *Theorem 3* the underlying details of no great importance. What really matters is the speed of the reduction in $|\varepsilon_{GS}(K)|$, as each successive $\text{erfc}(z_k^\pm/\sqrt{2})$ term is incorporated directly into (64). Is it fast enough to prove *Theorem* 3, without increasing the operational count beyond the limits imposed by *Theorem* 2? To answer this, suppose all the $\text{erfc}(z_k^\pm/\sqrt{2})$ terms for $k = M \pm \text{sgn}(\varepsilon), \ldots, M \pm l \times \text{sgn}(\varepsilon)$; $k = \pm\text{sgn}(\theta), \ldots \pm l \times \text{sgn}(\theta)$; are computed explicitly. Then the general remainder term will be given by (147) and the corresponding relative error will be indeed be bounded by (148). Notice that as $l$ increases in (148), the value of $P$ can also be increased, further reducing $|\varepsilon_{GS}(K)|$. The task now is to choose $l$ and $P$ to guarantee that $|\varepsilon_{GS}(K)|$ does not produce an increase in the relative error of $ZP(t)$, beyond that of the current $O(\varepsilon_t)$ level already present due to the various *derivation* errors discussed earlier.

To achieve this, first prescribe $P = \lceil \pi|2l - 1|^2/2e \rceil$, in which case

$$|\varepsilon_{GS}(K)| < \frac{2\sqrt{3}}{\sqrt{K}|2l-1|} e^{-\pi|2l-1|^2/e}. \tag{149}$$

Now in $ZP(t)$ there are $O(Qe^{-7/3}\log(t)[t/\varepsilon_t]^{1/3})$ Gaussian sums to evaluate, most of whom will be of $O(\sqrt{M_t^-})$ in size, but each of these must be divided by their appropriate amplitude factor $(\alpha_E^2 - a^2)^{1/4}$. The value for $\sqrt{M_t^-}(\alpha_E^2 - a^2)^{-1/4} \sim [\varepsilon_t pc/t(pc - 1)]^{1/6}$ from (24, 40), which means that if, in the *extremely unlikely event* that the phases of all the respective $\varepsilon_{GS}(K)$ were to correlate perfectly, the sum of the relative errors could reach a maximum $\varepsilon_{GS}(K) \times \log(t)[t/\varepsilon_t]^{1/3} \times [\varepsilon_t/t]^{1/6} = \varepsilon_{GS}(K)\log(t)[t/\varepsilon_t]^{1/6}$. To be certain that this term is no larger than the current $O(\varepsilon_t)$ relative error scale, one requires

$$\frac{\log(t)}{\sqrt{K}}\left(\frac{t}{\varepsilon_t}\right)^{1/6} e^{-\pi|2l-1|^2/e} \approx \varepsilon_t \Rightarrow 2l - 1 \approx \sqrt{\frac{e}{\pi}\log\left\{\frac{t^{1/6}\log(t)}{\varepsilon_t^{7/6}\sqrt{K}}\right\}},$$

$$\Rightarrow P = \frac{1}{2}\log\left\{\frac{t^{1/6}\log(t)}{\varepsilon_t^{7/6}\sqrt{K}}\right\}. \tag{150}$$

So the number of extra $\text{erfc}(z_k^\pm/\sqrt{2})$ terms which must be incorporated directly into (64) to ensure the relative error in $ZP(t)$ cannot exceed $O(\varepsilon_t)$ is $\approx \sqrt{\log(t/\varepsilon_t^7 K^3)}/5$ as $t \to \infty$. At the same time the corresponding cut-off integer $P$ increases as the square of this value. Implementing this strategy leads to an increase in the operational count over and above that predicted in *Theorem 2*. Currently a constant value of $\Omega^* = 0.59 \times 53$ has been employed to estimate the average operational count of computing the iteration count factor $y \times C(N, x, \theta, P)$, based on a value of cut-off integer $P = 3$ with $l$ set (implicitly) equal to one (see the discussion at the end of Section 3.6). This is used in (118-119) to obtain the upper bound on the operational count necessary to evaluate $ZP(t)$. Based on the estimates of Section 3.6, $\Omega^*$ would now become a linear function of $P$ and $l$, of the form $\Omega^* \sim 0.59 \times (12 + 5P + 9\{2 + l\}) ops_{S;c}$. Asymptotically, this means that in order to evaluate $ZP(t)$ and be *absolutely certain* of achieving an approximation to the partial sum (121) of the last $N \in [N_{\alpha_{Emax}}, N_t]$ terms of $Z(t)$ which is accurate to $O(\varepsilon_t)$ in the relative error, the overall operational count would have to rise from $O((t/\varepsilon_t)^{1/3}[log(t)]^2)$ of *Theorem 2* (assuming $K = O(10)$) to $O\left((t/\varepsilon_t)^{1/3}[log(t)]^2 log(t/\varepsilon_t^7 K^3)\right)$. This relatively small increase in the operational count guarantees the desired level of accuracy specified by *Theorem 3*. However, for actual computations (see Section 5) the implementation such a complex error refinement strategy is a proviso of only theoretical, as opposed to a practical, necessity.

Since both sources of error that arise when linking $ZP(t)$ to $Z(t)$ (from the *derivation* of the formulation of $ZP(t)$ and from its *computation* via *algorithm QGS*) have, with the refinements discussed above, been shown to give rise to corrections no larger than $O(\varepsilon_t)$, this completes the proof of *Theorem 3*.

*4.6 Corollary of Theorem 3*

This follows directly from the results of *Theorems 2 & 3* and the connection formula (15), linking terms of the RSF to the terms of the new (approximation) series (11) for $Z(t)$. From (8) and (121) one has

$$Z(t) = 2 \sum_{N=1}^{N_t} \frac{cos(\theta(t) - tlog(N))}{\sqrt{N}} + O\left(\frac{1}{t^{1/4}}\right) = 2 \sum_{N=1}^{N_c-1} \frac{cos(\theta(t) - tlog(N))}{\sqrt{N}} + ZP(t)[1 + O(\varepsilon_t)],$$

$$= Z(N_c, t) + ZP(t)[1 + O(\varepsilon_t)], \qquad (151)$$

with $N_c = \lfloor \alpha_E^c \{1 - \sqrt{1 - (a/\alpha_E^c)^2}\}/4 \rfloor$ and $\alpha_E^c = INT_E\left[(\varepsilon_t t^2)^{1/3}/\left(\sqrt{\pi} log(t)\right)\right]$. Clearly calculation of first sum on the right hand side of (151) requires just $(N_c - 1) ops_{S;c}$. Following the computation methodology of sections 4.3.3 & 4, *Theorem 2* proves that the calculation of second sum $ZP(t)$ requires no more than $(t/\varepsilon_t)^{1/3} log(t)[O(log(t)) + O(K)] ops_{S;c}$. In turn *Theorem 3* demonstrates (subject to the proviso discussed at the end of section 4.5.2b) that computing $ZP(t)$ by means of this methodology, yields an approximation to the sum of the latter $N \in [N_c, N_t]$ terms of $Z(t)$, which is accurate to $O(\varepsilon_t)$ in the relative error. Since the cut-off value $N_c \approx t^{1/3} log(t)/\varepsilon_t^{1/3} \sqrt{\pi}$ from (122), this means that is possible to compute an approximation for $Z(t)$ in just

$$\left(\frac{t}{\varepsilon_t}\right)^{1/3} \log(t) \left[\{O(\log(t)) + O(K)\}\log(t/\varepsilon_t^7 K^3) + \frac{1}{\sqrt{\pi}}\right] ops_{s;c}, \tag{152}$$

accurate to $|O(\varepsilon_t) \times ZP(t)|$ in the limit as $t \to \infty$. With the definitions for $K$ and $\varepsilon_t$ specified in the corollary, the final result follows automatically.

One is tempted to say that the precision specified is accurate to $O(\varepsilon_t)$ in the absolute relative error of $Z(t)$ itself. However, the analysis only supports an accuracy in the absolute relative error of $ZP(t)$. To say these two statements are equivalent requires the condition that $O(Z(t)) = O(ZP(t))$ to be satisfied. Employing the methodology of classical exponent pairs one can find a bound for the exponent $c_1 = 11081719/78340470 = 0.14145586\ldots$ in the estimate $Z(N_c, t) = O(|t|^{c_1+\varepsilon})$, less than the well-known Huxley bound $c = 32/205 = 0.15609\ldots$ for the estimate of $Z(t) = O(|t|^{c+\varepsilon})$ itself [21]. (The bound on $c_1$ is easily established from the exponent pair which results from the operations $A^2(BA^2)^2(BA^3)^4BA(\kappa, \lambda)$ starting from an initial exponent pair $(\kappa, \lambda) = (2/7, 4/7)$, The $A, B -$ processes are defined in [22 Ch. 2.]) So it would *appear* that $Z(N_c, t) = o(ZP(t))$, which would imply that $O(Z(t)) = O(ZP(t))$, as one would expect intuitively. In practical terms if one wishes to estimate for $Z(t)$ by computing $ZP(t)$ instead, then the calculation will be accurate to $|O(\varepsilon_t) \times ZP(t)|$.

## 5. Sample Computations

The computational methodology discussed in Section 4.3 in the proof of *Theorem 2* can be adapted to generate a suitable algorithm for practical calculations of $Z(t)$. However, the benefits (in terms of a significantly reduced operational count compared to the RSF) of estimating $Z(t)$ by computing $ZP(t)$ and then employing (151) will only start to manifest themselves for 'relatively' large $t$ values. The analysis leading to (119) shows that the computation of $ZP(t)$ requires $(t/\varepsilon_t)^{1/3}(\log(t))^{2+o(1)} ops_{s;c}$, compared to the $O(\sqrt{t})ops_{s;c}$ needed for (8). Hence significant savings should start to accrue beyond the value $t \sim 10^{20}$, when the latter operational count overtakes the former. To investigate further how this comes about, some sample calculations of $Z(t)$ for $t \in [10^{18} - 10^{23}]$ were carried out using the new methodology and compared to the results obtained using the RSF.

*5.1 A new algorithm for the computation of $Z(t)$ in just $O\left((t/\varepsilon_t)^{1/3}(\log(t))^{2+o(1)}\right)$ operations (algorithm ZT13) assuming storage capacity $O\left((\log(t))^2\right)$ bits.*

*Initialise $t > 10^{15}$, $0 < \varepsilon_t \ll 1$, $Y \in [1.05, 1.15]$, $ZP := 0$ and integer values for $P \in [2, 4]$ and $K \in [30, 80]$.*
*Compute $a = \sqrt{8t/\pi}$, $INT_O(a) + 2$, $ib = 2 * \lfloor t^{1/4} \rfloor$, $\alpha_E^c := INT_E\left[(\varepsilon_t t^2)^{1/3}/(\sqrt{\pi}\log(t))\right]$, $h_t := \log(\varepsilon_t)/\log(t)$ and $X := \exp[1/12 - h_t/3]$.*
*Initialise Block No $p := 0$. If $\{NINT_O(a) - a\} \in [\pm t^{-1/6},]$, then compute the transition term $T_{a+\varepsilon}$ (see eq. 13) and add to ZP. For summands $\alpha = [INT_O(a) + 2, INT_O(a) + 2 + ib]$ step 2,*

*compute first sum contribution to ZP(t) (see eq. 82) directly and add to ZP. Initialise first pivot* $\alpha_E := INT_O(a) + 2 + ib + 1$.

*while* $\alpha_E < \alpha_E^c$  *do Block No. p = 1 …*

  *if* $\alpha_E < Ya$ *then*

$$NQGS := \lfloor X^p * ib/2 \rfloor, M_t := NINT_O\left[\left\{\frac{\varepsilon_t a}{\pi}\left(\frac{\alpha_E}{a} - 1\right)\right\}^{\frac{1}{3}}\right],$$

  *else*

$$M_t := NINT_O\left[\frac{(\varepsilon_t^2 t)^{1/6}}{\sqrt{2\pi}}\left(\frac{2\alpha_E}{a} - \frac{a}{2\alpha_E}\right)\right],$$

  *endif*

  $M_t^- := (M_t - 1)/2, \quad \alpha_E := \alpha_E + M_t + 1,$

  *for* $j := 1$ *to NQGS do*

    *compute GS parameters* $x, \theta^\pm, \omega^\pm$ *and* $(\alpha_E^2 - a^2)^{1/4}$, *(see eq. 83)*

    *if* $M_t^- < K$ *then*

      *compute* $S_{M_t^-}^*(x, \theta^\pm)$ *directly*

    *else*

      *compute* $S_{M_t^-}^*(x, \theta^\pm)$ *using the refined version of algorithm QGS*

    *endif*

    $ZP := ZP + next\ term\ of\ second\ sum\ of\ ZP(t)$ *(see eq. 82)*

    $\alpha_E := \alpha_E + 2 * (M_t + 1)$

    *if* $\alpha_E > \alpha_E^c \quad \alpha_E^c := \alpha_E + M_t + 1 \quad exit\ p\ loop\ \ endif$

  *enddo (j loop)*

  $\alpha_E := \alpha_E - (M_t + 1)$

  *enddo (p loop)*

$N_c := \lfloor \alpha_E^c\{1 - \sqrt{1 - (\alpha/\alpha_E^c)^2}\}/4 \rfloor$ *(see eq. 122)*

  *Compute* $Z(N_c, t)$, *the cut-off Riemann-Siegel sum (see eq. 151)*

*Estimate for* $Z(t) := Z(N_c, t) + ZP \quad end.$

The computational structure of *algorithm ZT13* is based on the methodology discussed in detail in Section 4 with a few minor alterations. In practice it is difficult to predict exactly the integer values $L$ at which point the step up procedure can be safely initiated (see section 4.3.2) and $L_c$ (estimate 114), the number of blocks needed for the calculation to reach $\alpha_E^c$. In section 4.3.2 the choice of $L = \lfloor log(t) \rfloor - 28$ ensured $Y > 1.24$ as $t \to \infty$ (see 97). However, waiting for $\alpha_E > 1.24a$ before instigating the step up to larger $M_t(\alpha_E)$ values is somewhat inefficient for practical calculations. Instead it is better to initialise a value for the step up parameter $Y \in [1.05, 1.15]$. Fixing $Y$ to a slightly smaller value speeds up the algorithm (since the switch over to the larger $M_t(\alpha_E)$ values comes earlier in the calculation) with no detrimental effect on its accuracy (see section 5.2). The other point to note is that the computation of the two Gaussian sums $S_{MT^-}^*(x, \theta^\pm)$ in step 4) *does not require two separate calls to algorithm QGS*. Since the quadratic parameter $x$ in each Gaussian sum is identical, certain replicative computations can be eliminated. Utilising this fact it is easy to create a slightly modified version of *algorithm QGS* which computes both $S_{MT^-}^*(x, \theta^+)$ and $S_{MT^-}^*(x, \theta^-)$ simultaneously, at an operational cost of only about one and three-quarters times that of a single call to the standard version. The refinements to the basic version of *algorithm QGS* discussed in Section 3 go no further than the incorporation of the improved error terms given by (146) into (64).

*5.2 Computational performance and accuracy of algorithm ZT13*

In order to test out some the theoretical predictions made in the previous three sections, a Fortran code implementation of *algorithm ZT13* was developed to carry out some sample calculations of $Z(t)$ in the range $t \in [10^{18-23}]$. Across these calculations, fixed default settings of $\varepsilon_0 = 2/log(10^{18}) = 0.0482$ ..., for the relative error $\varepsilon_t$, $K = 30$ for the termination constant, cut-off integer $P = 3$ and step-up parameter $Y = 10/9$, were employed. Table 5.1 gives a detailed breakdown of all the intermediate calculations necessary to obtain an estimate for $Z(t = 10^{18})$, to give an insight into what is involved. The table breaks down each part of the overall calculation into the contributions given by the various blocks. The block size $\lfloor X^p t^{1/4} \rfloor$, which gives the number of pivots within the block, is shown first (here $X = 1.1137$ ... as defined in section 4.3.3), followed by the associated collection size $M_t$. The calculations encompassed by the $p = 0$ block, comprise the sum of those terms in (82) too large and distinctive to be effectively parcelled together into quadratic Gaussian sums (*step 3 in ZT13*). The size of this block was set to $2 \times INT_E(t^{1/4})$, double that suggested in section 4.3.1. Doubling the initial block size allows for significant parcelling to commence immediately from block $p = 1$, finessing the prescription set out in section 4.3.3 which postulated the parcelling should begin in pair-wise fashion with $M_{t,1} = 1$. In this example, the prescribed value of $\varepsilon_t$ means one can safely set $M_{t,1} = 9$, speeding up the computations a touch. The third column gives the range of integers encompassed by each block. So for instance, block 1 consists of all the (odd) integers $\alpha \in [1595832369, 1596536727]$, whilst the even integer $\alpha_{E,0} = 1595832368$ forms the boundary of blocks 0 and 1. Each set of Gaussian sum parameters in (82) is computed from the pivots $\alpha_E = \alpha_{E,0} + 10, \alpha_{E,0} + 30, \alpha_{E,0} + 50, ...$

that lie within block 1. This conforms to the prescription set out before (86-87), but with $M_{t,1} = 9$ rather than $M_{t,1} = 1$. The fourth column shows the total value of the partial sum over all the $\alpha$ integers for each block. The exact value of partial sum, derived from the corresponding terms of the RSF (summing backwards of course), appears in column six, followed by the relative error.

In the early blocks the quadratic Gaussian sums are relatively short and are computed exactly. However, upon reaching the sixth block one finds $M_{t,6} = 61 > 2K$, which means that the quadratic Gaussian sums are now sufficiently large to trigger the implementation of *algorithm QGS* to facilitate their computation (*step 4* in ZT13). The strict restrictions imposed on the size of the successive $M_{t,p}$, means that although these computations are no longer exact, the relative errors associated with the partial sums of the blocks are no larger than $\varepsilon_t$. By the time the computation reaches the end of block $p = 14$, the pivot values $\alpha_E > Ya \approx 1771303725$. At this stage it is safe to step-up $M_{t,p}$ to the larger of the two scales proposed in (83b). Hence the sharp jump in the designated $M_{t,p}$ sizes seen between blocks $p = 14$ & 15 in the tabulated results. The step-up occurs close to the asymptotically predicted block value of $L = \lfloor log(t) \rfloor - 28 = 13$. Beyond this point of the calculation all the block sizes remain fixed at $\lfloor X^{14} t^{1/4} \rfloor = 142854$ value reached by block $p = 14$. On reaching block 26 one encounters an anomaly in the form of a relative error value of 0.0850, somewhat larger than the upper bound $\varepsilon_t$. The cause of this anomaly is down to fact that magnitude of the partial sum itself just happens to lie close to zero, rather than any fault with the methodology. The corresponding absolute error ($5.24 \times 10^{-4}$) for this block is of the same order as the absolute errors found in the partial sums of the other blocks. Such anomalous behaviour will only be problematic if the value of $t$ is such that $Z(t) \approx 0$, when any discrepancy would give rise to a large relative error. In such instances one would need to repeat the calculation to greater precision by refining the choices of $\varepsilon_t$ and $K$ to establish estimates for the absolute error. From block 27 onwards the calculation continues in the prescribed manner until one reaches block 37. Close to the beginning of this block the pivots $\alpha_E$ exceed the cut off value $\alpha_E^c = INT_E[(\varepsilon_t t^2)^{1/3}/\sqrt{\pi} log(t)] = 4955842210$, which terminates the procedure. In practice termination always occurs after a complete Gaussian sum computation, which means that the final integer of the final block ($\alpha_0^c = 49558423767$ here) will always fractionally exceed the value of $\alpha_E^c$ (and will be odd because $M_{t,p}$ is odd). In this instance, termination at the 37[th] block corresponds to a value of $L_c = 37 - L = 23 \approx \left(\frac{1}{2}\right) log(t)$, somewhat smaller than the estimate (114) in the limit $t \to \infty$. As $t$ increases one would expect the $Q \approx \left(\frac{1}{2}\right)$ factor to move a little closer to the predicted upper bound (114) at around $Q \approx 2$.

Combining all the separate partial sums from each of the 37 blocks gives a value for $ZP(10^{18}) = -0.376110$. Since $\alpha_E^c = 4955843767$ corresponds to a cut-off integer $N_c = 65986402$ using (15 & 122), $ZP(10^{18})$ approximates the sum of the final $[N_c, N_t = 398942280]$ main terms of the RSF, which equals 0.375570 ... to six decimal places. The remaining $[1, N_c - 1]$ terms of the formula, plus the small scale $O(t^{-1/4})$ remainders specified in (8) amounts to 0.565274. Adding this to $ZP(10^{18})$ gives an estimate for $Z(10^{18})$ of 0.189164, compared to the actual value $Z(10^{18}) = 0.189704$ ... (to six d. p.). This equates to an overall relative (absolute) error of $2.85 \times 10^{-3}$ ($5.40 \times 10^{-4}$), an order of magnitude below prescribed valuation of $\varepsilon_t = 0.0482$. Indeed the smallness of this overall error is indicative that there is no sign of any correlations which might cause the computational errors associated with the multiple calls of *algorithm QGS* to aggregate detrimentally necessitating the implementation of the proviso discussed in section 4.5.2b.

Table 5.2 presents a much briefer summary of the performance of *algorithm ZT13* for values of $t = 10^{18,19,20,21,22\ \&\ 23}$, compared to the corresponding RSF calculations. All the relative errors lie an order of magnitude below the prescribed valuation of $\varepsilon_t = \varepsilon_0$, with agreement to about three significant figures. In the $t = 10^{23}$ calculation, *algorithm QGS* is utilised in all but blocks zero and one (i.e. $p^* = 1$), $M_{t,p}$ is stepped up after block number $L = 22$ and the computation terminates near the beginning of block 91. So in this case $L_c = 69 \equiv Qlog(t)$ corresponding to an increase in $Q \approx 1.3$, as anticipated. Across the six calculations, average parameter values for $\bar{Q} \approx 0.94, \overline{p^*} = 3, \bar{L} = 18\ \&\ \bar{L}_c = 45$ pertained. Table 5.3 gives an illustration the kind of accuracy attainable utilising *algorithm ZT13* to calculate $Z(10^{20})$ for three different values of $\varepsilon_t$. All three estimates fall well inside their prescribed error bounds, with the most accurate estimate pertaining for the smallest value of $\varepsilon_t$ as expected. This higher accuracy is achieved at the expense of a longer CPU run time. Changes to the settings of the termination constant $K$ and the step-up parameter $Y$ produce only minor changes to the accuracy and run times.

Table 5.2 also shows the (sequential) CPU time necessary to complete the $t = 10^{18-23}$ calculations. As one can see *algorithm ZT13* starts to overtake the RSF in terms of computational speed from $t \approx 10^{19}$ onwards. However, a better reflection of the performance of *algorithm ZT13*, in terms of the operational count framework utilised in Section 4, is obtained by considering the respective timing ratios. Fig. 2a shows the ratio of the timings of *algorithm ZT13* compared to the corresponding RSF calculation for $t^n$ with $n = 18 - 23$. An estimate for the upper bound on the operational count for the former was derived sections 4.3.1-6. To second order the bound (118) can be expressed as

$$< \left(\frac{t}{\varepsilon_t}\right)^{\frac{1}{3}} [log(t)]^2 \left\{A + \frac{B}{log(t)}\right\} ops_{s;c}, \qquad (153a)$$

with

$$A = \frac{\Omega^* Q}{\chi}\left\{\frac{1}{3} + Qlog(1 + \chi^{-1})\right\}, \qquad \chi = e^{(7-28h_t)/3},$$

$$B = \frac{\Omega^* Q}{\chi}\left\{\frac{4(\Omega + K)}{\Omega^*} + log\left(\frac{2\varepsilon_t^{2/3} Y^2}{\pi K^4 (1 + \chi^{-1})}\right) + \frac{e^{1/12 - h_t/3}}{Q(X - 1)}\left(B_0 log(X) + \frac{B_1}{4}\right)\right\}. \qquad (153b)$$

The constants $A$ and $B$ can be estimated from the prescribed values of $\varepsilon_t, K\ \&\ Y$ and the averages of $\bar{Q}, \overline{p^*}\ \&\ \bar{L}$ found above, provided the operational conversion factors $\Omega$ and $\Omega^*$ are known. Estimates from [27] discussed in section 4.2 show that $\Omega \approx 1.17 - 1.30$, whilst in section 3.6 it was established that $\Omega^* = (0.59 - 1.44) \times 53$ counts, using an efficient intrinsic erf routine. Equating $\bar{\Omega} = 1.24$ and $\overline{\Omega^*} = 1.01 \times 53$, one obtains the following mean values $\bar{A} = 1.012$ and $\bar{B} = -28.744$ for the constants. Dividing (153a) by $\sqrt{t/2\pi}\ ops_{s;c}$, the operational count of the RSF, gives an (almost) equivalent theoretical prediction for the timings ratio found from the data in Table 5.2. For the equivalence to be precise this prediction must include an extra parameter $\lambda \in (0, 1)$, to compensate for the fact that (153a) is an *upper bound* (rather than a direct estimate) on the operational count of *algorithm ZT13*. Fig. 2a shows a plot the theoretical ratio derived from (153) against the actual computational data of Table 5.2, with a fitted $\lambda = 0.312$ to minimize least squares errors. As one can the agreement between prediction and computations is excellent. The fitted value of $\lambda = 0.312$ suggests that (118,153) overestimates the operational count of *algorithm ZT13* by a factor of about 3.2 across this particular range of $t$.

Another way of comparing the performance of *algorithm ZT13* vis-a-vis the RSF is to compare the changes in the computational timings ratio $Z(10^n):Z(10^{n-1})$ as $n$ increases. These results are shown in Fig. 2b. In the case of (8), the ratio remains constant at $\sqrt{10} \approx 3.16$ for all $n$, as shown. By contrast for *algorithm ZT13* it *declines* as $n$ increases. The asymptotic prediction derived from (153a) is represented by the lower green line on the figure, whilst the upper red line shows the ratio's predicted behaviour based upon the values $\bar{A} = 1.012$ and $\bar{B} = -28.744$ established earlier. As can be seen, the corresponding data points for the actual calculations shown in Table 5.2 fall increasingly close to this latter result. This ratio analysis provides a somewhat better test of the predicted bounds (118, 153), since it is independent of the choice of $\lambda$ which cancels out. One concludes that in this computational range, increasing $t$ by a factor of ten requires only $\approx 2.6$ increase in the operational count using *algorithm ZT13*, compared to the $\approx 3.16$ required by the RSF. Asymptotically the ratio is predicted to fall to around $10^{1/3} \approx 2.15$, but this would require $t$ values well beyond anything currently computational feasible.

One final interesting feature of the performance of *algorithm ZT13* is the evolution of the computational time required to complete an individual block of calculations. Fig. 3, which pertains for the $Z(10^{23})$ calculation, is typical. As one can see from the figure, the longest computational times are required for those blocks numbered $p \approx L$ (marked by the vertical dashed line) close to where the step up of $M_{t,p}$ takes place. This is somewhat surprising since in the later blocks $M_{t,p}$ is much larger, and one might expect the associated quadratic Gaussian sums would take longer to estimate. However, this does not occur because the rate of decrease in the quadratic parameter $x = 1/(pc - 1)$ is faster than the corresponding increase in $M_{t,p}$, and by the time the cut-off point is reached $M_{t,p}x = O(1)$. This means that near the cut-off point *algorithm QGS* only requires a *single*, as opposed to $\sim log(M_{t,p}^-/K)$ iterations, to estimate each $S_{M_{t,p}^-}^*(x, \theta^\pm)$ sum. Hence the computations for these later blocks are somewhat faster than those for the intermediate numbered blocks. This is tied in with the fact that beyond the cut-off point the Gaussian sums start to lose their quadratic character and behave like ordinary geometric series. This is a distinctive feature of this particular representation of $Z(t)$ in terms of Gaussian sums, quite unlike the representation presented in [18].

All the detailed calculations presented here were carried out sequentially, using a single processor, in order to verify the operational estimates made in *Theorem 2*. However, very much faster computational speeds can be achieved by means of parallelisation. A parallel coding version of *algorithm ZT13* using the facilities of the *ARCHER UK National Supercomputing Service* (http://www.archer.ac.uk), is currently under development.

## 6 Conclusions

The paper [18] conclusively showed that it should be possible reduce the computational complexity of the RSF main sum down from $O_\varepsilon(t^{1/2-o(1)})$ proposal of [6] down to only $O(t^{1/3}(log(t))^\kappa)$ (where $\kappa$ is an absolute constant), by reformulating $Z(t)$ in terms generalised quadratic Gaussian sums typically of $O(t^{1/6})$ in length. The (approximate) quadratic reciprocity property (exemplified by eq. 58) of the latter, which means they can be estimated recursively using $\sim O(log(t))$ operations, provides the basis for the computational reduction. The algorithmic methodology of [18] was formulated directly from the terms of the RSF main sum. This paper demonstrates that the same kind of methodology can also be utilised to calculate $Z(t)$ based upon the approximate formulation given by (11 & 14) derived by [27 under the conditions set out in *Theorem 1*] from the sum of Kummer's functions (10). The analysis carried out in Section 2 shows that from this starting point $Z(t)$ can be rewritten in terms sub-sequences of generalised quadratic Gaussian sums (54), accurate to within a fixed bound on the relative error. Although structurally analogous to representation derived by [18], the detailed specification is distinctive and more subtle.

The remainder of this work is devoted in the development of a computational algorithm for the estimation $Z(t)$, utilising approximation (54). The first task is to formulate a practical methodology for the computation of a quadratic Gaussian sum of length $N$ at a logarithmical operational expense. The result is *algorithm QGS* of Section 3, which utilises the inherent quadratic reciprocity property [3] of Gaussian sums, to make such a computation estimate at a cost of just $O(\{2K + \Omega^* log(N/K)\})$ operations, where $K = O(100)$ termination constant. The work of [33] is crucial in this regard, both for the asymptotic quadratic reciprocity approximation (58) encapsulated by *Theorem [Paris]* and the initial error analysis which leads to error bounds formulated in section 3.5. Section 4 is devoted to two main results. In the first half is devoted to the sum $ZP(t)$ and its computation in no more than $(t/\varepsilon_t)^{1/3} log(t)\{O(log(t)) + O(K)\}$ operations (for $0 < \varepsilon_t \ll 1, K \in [10, 100])$, using *algorithm QGS* and the methodology established in the proof of *Theorem 2*. The second half is devoted to the accuracy of the approximation provided by $ZP(t)$ to a partial sum of the RSF for $Z(t)$. With some augmentations to the error analysis of [33], which leads to the formulation a *refined algorithm QGS*, one can prove that this approximation will be accurate to $O(\varepsilon_t)$ in the relative error in the limit as $t \to \infty$. A proviso to the analysis underlying the proof of *Theorem 3* specifies that in order to achieve this level of accuracy as $t \to \infty$, the operational count for the computation of $ZP(t)$ needs to rise, theoretically, by factor of $log(t/\varepsilon_t^7 \sqrt{K})$. However, the degree of correlation necessary, in what are essentially vast sequences of quasi-random error terms, to make such a refinement a necessity is so remote as to render it superfluous for all practical large scale calculations. The results of Sections 3 & 4 facilitate the development of the new *algorithm ZT13* (section 5.1), which produces an evaluation of the Hardy function in just $O\left((t/\varepsilon_t)^{1/3}(log(t))^{2+o(1)}\right)$ operations (excluding the proviso above). The detailed analysis of the behaviour of *algorithm QGS* enable one to make a much more specific estimate for the value of $\kappa$ than given by [18]. Sample computations utilising *algorithm ZT13* for values in the range $t = 10^{18-23}$ illustrate it behaves very much as predicted. Comparisons with exact computations, demonstrate that *algorithm ZT13* leads to estimates for $Z(t)$ which are well within the specified maximum relative error $\varepsilon_t$ (see Tables 5.2 & 3) tolerance. In fact for these specific examples the accuracy is an order of magnitude more precise. The speed of the computations conforms to the operational prediction of *Theorem 2* (Fig. 2b) and is faster than the RSF beyond $t \geq 10^{19}$, as shown in Fig. 2a.

It is interesting to speculate on future developments of this work. As mentioned in the Introductory remarks, [18, *Theorem 1.1*] formulates a procedure for the computation $Z(t)$ in just $O\left(t^{4/13}(log(t))^{\kappa}\right)$ operations in which "it is likely" that $\kappa$ "can be taken around 4". This formulation substitutes collections of $O(t^{5/26})$ cubic in place of the $O(t^{1/6})$ quadratic Gaussian sums, analogous to those utilised here. The former can then be estimated in much the same way as the latter, provided the cubic coefficient is sufficiently small. However, the implementation of any such a cubic procedure is problematic, since its computational speed would remain inferior to any $O\left(t^{1/3}(log(t))^3\right)$ quadratic procedure until $t \geq 10^{90}$, a figure well beyond anything currently practical. But this may be too pessimistic. In the methodology presented here, the analogous cubic Gaussian sum terms appear explicitly in (38), and the upper bound on the collection size (39) was based upon ensuring these terms remain sufficiently small. Some preliminary examination of the effects of including these terms directly, suggests that one might be able to formulate a representation for $Z(t)$ in which the cubic sums can take collection sizes as large as $O(t^{1/4})$ and still remain computable in just $O(log(t))$ operations, using analogous methods to those of *algorithm QGS*. Such a result would imply that one could compute $Z(t)$ using as little as $O\left(t^{1/4}(log(t))^{\kappa_1}\varepsilon_t^{-\kappa_2}\right)$ operations, accurate to some fixed *relative error* scale $\varepsilon_t$. One suspects that further advancements along these lines should prove possible in future.

## Tables

| $n$ | $L_n$ | $x_n$ | $s_n$ | $\theta_n$ |
|---|---|---|---|---|
| 1 | 129901233 | $1/\sqrt{45}$ | $+1$ | $1-\sqrt{(23/71)}$ |
| 2 | 19364532 | $7-\sqrt{45}$ | $+1$ | 0.390159303539095 |
| 3 | 5650494 | $(5-\sqrt{45})/4$ | $-1$ | 0.162904175231039 |
| 4 | 2413049 | $1-\sqrt{45}/5$ | $-1$ | $-0.381463061012124$ |
| 5 | 824395 | $(7-\sqrt{45})/4$ | $+1$ | 0.383438172238648 |
| 6 | 60139 | $x_2$ | $+1$ | 0.256248660552216 |
| 7 | 17548 | $x_3$ | $-1$ | $-0.378177224069896$ |
| 8 | 7493 | $x_4$ | $-1$ | $-0.114444787592628$ |
| 9 | 2559 | $x_5$ | $+1$ | 0.165014271963481 |
| 10 | 186 | $x_2$ | $+1$ | 0.262049291848319 |
| 11 | 54 | $x_3$ | $-1$ | $-0.398056283255953$ |
| 12 | 22 | $x_4$ | $-1$ | $-0.067895171805313$ |
| $n$ | Estimate of $S_{L_n}^{s_n}(x_n,\theta_n)$ | Exact Value of $S_{L_n}^{s_n}(x_n,\theta_n)$ | $|error|$ | Relative error |
| 12 | As exact value | $0.95384635 - 1.66611062i$ | 0.0 | 0.0 |
| 11 | $1.00607535 - 3.90371044i$ | $1.00333214 - 3.89043286i$ | 0.0135 | $3.374 \times 10^{-3}$ |
| 10 | $1.28630419 + 6.88855557i$ | $1.29322876 + 6.87130724i$ | 0.0185 | $2.658 \times 10^{-3}$ |
| 9 | $14.0775845 + 25.30169184i$ | $14.0772259 + 25.2328684i$ | 0.0688 | $2.382 \times 10^{-3}$ |
| 8 | $-6.12662881 - 48.9362647i$ | $-6.05515130 - 48.8420622i$ | 0.1182 | $2.402 \times 10^{-3}$ |
| 7 | $-27.6518893 + 69.5560559i$ | $-27.5007943 + 69.3982979i$ | 0.2184 | $2.926 \times 10^{-3}$ |
| 6 | $-62.3228952 + 123.849884i$ | $-62.0122200 + 123.581538i$ | 0.4105 | $2.969 \times 10^{-3}$ |
| 5 | $-479.859544 + 185.654054i$ | $-478.340786 + 185.691189i$ | 1.5192 | $2.961 \times 10^{-3}$ |
| 4 | $-531.440375 - 701.485722i$ | $-529.212066 - 700.177406i$ | 2.5840 | $2.944 \times 10^{-3}$ |
| 3 | $-78.9762197 + 1344.14844i$ | $-77.2560754 + 1340.59086i$ | 3.9516 | $2.943 \times 10^{-3}$ |
| 2 | $1777.43497 + 1746.87547i$ | $1774.57271 + 1740.14359i$ | 7.3151 | $2.943 \times 10^{-3}$ |
| 1 | $-4535.00190 - 4594.68373i$ | $-4527.85134 - 4577.13867i$ | 18.946 | $2.943 \times 10^{-3}$ |

Table 3.4A: Output generated by *algorithm QGS*, illustrating in detail the computation of the quadratic Gaussian sum $S_{N_0}(1/\sqrt{45}, 1-\sqrt{(23/71)})$. A set of recursive iterations are performed in which the original sum is rewritten in terms of a series of shorter, intermediate, quadratic Gauss sums, each with their own parameters. When $L_{n+1} < K$ the exact value of the $n$th (much shorter) sum is computed. From this starting point, the previous recursive iterations are reversed and an estimate of the original sum built up from equations (63-64). The exact values of the intermediate sums and the associated errors are also shown.

| $n$ | $L_n$ | $x_n$ | $s_n$ | $\theta_n$ |
|---|---|---|---|---|
| 1 | 129901233 | $1 - e/\pi$ | $+1$ | $1/e$ |
| 2 | 17503414 | $-0.42147960099874$ | $-1$ | $-0.230209768280677$ |
| 3 | 7377331 | $-0.37259406536018$ | $-1$ | $-0.453805669989322$ |
| 4 | 2748749 | $0.31611398846808$ | $+1$ | $0.282037310361788$ |
| 5 | 868918 | $-0.16341584517053$ | $-1$ | $-0.392201296527760$ |
| 6 | 141994 | $-0.11935763607527$ | $-1$ | $0.400019998785776$ |
| 7 | 16948 | $-0.37818201568035$ | $-1$ | $-0.351440359739468$ |
| 8 | 6409 | $0.35577061166962$ | $+1$ | $-0.429288927468500$ |
| 9 | 2279 | $0.18920009916778$ | $+1$ | $0.293354725243173$ |
| 10 | 431 | $-0.28540949184794$ | $-1$ | $-0.050499849278713$ |
| 11 | 122 | $0.49626228782619$ | $+1$ | $-0.176938226376921$ |
| 12 | 60 | $-0.01506345440907$ | $-1$ | $0.356541753660095$ |
| $n$ | Estimate of $S_{L_n}^{s_n}(x_n, \theta_n)$ | Exact Value of $S_{L_n}^{s_n}(x_n, \theta_n)$ | $|error|$ | Relative error |
| 12 | As exact value | $-5.99626094 - 6.21622453i$ | 0.0 | 0.0 |
| 11 | $-2.47725770 - 12.2276531i$ | $-2.46267551 - 12.2290382i$ | 0.0146 | $1.174 \times 10^{-3}$ |
| 10 | $11.8634890 + 20.0400796i$ | $11.8850320 + 20.02306601i$ | 0.0274 | $1.179 \times 10^{-3}$ |
| 9 | $48.22186610 + 20.4135977i$ | $48.2343029 + 20.3499866i$ | 0.0648 | $1.238 \times 10^{-3}$ |
| 8 | $80.8009685 - 38.3042272i$ | $80.721100 - 38.3815374i$ | 0.1111 | $1.244 \times 10^{-3}$ |
| 7 | $113.271534 + 92.1515455i$ | $113.109940 + 92.238325i$ | 0.1834 | $1.256 \times 10^{-3}$ |
| 6 | $-386.779631 + 164.145822i$ | $-386.400652 + 164.517221i$ | 0.5306 | $1.263 \times 10^{-3}$ |
| 5 | $876.083627 - 557.413123i$ | $876.314137 - 556.113049i$ | 1.3203 | $1.272 \times 10^{-3}$ |
| 4 | $1552.63044 - 998.33174$ | $1553.04478 - 996.036357i$ | 2.3324 | $1.264 \times 10^{-3}$ |
| 3 | $147.310219 + 3021.03492i$ | $150.742369 + 3019.38816i$ | 3.8067 | $1.259 \times 10^{-3}$ |
| 2 | $-1559.83456 - 4388.23283i$ | $-1553.97933 - 4388.90792i$ | 5.8642 | $1.259 \times 10^{-3}$ |
| 1 | $-5290.66224 + 11536.2495i$ | $-5301.47806 + 11524.4924i$ | 15.975 | $1.259 \times 10^{-3}$ |

Table 3.4B: As Table 3.4A, but for the quadratic Gaussian sum $S_{N_0}(1 - e/\pi, 1/e)$.

| $n$ | $L_n$ | $x_n$ | $s_n$ | $\theta_n$ |
|---|---|---|---|---|
| 1 | 129901233 | $\sqrt{2}/10$ | $+1$ | $\sqrt{10/71}$ |
| 2 | 18370808 | $7 - 10/\sqrt{2}$ | $-1$ | $-0.153724462171376$ |
| 3 | 1305572 | $x_2$ | $-1$ | $0.163067331555985$ |
| 4 | 92784 | $x_2$ | $-1$ | $-0.294531480224220$ |
| 5 | 6593 | $x_2$ | $-1$ | $0.144372430964127$ |
| 6 | 468 | $x_2$ | $-1$ | $-0.031474266260106$ |
| 7 | 33 | $x_2$ | $-1$ | $-0.442876534874666$ |
| $n$ | Estimate of $S_{L_n}^{s_n}(x_n, \theta_n)$ | Exact Value of $S_{L_n}^{s_n}(x_n, \theta_n)$ | $|error|$ | Relative error |
| 7 | As exact value | $6.36379881 - 0.158944615i$ | $0.0$ | $0.0$ |
| 6 | $19.1558359 - 13.4017658i$ | $19.1558145 - 13.4017626i$ | $2.163 \times 10^{-5}$ | $9.214 \times 10^{-7}$ |
| 5 | $64.1769995 + 59.1391565i$ | $64.176950 + 59.139096i$ | $7.732 \times 10^{-5}$ | $8.859 \times 10^{-7}$ |
| 4 | $-222.241162 + 241.875551i$ | $-222.241036 + 241.875239i$ | $3.373 \times 10^{-4}$ | $1.027 \times 10^{-6}$ |
| 3 | $-428.153221 - 1155.83662i$ | $-428.153209 - 1155.83534i$ | $1.280 \times 10^{-3}$ | $1.039 \times 10^{-6}$ |
| 2 | $-2663.73406 + 3781.54500i$ | $-2663.73285 + 3781.54032i$ | $4.484 \times 10^{-3}$ | $1.046 \times 10^{-6}$ |
| 1 | $12144.44498 - 1943.67131i$ | $12144.43440 - 1943.66515i$ | $1.224 \times 10^{-2}$ | $9.958 \times 10^{-7}$ |

Table 3.4C: As Table 3.4A, but for the quadratic Gaussian sum $S_{N_0}(\sqrt{2}/10, \sqrt{10/71})$.

| $n$ | $L_n$ | $x_n$ | $s_n$ | $\theta_n$ |
|---|---|---|---|---|
| 1 | 129901233 | $0.332613390928$ | $+1$ | $1/2e$ |
| 2 | 43206889 | $-0.00649350649$ | $-1$ | $-0.05301357552$ |
| 3 | 280564 | $3.628491 \times 10^{-12}$ | $+1$ | $-0.16409063118$ |
| 4 | $-1$ | irrelevant | irrelevant | irrelevant |
| $n$ | Estimate of $S_{L_n}^{s_n}(x_n, \theta_n)$ | Exact Value of $S_{L_n}^{s_n}(x_n, \theta_n)$ | $|error|$ | Relative error |
| 4 | As exact value | $0.5 + 0.0i$ | $0.0$ | $0.0$ |
| 3 | $-0.865720322 - 0.711435501i$ | $-0.865720322 - 0.711435501i$ | $< 10^{-10}$ | † $< 10^{-10}$ |
| 2 | $-4.920381960 + 3.937306408i$ | $-4.920381970 + 3.937306410i$ | $9.728 \times 10^{-8}$ | $1.542 \times 10^{-9}$ |
| 1 | $-10.05737802 + 2.726908802i$ | $-10.05611070 + 2.724765960i$ | $2.489 \times 10^{-3}$ | $2.389 \times 10^{-4}$ |

Table 3.4D: As Table 3.4A, but for the quadratic Gaussian sum $S_{N_0}(0.3326339\ldots, 1/2e)$.

| $n$ | $L_n$ | $x_n$ | $s_n$ | $\theta_n$ |
|---|---|---|---|---|
| 1 | 129901233 | $1/2 - \sqrt{\pi}/N_0^2$ | $+1$ | $1/\pi N_0$ |
| 2 | 64950616 | $4.201538 \times 10^{-16}$ | $-1$ | $-4.900798 \times 10^{-9}$ |
| 3 | 0 | irrelevant | irrelevant | irrelevant |
| $n$ | Estimate of $S_{L_n}^{s_n}(x_n, \theta_n)$ | Exact Value of $S_{L_n}^{s_n}(x_n, \theta_n)$ | $|error|$ | Relative error |
| 3 | As exact value | $0.5 + 0.0i$ | $0.0$ | $0.0$ |
| 2 | $(2.95151374 - 1.90568655i) \times 10^7$ | $(2.95151374 - 1.90568655i) \times 10^7$ | $4.432 \times 10^{-3}$ | $1.261 \times 10^{-10}$ |
| 1 | $(4.85720022 + 1.04582716i) \times 10^7$ | $(4.85720022 + 1.04582716i) \times 10^7$ | $7.052 \times 10^{-2}$ | $1.419 \times 10^{-9}$ |

Table 3.4E: As Table 3.4A, but for the quadratic Gaussian sum $S_{N_0}(1/2 - \sqrt{\pi}/N_0^2, 1/\pi N_0)$.

|  | $N = N_0$ | $x_0$ | $s_0$ | $\theta_0$ |
|---|---|---|---|---|
| 3.4A | 129901233 | $1/\sqrt{45}$ | $+1$ | $1 - \sqrt{(23/71)}$ |
| 3.4B | 129901233 | $1 - e/\pi$ | $+1$ | $1/e$ |
| 3.4C | 129901233 | $\sqrt{2}/10$ | $+1$ | $\sqrt{10/71}$ |
| 3.4D | 129901233 | $0.332613390928$ | $+1$ | $1/2e$ |
| 3.4E | 129901233 | $1/2 - \sqrt{\pi}/N_0^2$ | $+1$ | $1/\pi N_0$ |
|  | Estimate of $S_{N_0}^{s_0}(x_0, \theta_0)$ | Exact Value of $S_{L_n}^{s_n}(x_n, \theta_n)$ | $|error|$ | Relative error |
| 3.4A | $-4529.38172 - 4579.59925i$ | $-4527.85134 - 4577.13867i$ | $2.897$ $(18.946)$ | $4.501 \times 10^{-4}$ $(2.943 \times 10^{-3})$ |
| 3.4B | $-5299.74860 + 11526.7988i$ | $-5301.47806 + 11524.4924i$ | $2.883$ $(15.975)$ | $2.273 \times 10^{-4}$ $(1.259 \times 10^{-3})$ |
| 3.4C | $12144.43639 - 1943.66622i$ | $12144.43440 - 1943.66515i$ | $2.264 \times 10^{-3}$ $(1.224 \times 10^{-2})$ | $1.841 \times 10^{-7}$ $(9.958 \times 10^{-7})$ |
| 3.4D | $-10.0563258 + 2.72513736i$ | $-10.05611070 + 2.724765960i$ | $4.292 \times 10^{-4}$ $(2.489 \times 10^{-3})$ | $4.120 \times 10^{-5}$ $(2.389 \times 10^{-4})$ |
| 3.4E | $(4.85720022 + 1.04582716i) \times 10^7$ | $(4.85720022 + 1.04582716i) \times 10^7$ | $1.283 \times 10^{-2}$ $(7.052 \times 10^{-2})$ | $2.581 \times 10^{-10}$ $(1.419 \times 10^{-9})$ |

Table 4.5.2: Higher precision Gaussian sum computations for the examples shown in Tables 3.4A-E, incorporating the improvements (146) to the remainder terms. The resulting *refined algorithm QGS* produces reduced errors compared to those shown in parentheses found previously.

| | Computation of $Z(10^{18})$, for $[a] = 1595769123$, $N_t = 398942280$, $K = 30$, $P = 3$ & $\varepsilon_t = \varepsilon_0$. | | | | | |
|---|---|---|---|---|---|---|
| $p$ | Block & collection sizes | Final $\alpha$ value of each block. | Partial Sum | Corresponding final $N$ value. | Partial Sum | Relative Error |
| 0 | 63244 & 0 | 1595832367 | −0.033078 | 395406152 | −0.033033 | $1.37 \times 10^{-3}$ |
| 1 | 35218 & 9 | 1596536727 | 0.223712 | 386758725 | 0.226116 | $1.06 \times 10^{-2}$ |
| 2 | 39223 & 21 | 1598262539 | 0.226875 | 377255248 | 0.227283 | $1.80 \times 10^{-3}$ |
| 3 | 43684 & 33 | 1601233051 | −0.233534 | 367266439 | −0.233831 | $1.27 \times 10^{-3}$ |
| 4 | 48652 & 43 | 1605514427 | −0.028613 | 357221582 | −0.029194 | $1.99 \times 10^{-2}$ |
| 5 | 54185 & 53 | 1611366407 | 0.366864 | 346927266 | 0.366501 | $9.90 \times 10^{-4}$ |
| 6† | 60347 & 61 | 1618849435 | 0.085892 | 336615735 | 0.085828 | $7.46 \times 10^{-4}$ |
| 7 | 67211 & 69 | 1628258975 | −0.058909 | 326152678 | −0.058845 | $1.09 \times 10^{-3}$ |
| 8 | 74854 & 79 | 1640235615 | −0.132097 | 315225569 | −0.132370 | $2.06 \times 10^{-3}$ |
| 9 | 83367 & 87 | 1654908207 | 0.195195 | 304113503 | 0.195779 | $2.98 \times 10^{-3}$ |
| 10 | 92849 & 95 | 1672735215 | 0.075149 | 292793411 | 0.075219 | $9.31 \times 10^{-4}$ |
| 11 | 103408 & 105 | 1694657711 | −0.391719 | 281057961 | −0.391807 | $2.24 \times 10^{-4}$ |
| 12 | 115168 & 113 | 1720916015 | −0.120074 | 269163637 | −0.119802 | $2.27 \times 10^{-3}$ |
| 13 | 128266 & 123 | 1752725983 | 0.058025 | 256941103 | 0.058231 | $3.53 \times 10^{-3}$ |
| 14 | 142854 & 133 | 1791010855 | −0.138548 | 244460930 | −0.138875 | $2.35 \times 10^{-3}$ |
| 15* | 142854 & 261 | 1865866351 | −0.641492 | 224730918 | −0.642255 | $1.19 \times 10^{-3}$ |
| 16 | 142854 & 277 | 1945293175 | −0.116383 | 208194901 | −0.117073 | $5.89 \times 10^{-3}$ |
| 17 | 142854 & 293 | 2029291327 | −0.317633 | 193919717 | −0.319356 | $5.40 \times 10^{-3}$ |
| 18 | 142854 & 311 | 2118432223 | −0.127847 | 181283854 | −0.127376 | $3.70 \times 10^{-3}$ |
| 19 | 142854 & 329 | 2212715863 | −0.043415 | 169966199 | −0.042267 | $2.72 \times 10^{-2}$ |
| 20 | 142854 & 349 | 2312713663 | 0.444972 | 159686675 | 0.445641 | $1.50 \times 10^{-3}$ |
| 21 | 142854 & 369 | 2418425623 | 0.192445 | 150300383 | 0.192441 | $2.07 \times 10^{-5}$ |
| 22 | 142854 & 391 | 2530423159 | −0.098546 | 141652491 | −0.098288 | $2.62 \times 10^{-3}$ |
| 23 | 142854 & 413 | 2648706271 | −0.255487 | 133666527 | −0.255414 | $2.86 \times 10^{-4}$ |
| 24 | 142854 & 437 | 2773846375 | 0.244961 | 126245552 | 0.245158 | $8.04 \times 10^{-4}$ |
| 25 | 142854 & 463 | 2906414887 | 0.145447 | 119316290 | 0.144246 | $8.33 \times 10^{-3}$ |
| 26 | 142854 & 489 | 3046411807 | −0.005642 | 112847133 | −0.006166 | $8.50 \times 10^{-2}$ |
| 27 | 142854 & 515 | 3193837135 | 0.189155 | 106807432 | 0.189422 | $1.41 \times 10^{-3}$ |
| 28 | 142854 & 545 | 3349833703 | −0.059023 | 101128560 | −0.059113 | $1.52 \times 10^{-3}$ |
| 29 | 142854 & 575 | 3514401511 | 0.422062 | 95795370 | 0.421981 | $1.92 \times 10^{-4}$ |
| 30 | 142854 & 605 | 3687540559 | −0.240625 | 90791122 | −0.240133 | $2.05 \times 10^{-3}$ |
| 31 | 142854 & 639 | 3870393679 | −0.389069 | 86070339 | −0.389131 | $1.59 \times 10^{-4}$ |
| 32 | 142854 & 673 | 4062960871 | 0.262033 | 81623926 | 0.261562 | $1.80 \times 10^{-3}$ |
| 33 | 142854 & 711 | 4266384967 | 0.270332 | 77418503 | 0.270616 | $1.05 \times 10^{-3}$ |
| 34 | 142854 & 749 | 4480665967 | −0.524108 | 73448753 | −0.523929 | $3.42 \times 10^{-4}$ |
| 35 | 142854 & 789 | 4706375287 | −0.051565 | 69698128 | −0.051865 | $5.78 \times 10^{-3}$ |
| 36 | 142854 & 831 | 4944084343 | 0.367282 | 66152213 | 0.367593 | $8.46 \times 10^{-4}$ |
| 37 | 6712 & 875 | 4955843767: $\alpha_E^c$ | −0.139107 | 65986402: $N_C$ | −0.139114 | $5.03 \times 10^{-5}$ |
| | Total $ZP(10^{18})$ | | −0.376110 | $2 \sum_{N=N_c}^{N_t} \frac{\cos\{\theta_C(t) - t\log(N)\}}{\sqrt{N}} = -0.375570$ | | |
| $2 \sum_{N=1}^{N_c-1} \frac{\cos\{\theta_C(t) - t\log(N)\}}{\sqrt{N}} + O(t^{-1/4})$ | | | 0.565274 … | 0.565274 … | | |
| | Estimate for $Z(10^{18})$ | | 0.189164 | $Z(10^{18})$ | 0.189704 | $2.85 \times 10^{-3}$ |

Table 5.1: Detailed breakdown of calculations leading to the hybrid estimate of $Z(10^{18})$

| $Z(t)$ | $t = 10^{18}$ | $t = 10^{19}$ | $t = 10^{20}$ | $t = 10^{21}$ | $t = 10^{22}$ | $t = 10^{23}$ |
|---|---|---|---|---|---|---|
| RSF | 0.189704 | −28.270243 | 3.345199 | 2.610424 | −5.227095 | −1.608632 |
| (CPU time) | $(5.27 \times 10^3)$ | $(1.67 \times 10^4)$ | $(5.36 \times 10^4)$ | $(1.67 \times 10^5)$ | $(5.27 \times 10^5)$ | $(1.67 \times 10^6)$ |
| Hybrid Estimate | 0.189164 | −28.266343 | 3.339213 | 2.613327 | −5.221425 | −1.617695 |
| (CPU time) | $(5.93 \times 10^3)$ | $(1.65 \times 10^4)$ | $(4.46 \times 10^4)$ | $(1.16 \times 10^5)$ | $(2.97 \times 10^5)$ | $(7.50 \times 10^5)$ |
| Absolute Error | 0.000540 | −0.003900 | 0.005986 | −0.002903 | −0.005670 | −0.009065 |
| Relative Error | $2.85 \times 10^{-3}$ | $1.38 \times 10^{-4}$ | $1.79 \times 10^{-3}$ | $1.11 \times 10^{-3}$ | $1.08 \times 10^{-3}$ | $5.64 \times 10^{-3}$ |

Table 5.2: Hybrid estimates of $Z(10^{18-23})$ utilising *algorithm ZT13* with settings for the relative error $\varepsilon_t = \varepsilon_0$, termination constant $K = 30$ and step up parameter $Y = 10/9$. Comparisons of accuracy and sequential CPU time are made to the RSF.

| $\varepsilon_t$ | $K$ ; $Y$ | Hybrid Estimate of $Z(10^{20})$ | CPU time | Absolute Error | Relative Error |
|---|---|---|---|---|---|
| $\varepsilon_0/2$ | 30 ; 1.11 | 3.340945 | $4.90 \times 10^4$ | 0.004254 | $1.27 \times 10^{-3}$ |
| $\varepsilon_0$ | 30 ; 1.11 | 3.339213 | $4.46 \times 10^4$ | 0.005986 | $1.79 \times 10^{-3}$ |
| $2\varepsilon_0$ | 30 ; 1.11 | 3.334418 | $4.01 \times 10^4$ | 0.010781 | $3.22 \times 10^{-3}$ |
| $\varepsilon_0$ | 20 ; 1.11 | 3.339977 | $4.53 \times 10^4$ | 0.005221 | $1.53 \times 10^{-3}$ |
| $\varepsilon_0$ | 50 ; 1.11 | 3.339812 | $4.61 \times 10^4$ | 0.005387 | $1.61 \times 10^{-3}$ |
| $\varepsilon_0$ | 30 ; 1.08 | 3.346098 | $4.37 \times 10^4$ | −0.000090 | $2.68 \times 10^{-4}$ |
| $\varepsilon_0$ | 30 ; 1.14 | 3.346866 | $4.49 \times 10^4$ | −0.000166 | $4.98 \times 10^{-4}$ |

Table 5.3: Hybrid estimates of $Z(10^{20})$ utilising *algorithm ZT13* for different values of the relative error $\varepsilon_t$, termination constant $K$ and step up parameter $Y$.

**Figures**

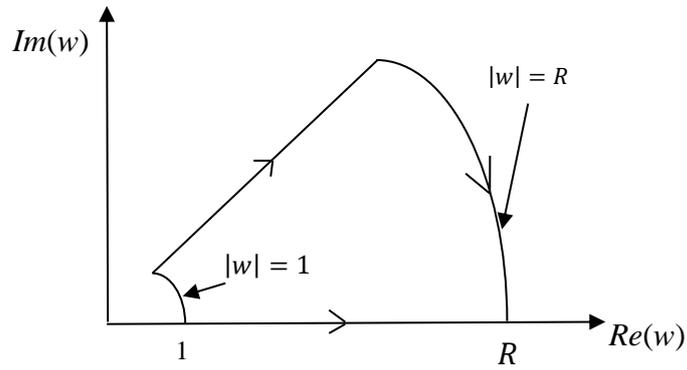

Fig. 1a

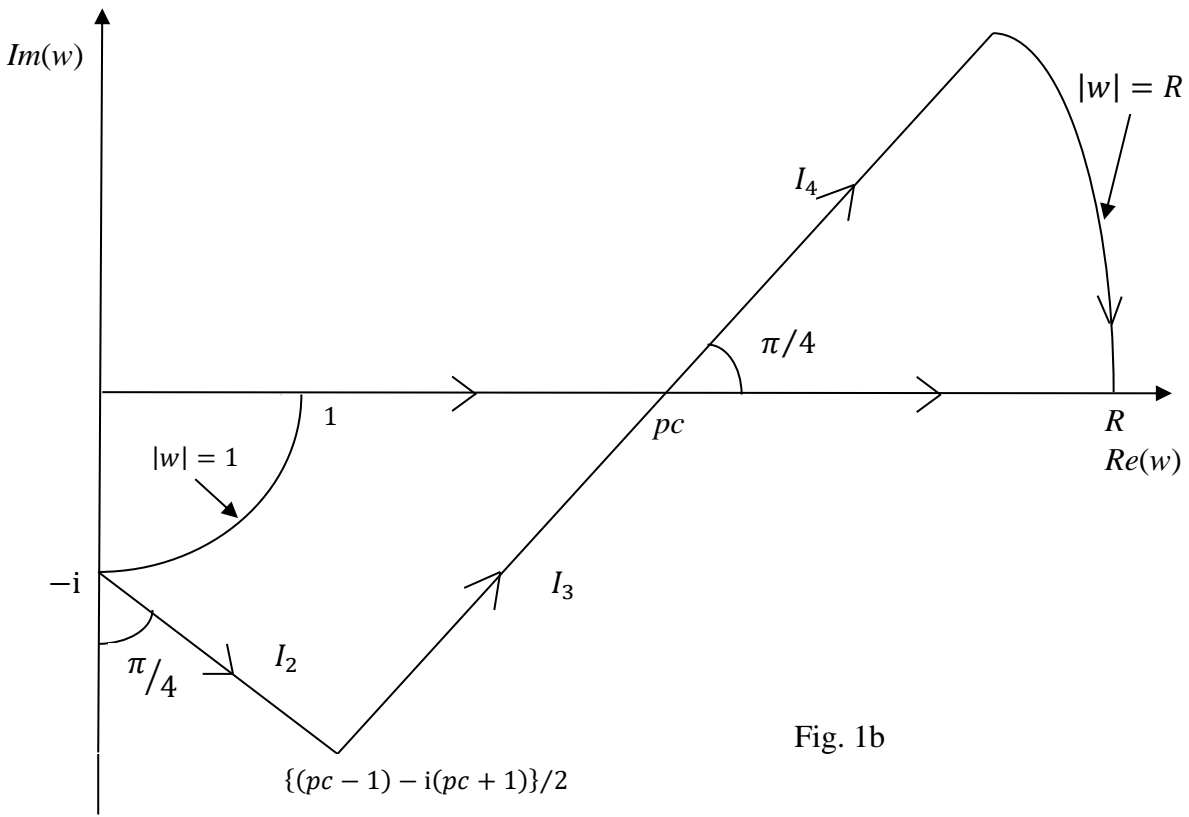

Fig. 1b

Figures 1a & b. Schematics of the general contours of integration in $w$ space used to estimate integral $B(\alpha, t)$ given by eq. (20). Fig. 1a is for the case when $\alpha < a$, whilst Fig. 1b is when $\alpha > a$. The integration through the saddle point situated at $w = pc$ in Fig. 1b, is what gives rise to the main terms in the formulation for $Z(t)$ specified by eq. (11).

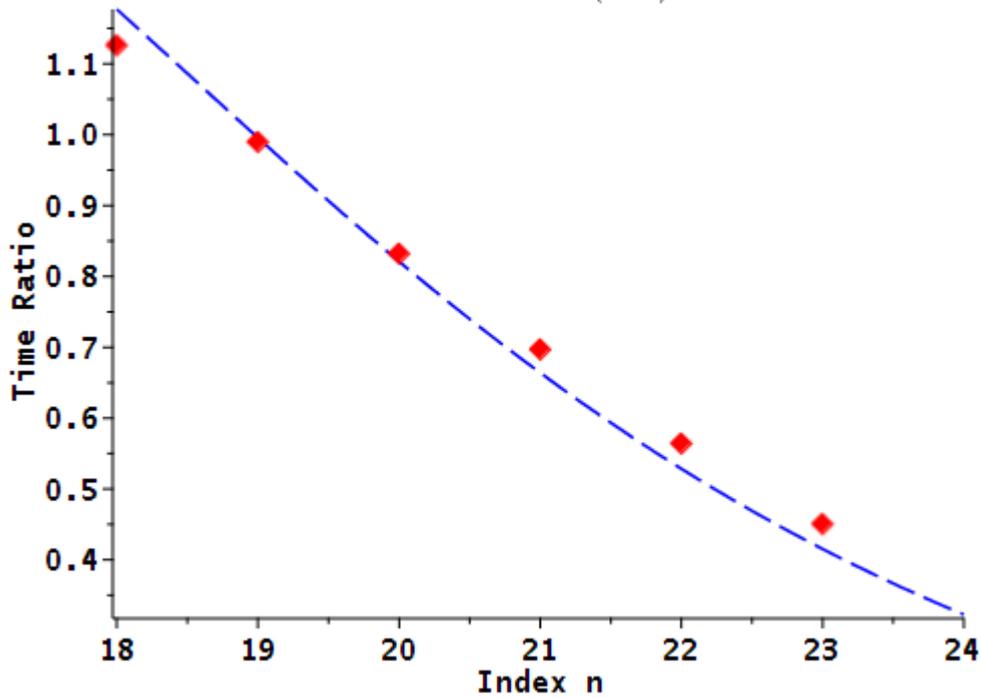

Fig.2a. Speed of hybrid formula to RS formula for calculation of $Z(10^n)$.

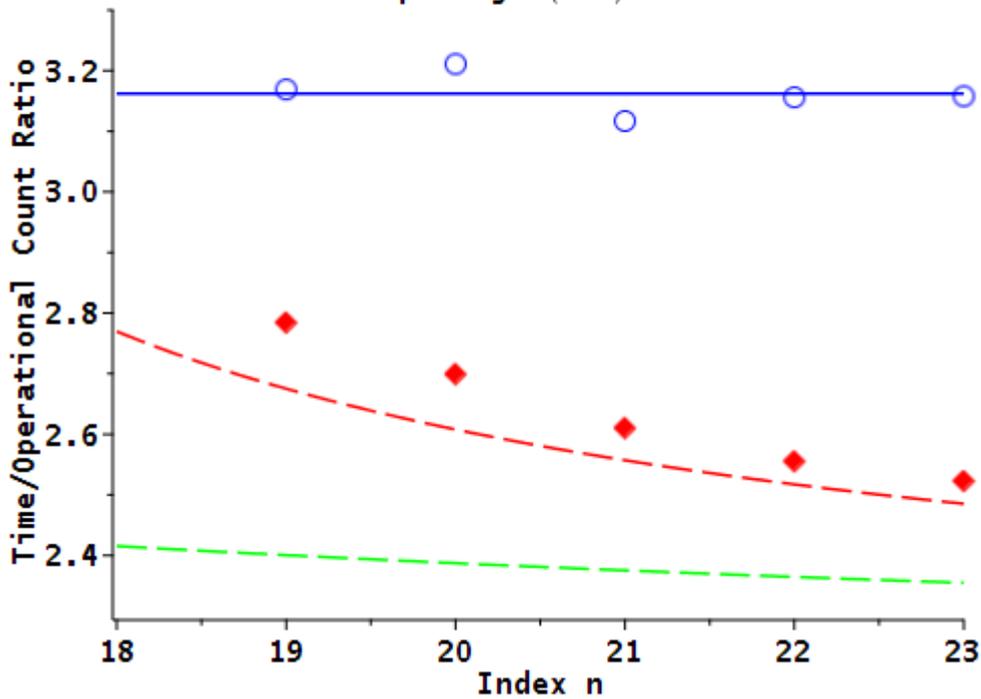

Fig.2b. Time/Operational count ratio for Computing $Z(10^n)$.

Figure 2a) Timings ratio of *algorithm ZT13* compared to the RSF. The red diamonds represent actual calculations; the blue dashed line is the fit obtained from the predicted upper bound operational count (153) with compensation factor $\lambda = 0.312$: & 2b) Timings ratio for the computation of $Z(10^n):Z(10^{n-1})$ using *algorithm ZT13*. The red diamonds represent actual calculations; the red dashed curve is the predicted behaviour of (153) across the range of $n = 19 - 23$; the green dashed curve is the asymptotic prediction as $n \to \infty$.

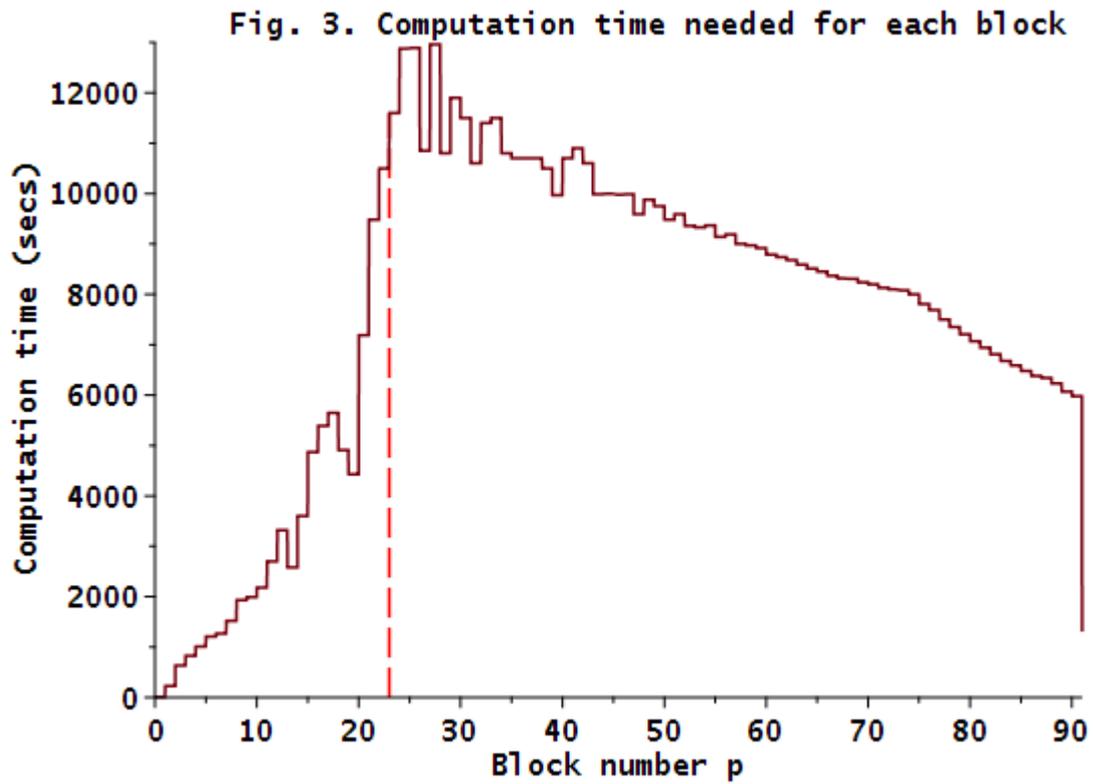

Figure 3. (Sequential) computational times necessary to complete the calculations for each pivot Block $p$. Data taken from the $Z(10^{23})$ computation, which is typical. The vertical dashed line represents Block no. 23, the point at which the step-up of the collection size $M_{t,p}$ is initiated in this case.